\documentclass[11pt, reqno]{amsart}
\usepackage[dvipsnames,svgnames,x11names]{xcolor}
\usepackage[tickmarkheight=.5em,textwidth=\marginparwidth,textsize=tiny]{todonotes}
\usepackage[mono=false]{libertine}
\usepackage[british]{babel}
\usepackage[margin=3cm]{geometry}
\usepackage[markup=underlined]{changes}
\definechangesauthor[color=teal]{A}
\definechangesauthor[color=BlueGreen]{B}
\usepackage{amsfonts}
\usepackage{amsthm,amsmath}
\usepackage{amsbsy}
\usepackage{bm}
\usepackage{amssymb} % Without this \qed does not get filled square
\usepackage{verbatim}
\usepackage{wrapfig}

\usepackage{times}
\usepackage{pgf,tikz}
\usepackage{graphicx}
\usepackage{tgpagella}
\usepackage{mathpazo}
\usetikzlibrary{trees,shadows}
\usetikzlibrary{arrows}
\usepackage{units}
\usepackage{enumerate,enumitem}
\usepackage{bbm}
\usepackage{exscale}
\usepackage[hyperindex,breaklinks]{hyperref}
\hypersetup{ colorlinks=true, linkcolor=blue, citecolor=blue,
  filecolor=blue, urlcolor=blue}
  \numberwithin{equation}{section} 

 \makeatletter
\def\paragraph{\@startsection{paragraph}{4}%
  \z@\z@{-\fontdimen2\font}%
  {\normalfont\itshape}}
  \makeatother
\DeclareMathOperator{\one}{\mathbbm{1}} %%%% indicator function

\usepackage[comma,sort&compress]{natbib}

\theoremstyle{plain} %documentation says there are only three styles
    \newtheorem{theorem}{Theorem}[section]
    \newtheorem{lemma}[theorem]{Lemma}
    \newtheorem{proposition}[theorem]{Proposition}
    
    \newtheorem{claim}[theorem]{Claim}

\theoremstyle{definition} % For roman text in the body
    \newtheorem{definition}[theorem]{Definition}

    \newtheorem{remark}[theorem]{Remark}

\DeclareMathOperator{\R}{\mathbb{R}}

\DeclareMathOperator{\Z}{\mathbb{Z}}

\DeclareMathOperator{\cI}{\mathcal{I}}
\usepackage{mathtools}
\DeclarePairedDelimiter\ceil{\lceil}{\rceil}
\DeclarePairedDelimiter\floor{\lfloor}{\rfloor}

\DeclareMathOperator{\De}{d}

\newcommand{\e}{\mathrm{e}}

\newcommand{\cF}{\mathcal{F}}
\newcommand{\cN}{\mathcal{N}}

\newcommand{\cS}{\mathcal{S}}
\newcommand{\cE}{\mathcal{E}}
\newcommand{\cQ}{\mathcal{Q}}
\newcommand{\cG}{\mathcal{G}}

\DeclareMathOperator{\diam}{\mathrm{diam}}

\DeclareMathOperator{\dist}{\mathrm{dist}}
\DeclareMathOperator{\supp}{\mathrm{supp}}

\newcommand{\g}{\gamma}
\renewcommand{\a}{\alpha}
\renewcommand{\b}{\beta}
\newcommand{\eps}{\varepsilon}
\newcommand{\vr}{\varphi}
\newcommand{\var}{\mathbf{Var}}
\newcommand{\cov}{\mathbf{Cov}}
\renewcommand{\r}{\rho}
\newcommand{\Lda}{\Lambda}

\newcommand{\prob}{\mathbb{P}}
\newcommand{\E}{\mathbb{E}}
\newcommand{\sgn}{\mathrm{sgn}}

\newcommand{\la}{\left\langle}
\newcommand{\ra}{\right\rangle}
\newcommand{\diver}{\operatorname{div}}

\date{\today}
\bibfont{\footnotesize}

\usepackage{etoolbox}
\newif\ifrevisionmarking
\revisionmarkingtrue
\newcommand{\revisioncolour}{\ifrevisionmarking\color{red}%
  \hypersetup{linkcolor=red,citecolor=red,urlcolor=red}\fi}

\newcommand{\HH}{\mathsf H}
\newcommand{\bx}{\mathsf Q}
\newcommand{\mom}{\mathfrak m}
\newcommand{\ip}[2]{\la #1,#2\ra}
\allowdisplaybreaks[2]
\begin{document}
\title[{Integer-Valued Membrane Model}]{Localisation and Delocalisation of the Integer-Valued Membrane Model}

%    author one information
\author[A. Cipriani]{Alessandra Cipriani}
\address{Department of Statistical Science, UCL, 1-19 Torrington Place, London, WC1E 7HB, UK }
\email{a.cipriani@ucl.ac.uk}
%%    author two  information
\author[B. Dan]{Biltu Dan}
\address{Department of Mathematics \& Statistics, Indian Institute of Science Education and Research Kolkata, Mohanpur 741246, West Bengal, India}
\email{bdan@iiserkol.ac.in}
%%    author three  information
\author[R.~S.~Hazra]{Rajat Subhra Hazra}
\address{Mathematical Institute, Leiden University, Einsteinweg 55, 2333
CC Leiden, The Netherlands}
\email{r.s.hazra@math.leidenuniv.nl}

\date{\today}

\begin{abstract}

We study the integer-valued membrane model on boxes
$[-L,L]^d\cap\Z^d$ with zero boundary conditions. We show that in dimension two the model
is delocalised at sufficiently high temperature. More precisely,
the variance at the origin is of order $L^2/\b$, where $\b$ is the
inverse temperature. More generally, for every $\eps>0$ and all
sufficiently small $\b>0$, the logarithm of each Laplace transform
lies between $(1+\eps)^{-1}$ and $1$ times its Gaussian counterpart.
The temperature threshold is independent of $L$ and the test
function. In dimensions three and
four, we prove that the variance at the origin is bounded uniformly
in the box size at every positive inverse temperature. Thus the
integer constraint changes the order of fluctuations in these
dimensions. The proofs use a modification of the Fr\"ohlich--Spencer
argument in dimension two and a comparison with integer-valued
Gaussian free fields in higher dimensions.

\end{abstract}

\maketitle
\tableofcontents

\section{Introduction}
In this article, we study the integer-valued membrane model (IVMM),
a random interface with a quadratic interaction involving the
discrete Laplacian. The model is defined on
$\Lda_L=[-L,L]^d\cap\Z^d$, with heights fixed to zero outside
$\Lda_L$. A configuration $m:\Z^d\to\Z$ has probability
proportional to $\exp[-\b H(m)]$, where $\b>0$ is the inverse
temperature and
\[
 H(m)=\frac12\sum_{x\in\Z^d}(\Delta m(x))^2
     =\frac12\la m,\Delta^2m\ra.
\]
Here $\Delta u(x)=\sum_{y:y\sim x}(u(y)-u(x))$ is the
unnormalised nearest-neighbour Laplacian, the symbol $\sim$ indicates graph nearest neighbours, and $\la\cdot,\cdot\ra$
denotes the usual inner product on finitely supported functions on the lattice. We are interested in how
the integer constraint affects the fluctuations of the interface
as $L$ tends to infinity.
Allowing the heights to take real values gives the membrane model
(MM), also known as the discrete bilaplacian model. This is a
centred Gaussian field whose covariance is $\b^{-1}$ times the
inverse of the bilaplacian matrix with zero boundary conditions. Its behaviour depends
on the dimension. At the centre of a box, for $L\ge2$, one has
\[
 \E^{MM}_{\b,\Lda_L}[\vr(0)^2]\asymp\frac1\b
 \begin{cases}
  L^2, & d=2,\\
  L, & d=3,\\
  \log L, & d=4.
 \end{cases}
\]
In dimensions $d\ge5$, the variance is bounded uniformly in $L$.
The constants implicit in $\asymp$ depend only on the dimension.
Thus four is the critical dimension for the real-valued model.
Entropic repulsion for this model was studied in
\citet{Sakagawa2003,Kurt2007,Kurt2009}.
Sharp estimates for the Green's function in dimensions two and
three were obtained by \citeauthor{muller2019estimates}.
The scaling limit of the real-valued model was identified by
\citeauthor{CiprianiDanHazra2019}: in dimensions two and three, it is a
H\"older-continuous Gaussian field, whereas in dimensions four and above
it is a random distribution.

The effect of restricting heights to integers is already
apparent for interfaces with a gradient interaction. As an example, the
two-dimensional discrete Gaussian free field has a variance which
grows logarithmically with the box size. Its integer-valued
counterpart is localised at low temperature and delocalised at
high temperature (that is, it passes from a ``flat'' phase to a ``rough'' phase, exhibiting a roughening transition). While the low-temperature result follows from a standard contour argument, the high-temperature result is due to
\citeauthor{F-S_1981}. Their complex proof uses an expansion into charges,
a multiscale decomposition and suitable changes of the field,
usually called spin waves. We refer to \citet{K-P_2017} for a
detailed presentation of this argument and to \citet{Wirth_2019}
for a treatment with zero boundary conditions and a study of the maximum of the field in that setting. 
It is therefore natural to ask whether the integer-valued
membrane model, in the same way, retains the fluctuation orders of the real-valued
model. The answer is different in dimension two and in dimensions
three and four. We use localisation to mean that the variance at
the origin is bounded uniformly in the box size. Delocalisation
means that this variance diverges as the box grows. Note that, since the
measure is invariant under $m\mapsto-m$, the field is centred, and therefore the variance offers the best observable to understand the behaviour of the field.
Our main results are as follows.
\begin{enumerate}[label=\textup{(\roman*)},leftmargin=*]
\item In dimension two, the model is delocalised for all
sufficiently small positive $\b$. More precisely, there exist
$\b_0>0$ and constants $0<c<C<\infty$, independent of $L$ and
$\b$, such that
\[
 \frac{c}{\b}L^2
 \le \E^{IMM}_{\b,\Lda_L}[m(0)^2]
 \le \frac{C}{\b}L^2,
 \qquad L\ge1,\quad 0<\b<\b_0
\]
(Theorem~\ref{thm:ivmm2d}). The upper bound holds for every
$\b>0$. The lower bound follows from a comparison of Laplace
transforms which is uniform in both the box size and the test
function.
\item In dimensions three and four, the model is localised at
every positive inverse temperature. For each $\b>0$ there is a
finite constant $c_\b$ such that
\[
 \sup_{L\ge1}\E^{IMM}_{\b,\Lda_L}[m(0)^2]\le c_\b
\]
(Theorem~\ref{thm:dim=3,4}). Thus the integer constraint surprisingly changes
the order of fluctuations in these dimensions.
\end{enumerate}
Before continuing, we complete the picture by recalling that, for every dimension and every $\b>0$, the Laplace transform
of the integer-valued model is bounded above by that of the
real-valued membrane model at the same inverse temperature
(Proposition~\ref{thm:upper_bound}). In particular, every linear
functional has variance at most that of its Gaussian counterpart. This fact follows from applying the Poisson summation formula to integer-valued Gaussians and is a feature of any discrete multivariate Gaussian.

Let us now recap our two-dimensional results. Theorem~\ref{thm:ivmm2d} and Proposition~\ref{thm:upper_bound} actually show something stronger than a variance bound. Writing $\Lda=\Lda_L$, let $\sigma_f$ solve
$\Delta^2\sigma_f=f$ on $\Lda$, with $\sigma_f=0$ outside
$\Lda$. For every $\eps\in(0,1)$, there exists $\b_\eps>0$,
independent of $L$, such that for $0<\b<\b_\eps$ and every
$f:\Z^2\to\R$,
\[
 \exp\!\left[\frac{1}{2(1+\eps)\b}\la\sigma_f,f\ra\right]
 \le \E^{IMM}_{\b,\Lda}\!\left[\e^{\la m,f\ra}\right]
 \le \exp\!\left[\frac{1}{2\b}\la\sigma_f,f\ra\right].
\]
The bounds show that the Gaussian comparison becomes arbitrarily
sharp as $\b$ decreases to zero, uniformly in $L$ and $f$. We note that the lower bound disproves a conjecture in~\cite{cuerno2001dynamic} stating that at high temperature the model would exhibit a GFF-type scaling limit.

We briefly explain the main ideas of the proofs. In dimension
two, we adapt the Fr\"ohlich--Spencer argument to the membrane
energy. The main difficulty is that the cancellation of the total charge (as defined in~\cite{K-P_2017}),
which is sufficient for a gradient energy, does not give the
required estimate in terms of second differences. Indeed,
controlling first differences by the membrane energy introduces
a factor depending on the size of the box. We therefore keep
track of both the total charge and its first moment. A charge
with both quantities equal to zero annihilates every affine
function. Subtracting a suitable affine function then allows us
to estimate its action directly in terms of second differences.
The spin waves must be modified accordingly. For this, we use a
local modification estimate of
\citet[Lemma~4.3]{schweiger2022pinning} and construct compatible
modifications for the separated charges. This gives a temperature
threshold which is independent of the box size.
The argument in dimensions three and four takes a different
form. We first express the membrane energy as a sum of gradient
energies of the slope components. A covariance comparison for
discrete Gaussians allows us to remove the constraints which
make these components the slopes of a single height function, thus making them independent.
This particular comparison follows from the lattice-Gaussian inequality of
\citet{RSD17-discrete-gaussian-ineq}. In dimension three, we then
use the mass-gap theorem of \citet{GM1982} for the integer-valued
Gaussian free field. The resulting summability of covariances,
together with a flow of finite energy, gives a uniform bound on
the height variance. In dimension four, a further covariance
comparison reduces the problem to the three-dimensional model.
%\paragraph{Open problems}Our results leave several natural questions open. 
%\begin{itemize}
%    \item Does the
%two-dimensional model localise at sufficiently low temperature?
% \item What is its scaling limit in the delocalised regime?
% \end{itemize}
%  The
% Laplace transform comparison gives uniform control of
% fluctuations, but identifying a limit at a fixed positive inverse
% temperature requires further arguments.
\paragraph{Structure of the paper.}
Section~\ref{sec:mainresults} gives the definitions and the precise
statements of the results. In
Section~\ref{sec:localisation} we prove localisation in dimensions
three and four. The two-dimensional lower bound is proved in
Section~\ref{sec:lower}. In Appendix~\ref{sec:proof} we prove the
Gaussian upper bound and deduce the two-dimensional variance
estimate from the Laplace transform lower bound. Appendix~\ref{app:discrete-gaussians}
contains the discrete-Gaussian comparison results used in the
localisation proof. The proofs of the covering bound and convex expansion (Proposition~\ref{prop:2.1} and Theorem~\ref{thm:2.2}) are given in Appendix~\ref{d2:sec:regrouping}.

\paragraph{Acknowledgements.} BD acknowledges support from the Indian Institute of Science Education and Research Kolkata through a start-up grant and from the Department of Science and Technology (DST), India, through the INSPIRE Faculty Fellowship, IFA22-MA176. Part of the work was carried out when BD was at the Indian Institute of Science. He thanks the institute for its generous hospitality. AC is supported by the EPSRC New Investigator Grant UKRI3308. The authors would like to thank Christophe Garban, Manjunath Krishnapur and Florian Schweiger for discussions and suggestions regarding the problem. Florian Schweiger also gave valuable comments on earlier results and a preliminary draft of this paper.

\paragraph{Use of AI} 
AI tools assisted the development of the proofs in dimensions
$d=2,3,4$. In particular, large language models (LLMs) suggested the comparison
with the three-dimensional integer-valued Gaussian free field used
in the localisation argument for $d=3,4$. For $d=2$, LLMs provided
an estimate missing from an earlier version of the proof and
suggested the use of affine functions to handle dipole charges. LLMs were also used for proofreading at the final stage of the project.
The authors take full responsibility for the mathematical arguments and
the presentation of the paper.
\section{Main results}\label{sec:mainresults}
Let $L\ge 1$ be an integer. Let $\Lda_L = [-L, L]^d \cap \Z^d$, a hypercube of side length $2L$ in the lattice $\Z^d$, for $d\ge 1$. For simplicity, we write $\Lda$ for $\Lda_L$ .

For any two vertices $x,y\in \Z^d$, we write $x\sim y$ if there is an edge between $x$ and $y$. We define the discrete Laplacian operator $\Delta$ as follows:
\begin{align*}
\Delta f(x):=\sum_{y\in\Z^d:y\sim x}(f(y)-f(x)), \qquad f:\Z^d\to \R.
\end{align*}
The discrete bilaplacian operator $\Delta^2$ is defined by $\Delta^2:=\Delta(\Delta)$.
%Also, we define $\Delta^2_{\Lda}=(\Delta^2(x,y))_{x,y\in\Lda}$ to be restriction of the matrix $\Delta^2$ onto $\Lda$. 
For any two functions $f, g:\Z^d\to\R$, at least one of which is finitely supported, we define the inner product
\begin{align*}
\la f, g\ra:=\sum_{x\in\Z^d} f(x) g(x).
\end{align*}
\begin{definition}[Integer-valued membrane model]
    We say that a random field $m:\Z^d\to\Z$ is an \textit{integer-valued membrane model} (IVMM) on $\Lda$ at inverse temperature $\b$ with zero boundary conditions if its law $\prob^{IMM}_{\b, \Lda}$ is given by
\begin{align}\label{eq:membrane-measure}
\prob^{IMM}_{\b, \Lda}\left(m=g\right) = \frac{1}{Z^{IMM}_{\b, \Lda}} \exp\left[-\frac{\b}{2} \la g, \Delta^2 g\ra\right] \prod_{x\in\Z^d\setminus\Lda}\one_{\{0\}}(g(x)),
\end{align}
where $\one_A$ is the indicator function of the set $A$ and $Z^{IMM}_{\b, \Lda}>0$ is the normalising constant.
\end{definition}
\begin{definition}[Membrane model] A random field $\vr:\Z^d\to\R$ is a (real-valued) \textit{membrane model} (MM) on $\Lda$ at inverse temperature $\b$ with zero boundary conditions if its law $\prob^{MM}_{\b, \Lda}$ is given by
\begin{align*}
\prob^{MM}_{\b, \Lda}(d\vr) = \frac{1}{Z^{MM}_{\b, \Lda}} \exp\left[-\frac{\b}{2} \la \vr, \Delta^2 \vr\ra\right] \prod_{x\in\Lda} d\vr(x)\prod_{x\in\Z^d\setminus\Lda}\delta_0(d\vr(x)),
\end{align*}
where $Z^{MM}_{\b, \Lda}$ is the normalising constant. 
\end{definition}
Also, we recall the definitions of the integer-valued discrete Gaussian free field (IVDGFF) and the (real-valued) discrete Gaussian free field. 
\begin{definition}[Integer-valued Gaussian free field]
    A random field $h:\Z^d \to\Z$ is an \textit{integer-valued discrete Gaussian free field} on $\Lda$ at inverse temperature $\b$ with zero boundary conditions if its law $\prob^{IGFF}_{\b, \Lda}$ is given by
\begin{align}\label{eq:def_IVGFF}
\prob^{IGFF}_{\b, \Lda}\left(h=g\right) = \frac{1}{Z^{IGFF}_{\b, \Lda}} \exp\left[-\frac{\b}{2} \la g,\, -\Delta g \ra\right] \prod_{x\in\Z^d\setminus\Lda}\one_{\{0\}}(g(x)).
\end{align}
Here $Z^{IGFF}_{\b, \Lda}$ is the normalising constant. 
\end{definition}
\begin{definition}[Discrete Gaussian free field]
    A random field $\vr:\Z^d\to\R$ is a (real-valued) \textit{discrete Gaussian free field} on $\Lda$ at inverse temperature $\b$ with zero boundary conditions if its law $\prob^{GFF}_{\b, \Lda}$ is given by
\begin{align*}
\prob^{GFF}_{\b, \Lda}(d\vr) = \frac{1}{Z^{GFF}_{\b, \Lda}} \exp\left[-\frac{\b}{2} \la \vr,\, -\Delta \vr \ra\right] \prod_{x\in\Lda} d\vr(x)\prod_{x\in\Z^d\setminus\Lda}\delta_0(d\vr(x)),
\end{align*}
where $Z^{GFF}_{\b, \Lda}$ is the normalising constant. 
\end{definition}

Throughout, functions supported in $\Lda$ or any other finite subset of $\Z^d$ are extended by zero to the full lattice. We write
\[
 \Delta_\Lda^2=(\Delta^2(x,y))_{x,y\in\Lda},
 \qquad G_\Lda=(\Delta_\Lda^2)^{-1}.
\]
The first matrix is the restriction of the full-lattice bilaplacian;
it is not the square of the restricted Laplacian.
For finitely supported $u$,
$\la u,\Delta^2u\ra=\|\Delta u\|_2^2$, with the norm taken on
the full lattice. This quadratic form is positive definite on
functions supported in a fixed finite set.
For a finite $D\subset\Z^d$ we also use
$\partial_2D=\{x\notin D:\dist(x,D)\le2\}$, where this
distance $\dist$ is the nearest-neighbour graph distance.

We now state precisely the main results of this paper.
%\section{Main results}
For $f:\Lda\to\R$, let $\sigma_f:\Z^d\to\R$ be the solution of 
\begin{equation}\label{eq:bvp_sigma}
\begin{cases}
    \Delta^2\sigma=f & \mbox{on } \Lda\\
    \sigma=0 &\mbox{on } \Z^d\setminus\Lda.
\end{cases}
\end{equation}
The following result, which is dimension-independent, is analogous to Proposition~$1.2$ of \cite{K-P_2017}. 
%It follows from \citet[Proposition 2.3]{garban2024invisibilityintegersdiscretegaussian}.
\begin{proposition}\label{thm:upper_bound}
For any $\b>0$ and $f:\Lda\to\R$,
\begin{align*}
\E^{IMM}_{\b, \Lda}\left[ \e^{\la m, f\ra }\right] \le \E^{MM}_{\b, \Lda} \left[ \e^{\la\varphi, f\ra}\right]= \exp\left[ \frac{1}{2\b}\la \sigma_f, f\ra\right],
\end{align*}
where $\sigma_f$ is the solution of \eqref{eq:bvp_sigma} corresponding to $f$.
\end{proposition}
The proof is a consequence of
\citet[Proposition~2.3]{garban2024invisibilityintegersdiscretegaussian}. For completeness, we provide a proof in Subsection~\ref{subsec:proof_Garban}. As a corollary of this, we know that the IVMM is localised whenever the MM is, that is, in dimensions greater than or equal to five.

\paragraph{Localisation in higher dimensions.}For $d=3,4$, we prove that the variance at the origin is bounded for any $\b>0$.
\begin{theorem}\label{thm:dim=3,4}
Let $d\in\{3,4\}$. For every $\b>0$ there exists a finite constant
$c_\b$ such that for all $L\geq 1$ one has that
 \begin{align*}
\E^{IMM}_{\b,\Lda }\left[ m(0)^2\right]  \le c_\b.
\end{align*}
\end{theorem}
This is quite striking, since, in the same dimensions, the MM is not localised. The crucial step in the proof where the integer-valued constraint enters is the comparison with the Gaussian free field: in its integer-valued version, the GFF is localised with finite susceptibility, which is not the case for the real-valued field. This is why the argument we are giving cannot work for the MM.

\paragraph{Delocalisation at high temperature in $d=2$.} For the lower bound corresponding to Proposition~\ref{thm:upper_bound}, we have the following result, analogous to Theorem~$1.1$ of \cite{K-P_2017} and Proposition~$8$ of \cite{Wirth_2019}, which prove a lower bound for the IVDGFF.

\begin{theorem}\label{thm:lower_bound_res}
Let $d=2$. For every $\eps\in(0,1)$ there exists $\b_\eps>0$ such that, for every $0<\b<\b_\eps$,
$L\ge1$, and $f:\Lda\to\R$,
\begin{equation}\label{eq:lower-uniform}
 \E^{IMM}_{\b,\Lda}\!\left[\e^{\la m,f\ra}\right]
 \ge
 \exp\!\left[\frac{1}{2(1+\eps)\b}\la\sigma_f,f\ra\right].
\end{equation}
Consequently,
\begin{equation}\label{eq:variance-uniform}
 \E^{IMM}_{\b,\Lda}\!\left[\la m,f\ra^2\right]
 \ge \frac{1}{(1+\eps)\b}\la\sigma_f,f\ra.
\end{equation}
\end{theorem}
As a consequence of Proposition~\ref{thm:upper_bound} and Theorem~\ref{thm:lower_bound_res}, we conclude the following.
\begin{theorem}\label{thm:ivmm2d}
Let $d=2$. There exist $\b_0>0$ and constants $0<c<C<\infty$ such that for every $0<\b<\b_0$ and $L\ge1$,
\[
 \frac{c}{\b}L^2
 \le \E^{IMM}_{\b,\Lda}[m(0)^2]
 \le \frac{C}{\b}L^2.
\]
The upper bound holds for every $\b>0$.
\end{theorem}
\begin{remark}
In the above theorem, we obtain a lower bound of the optimal order. In any case, one can obtain the following lower bound easily: for any $\eps>0$ there exists $\b_{1,\eps}>0$ such that for all $0<\b<\b_{1,\eps}$
\begin{align*}
\E^{IMM}_{\b, \Lda}\left[ m(0)^2\right] \ge \frac{c}{(1+\eps)\b} \log L.
\end{align*} 
Indeed, by a correlation inequality between two integer-valued Gaussians which is proved in \cite{RSD17-discrete-gaussian-ineq} (see also Theorem 2.4 in \cite{garban2024invisibilityintegersdiscretegaussian}) one can show that
\begin{align*}
\E^{IMM}_{\b, \Lda}\left[ m(0)^2\right] \ge \E^{IGFF}_{8\b, \Lda}\left[ m(0)^2\right].   
\end{align*}
We now only need to find the order of the right-hand side above. By Proposition~$8$ of \cite{Wirth_2019}, the discussion in Section~$4.2.2$ of \cite{garban2024invisibilityintegersdiscretegaussian} and \citet[Eq. (8.30)]{friedli_velenik_2017} one can conclude that for any $\eps>0$ there exists $\b_{2,\eps}>0$ such that for all $0<\b<\b_{2,\eps}$
\begin{align*}
\E^{IGFF}_{\b, \Lda}\left[ m(0)^2\right] \ge \frac{c}{(1+\eps)\b} \log L.
\end{align*}
\end{remark}

\section{Proof of Theorem~\ref{thm:dim=3,4}}\label{sec:localisation}
 \subsection{The \texorpdfstring{$d=3$}{d=3} case}
We begin by showing a few properties of the IVDGFF as introduced in~\eqref{eq:def_IVGFF}.

The following lemma holds in any dimension, but we state and use it only in $d=3$.
\begin{lemma}[Domain monotonicity of the IVDGFF]\label{lem:domain-monotonicity}
Let $D\subset D'\subset\Z^3$ be finite. 
%\[
%  \cov_{\prob^{IGFF}_{\b, D}}(h)
%  \preceq
%  \cov_{\prob^{IGFF}_{\b, D'}}(h).
%\]
For every function $f:\Z^3 \to \R$ with $\supp(f)\subset D$,
\begin{equation}\label{eq:domain-var}
  \var_{\prob^{IGFF}_{\b, D}}\left(\la f,h\ra\right)
  \le
  \var_{\prob^{IGFF}_{\b, D'}}\left(\la f, h\ra \right).
\end{equation}
\end{lemma}
\begin{proof}
Note that $\Z^{D'} \subseteq \R^{D'}$ is a lattice. One can view $\mathbb{Z}^{D}$ as the sublattice of $\Z^{D'}$ consisting of configurations that vanish on $D'\setminus D$.

Consider the vector space
\begin{align*}
    V= \left\{ f:D'\cup\partial_2D' \to \R \;\middle|\; f\equiv 0\text{ on }\partial_2D'
\right\} \subseteq
\R^{D'\cup\partial_2D'}.
\end{align*}
Define
\begin{align*}
   Q(v)=\la v,-\Delta v\ra, \qquad v\in V. 
\end{align*}
Then $Q$ is positive definite. Let 
\begin{align*}
    M=\left\{f:D'\cup\partial_2D'\to\Z \;\middle|\;f\equiv 0\text{ on } ( D'\cup\partial_2D')\setminus D \right\},
\end{align*}
and
\begin{align*}
    L=\left\{f:D'\cup\partial_2D'\to\Z \;\middle|\;f\equiv 0\text{ on } \partial_2D' \right\}.
\end{align*}
Then $M\subset L\subset V$ are discrete additive subgroups of $V$. In fact, $M\subset L \subset \R^{D'\cup\partial_2D'}$ are lattices. Hence \eqref{eq:domain-var} follows from Proposition~ \ref{prop:sublattice-cov} with the above $Q$.
\end{proof}

{We use the following classical theorem in its height-function form. The exponential decay of correlations for the three-dimensional
integer-valued Gaussian free field at every positive inverse
temperature is due to~\citet{GM1982}.
See also the discussions in
\citet[Section~1.2]{LammersOtt2024} and
\citet[Section~1.4.2]{adhikari2026rigidity}. The zero-boundary height two-point statement below is also recorded explicitly in \citet[p.~1]{RFN2024}.
\begin{theorem}[{Mass gap; \cite{GM1982}}]\label{thm:mass-gap}
For every $\b>0$, let $B_n=[-n,n]^3\cap\Z^3$ and let
$\prob^{IGFF}_{\b,B_n}$ be the zero-boundary measure in \eqref{eq:def_IVGFF}.
For each fixed $z\in\Z^3$, the limit
\[
  C_\b(z)
  :=\lim_{n\to\infty}
       \E_{\prob^{IGFF}_{\b,B_n}}[h(0)h(z)]
\]
exists.  Moreover, there are constants $A_\b,a_\b>0$ such that
\begin{equation}\label{eq:mass-gap}
  0\le C_\b(z)\le A_\b \e^{-a_\b|z|} ,
  \qquad z\in\Z^3.
\end{equation}
\end{theorem}
The consequence of this result that we will use is the summability of the covariance
kernel:
\[
    \chi_\beta
    :=\sum_{z\in\mathbb{Z}^3}C_\beta(z)<\infty.
\]
}

The following result is a consequence of the FKG property of the IVDGFF. For completeness, we give a proof in Appendix~\ref{sec:appendix_B.2}.

\begin{lemma}[Positive association of the IVDGFF]\label{lem:FKG}
Let $D\subset\Z^3$ be finite. Then, for every $x,y\in D$,
\[
  \cov_{\prob^{IGFF}_{\b,D}}(h(x),h(y))\ge0.
\]
\end{lemma}
\begin{remark}
Lemma~\ref{lem:FKG} gives an alternative proof that $C_\b(z)\geq 0.$
\end{remark}
We now show a crucial bound of the variance in terms of the energy of a function and susceptibility.
\begin{proposition} \label{prop:susceptibility}
For every finite $D\subset\Z^3$ and every function $f:\Z^3 \to \R$ with $\supp(f)\subset D$,
\begin{equation}\label{eq:susceptibility-bound}
  \var_{\prob^{IGFF}_{\b,D}}\left(\la f, h\ra\right)
  \le \chi_\b \la f, f\ra,
\end{equation}
where
\begin{equation}\label{eq:chi-def}
  \chi_\b=\sum_{z\in\Z^3}C_\b(z)<\infty.
\end{equation}
\end{proposition}
\begin{proof}
Let $B_n=[-n,n]^3\cap\Z^3$ for $n\ge 1$. For a finitely supported $f:\Z^3\to\R$, define
\begin{align*}
   V_n(f):= \var_{\prob^{IGFF}_{\b,B_n}}\left(\la f, h\ra\right),  
\end{align*}
for all sufficiently large $n$ such that $\supp (f) \subset B_n$.  By Lemma~\ref{lem:domain-monotonicity}, $V_n(f)$ is increasing in $n$. We claim that it is bounded. 
%To prove the claim, we first show that
%\begin{align*}
% \sup_n\var_{\prob^{IGFF}_{\b,B_n}}(h_x)\le C_\b(0).   
%\end{align*}

Fix $x\in\supp(f)$. Let $r=\|x\|_\infty$. Then $B_n\subset x+B_{n+r}.$
Using Lemma~\ref{lem:domain-monotonicity} and then translating the domain, we have
\begin{align*}
  \var_{\prob^{IGFF}_{\b,B_n}}(h(x))
  \le
  \var_{\prob^{IGFF}_{\b,x+B_{n+r}}}(h(x))
  =
  \var_{\prob^{IGFF}_{\b,B_{n+r}}}(h(0)).
\end{align*}
The right-hand side is increasing in $n$ and converges, by
Theorem~\ref{thm:mass-gap} with $z=0$, to $C_\b(0)<\infty$.  Hence
\begin{align*}
  \sup_n\var_{\prob^{IGFF}_{\b,B_n}}(h(x))\le C_\b(0).
\end{align*}
Now, using the Cauchy--Schwarz inequality, we have
\begin{align*}
  V_n(f)= \sum_{x,y}f(x)f(y)\E_{\prob^{IGFF}_{\b,B_n}}\left(h(x) h(y)\right)
  \le C_\b(0)\left(\sum_x|f(x)|\right)^2.
\end{align*}
Thus we have proved that $V_n(f)$ is increasing and bounded. Therefore 
\begin{align*}
  V_\b(f):=\lim_{n\to\infty}V_n(f)
\end{align*}
exists and is finite.

We next prove that $V_\b$ is translation invariant, i.e., for any $a\in\Z^3$
\begin{equation}\label{eq:translation-V}
  V_\b(f)=V_\b(\tau_af),
\end{equation}
 where $(\tau_af)(y):=f(y+a)$ for $y\in\Z^3$.
Let $a\in\Z^3$
and $r=\|a\|_\infty$.  Since $B_n\subset a+B_{n+r}$, by Lemma~\ref{lem:domain-monotonicity}
\begin{align*}
  V_n(f)
\le\var_{\prob^{IGFF}_{\b,a+B_{n+r}}}\left(\la f, h\ra\right).
\end{align*}
Now, translating the whole configuration by $-a$, we get 
\begin{align*}
\var_{\prob^{IGFF}_{\b,a+B_{n+r}}}\left(\la f, h\ra\right)= V_{n+r}(\tau_af),  
\end{align*}

Thus $V_n(f) \le V_{n+r}(\tau_af)$ and hence $V_\b(f)\le V_\b(\tau_af).$
Applying the same argument to $-a$ yields the reverse inequality. Hence \eqref{eq:translation-V} holds.

Next, define the symmetric bilinear form
\begin{align*}
  \mathcal C_\b(f,g)
  :=\frac12\bigl(V_\b(f+g)-V_\b(f)-V_\b(g)\bigr),
\end{align*}
on functions with finite support.
Since $V_\b(f)$ is a quadratic form, we have
\begin{align*}
\mathcal C_\b(f,f)=V_\b(f).
\end{align*}
By \eqref{eq:translation-V}, $\mathcal C_\b(f,\,g)=\mathcal C_\b(\tau_a f,\,\tau_a g)$; therefore, the kernel 
  $K_\b(x,y):=
  \mathcal C_\b(\one_{\{x\}},\one_{\{y\}})$
depends only on $y-x$; namely, 
\begin{equation}\label{eq:three_here}
K_\b(x,y)=K_\b(0,\,y-x).    
\end{equation} 
For a function $f$ with finite support, we have
\begin{equation}\label{eq:one_here}
V_\b(f)= \mathcal C_\b\left(\sum_{x}f(x)\one_{\{x\}},\, \sum_{y}f(y)\one_{\{y\}}\right)=\sum_{x,\,y}f(x) f(y) K_\b(x,
\,y).
\end{equation}
Also, we have
\begin{align}
   K_\b(0,\,z)&=\lim_{n\to\infty}\frac12\left[V_{n}(\delta_0+\delta_z)-V_n(\delta_0)-V_n(\delta_z)\right]\nonumber\\
   &=\lim_{n\to\infty}\mathrm{Cov}_{\prob^{IGFF}_{\b,B_n}}(h(0),\,h(z))=\lim_{n\to\infty}\mathbb{E}_{\prob^{IGFF}_{\b,B_n}}[h(0)h(z)]=C_\b(z).\label{eq:two_here}
\end{align}
Consequently, \eqref{eq:three_here}, \eqref{eq:one_here} and \eqref{eq:two_here} give
\begin{equation}\label{eq:V-kernel}
  V_\b(f)=\sum_{x,y}f(x)f(y)C_\b(y-x).
\end{equation}
By Lemma~\ref{lem:FKG}, we have that for any $z\in\Z^3$, $\mathbb{E}_{\prob^{IGFF}_{\b,B_n}}[h(0)h(z)] \ge 0$ for all large $n$ such that $z\in B_n$. Hence $C_\b(z)\ge0$.  By Theorem~\ref{thm:mass-gap},
\begin{align*}
  \chi_\b:=\sum_zC_\b(z)
  \le A_\b\sum_z\e^{-a_\b|z|}<\infty.
\end{align*}
Now, using $2|ab|\le a^2+b^2$ and non-negativity of $C_\b$,
\begin{align*}
  V_\b(f)
  &\le\sum_{x,y}|f(x)||f(y)|C_\b(y-x)\\
  &\le\frac12\sum_{x,y}
       (f(x)^2+f(y)^2)C_\b(y-x)\\
  &=\chi_\b\sum_xf(x)^2.
\end{align*}
Finally, if $D$ is a finite subset of $\Z^3$, then $D\subset B_n$ for all large $n$. Hence Lemma~\ref{lem:domain-monotonicity} gives
\begin{align*}
  \var_{\prob^{IGFF}_{\b,D}}\left(\la f,h\ra\right)
  \le V_n(f)\le V_\b(f),
\end{align*}
and \eqref{eq:susceptibility-bound} follows.

\end{proof}

\subsubsection{The membrane energy as a gradient energy of slopes}

For $1\le i\le d$, define the forward and backward differences
\begin{align*}
  &D_i^+f(x):=f(x+e_i)-f(x)\\
  &D_i^-f(x):=f(x)-f(x-e_i).
\end{align*}
So
\begin{equation}\label{eq:laplacian-factor}
  \Delta=\sum_{i=1}^dD_i^+D_i^-.
\end{equation}

\begin{lemma}[Discrete Hessian identity,{~\cite[Lemma 3.2]{bolthausen2016fastdecaycovariancesdeltapinning}}]\label{lem:hessian}
For every finitely supported real-valued function $f$ on $\Z^d$,
\begin{equation}\label{eq:hessian-id}
  \|\Delta f\|_2^2
  =\sum_{i,j=1}^d\|D_i^+D_j^+f\|_2^2.
\end{equation}
\end{lemma}

\begin{proof}
Because $\Delta$ is self-adjoint,
\begin{align*}
  \|\Delta f\|_2^2
  =\la f, (-\Delta)^2f\ra.
\end{align*}
On the other hand, using commutation of differences and
\eqref{eq:laplacian-factor},
\begin{align*}
  \sum_{i,j}\|D_i^+D_j^+f\|_2^2
  &=\sum_{i,j}
     \la f, D_j^-D_i^-D_i^+D_j^+f\ra\\
  &=\la f, 
      \left(\sum_iD_i^-D_i^+\right)
      \left(\sum_jD_j^-D_j^+\right)f\ra\\
  &=\la f, (-\Delta)^2f\ra.
\end{align*}
This proves \eqref{eq:hessian-id}.
\end{proof}

The idea now is that each slope component carries the energy of an ordinary scalar integer-valued GFF.  The only coupling among the slope components is the condition that they come from one height field.  That condition means that the slope configurations form a sublattice of the full product slope space, and Lemma~\ref{lem:slope-relaxation} shows that imposing it can only decrease variances of linear slope observables. We now make this precise.

For $j=1,2,3$, define the finite support set
\begin{equation}\label{eq:Sj}
  \Lda_{L,j}:=\{x\in\Z^3:\ x\in\Lda
             \text{ or }x+e_j\in\Lda\}.
\end{equation}
If $f:\Z^3\to\Z$ is supported in $\Lda$, then $D_j^+f$ is supported in $\Lda_{L,j}$.
Let
\begin{align*}
  P_L:=\bigoplus_{j=1}^3\Z^{\Lda_{L,j}}
\end{align*}
be the full product slope lattice.

\begin{lemma}\label{lem:slope-relaxation}
For every linear functional $f$ on $P_L$,
\begin{equation}\label{eq:slope-relax}
\var_{\prob^{IMM}_{\b,\Lda}}\left(f(D_1^+m,D_2^+m,D_3^+m)\right)\le \var_{\bigotimes_{j=1}^3\prob^{IGFF}_{\b,\Lda_{L,j}}}
       \left(f(h_1,h_2,h_3)\right).
\end{equation}
\end{lemma}
\begin{proof}
      Define on $P_L$ the positive-definite
quadratic form
\begin{equation}\label{eq:product-slope-Q}
  Q(u):=
  \sum_{j=1}^3\sum_{x\in\Z^3}\sum_{i=1}^3
       \bigl(D_i^+u_j(x)\bigr)^2.
\end{equation}
It is positive definite because a finitely supported function with zero
gradient is zero.

Define the slope map $T:\Z^{\Lda}\longrightarrow P_L$ by
\begin{align*}
  Tf=(D_1^+f,D_2^+f,D_3^+f).
\end{align*}
Then $T$ is injective, as a finitely supported function with zero
gradient is zero.  Define
\begin{align*}
  H_L:=T(\Z^{\Lda}).
\end{align*}
Then $H_L$ is a sublattice of $P_L$.  

For any function $g$ supported in $\Lda$, set
\begin{equation}\label{eq:slope-def}
  s_j:=D_j^+g,\qquad j=1,2,3.
\end{equation}
Then Lemma~\ref{lem:hessian} gives the exact decomposition
\begin{equation}\label{eq:slope-energy}
 \la g, \Delta^2 g \ra = \|\Delta g\|_2^2
  =\sum_{j=1}^3\sum_{i=1}^3\|D_i^+s_j\|_2^2.
\end{equation}

By \eqref{eq:slope-energy}, the pushforward of the integer-valued membrane law \eqref{eq:membrane-measure} under $T$ is precisely the centred discrete Gaussian on $H_L$ with weight $\exp[-\b Q/2]$. Also, let $Y=(Y_1,Y_2,Y_3)$, where $Y_j\in \Z^{\Lda_{L,j}}, \, j=1,2,3$, be a centred discrete Gaussian on $P_L$ such that
\begin{align*}
  \prob(Y=y)\propto e^{-\b Q(y)/2},\qquad y\in P_L.
\end{align*}  
Then, from the definition of $Q$, it follows that the $Y_j$ are independent and the law of $Y_j$ is $\prob^{IGFF}_{\b,\Lda_{L,j}}$.

Now, \eqref{eq:slope-relax} is an immediate consequence of the above observation together with Proposition~\ref{prop:sublattice-cov}.
\end{proof}

Theorem~\ref{thm:dim=3,4} for $d=3$ follows from the following proposition.
\begin{proposition}[Uniform variance bound in dimension $3$]\label{prop:d3}
For every $\b>0$ and every $L$,
\begin{equation}\label{eq:d3-bound}
  \var_{\prob^{IMM}_{\b,\Lda}}(m(0))
  \le \chi_\b R_3,
\end{equation}
where $\chi_\b$ is defined in \eqref{eq:chi-def} and
$R_3=G(0)/6$.
\end{proposition}

\begin{proof}
First we want to write $m(0)$ as a linear functional of the slopes $D_j^+m$ with a coefficient vector whose norm is bounded uniformly in $L$. 
%The required coefficients are a unit flow of finite energy.
Let $(X_n)_{n\ge0}$ be a simple random walk on $\Z^3$ and let
\begin{equation}\label{eq:Green}
  G(x):=\sum_{n=0}^{\infty}\prob_0(X_n=x)
\end{equation}
be its Green's function. We know that $G(0)<\infty$.
The Green's function satisfies
\begin{equation}\label{eq:Green-equation}
  -\Delta G=6\,\delta_0.
\end{equation}
Put
\begin{align*}
  V(x):=\frac16G(x),
\end{align*}
and define
\begin{equation*}\label{eq:global-current}
  J_j(x):=V(x)-V(x+e_j),\qquad j=1,2,3,
\end{equation*}
and
\begin{equation*}\label{eq:divergence}
  (\diver J)(x)
  :=\sum_{j=1}^3\bigl(J_j(x)-J_j(x-e_j)\bigr).
\end{equation*}
Then \eqref{eq:Green-equation} gives
\begin{equation*}\label{eq:current-div}
  \diver J=\delta_0.
\end{equation*}
On $[-\pi,\pi]^3$, put
\[
 a(t)=2\sum_{i=1}^3(1-\cos t_i).
\]
The Fourier representation of $V$ is $\widehat V(t)=a(t)^{-1}$.
The Fourier transform of $J_i$ is $(1-e^{it_i})/a(t)$, which is square
integrable because $a(t)\asymp|t|^2$ near zero. Plancherel therefore gives
\[
 \sum_{x\in\Z^3}\sum_{i=1}^3 J_i(x)^2
 =\frac1{(2\pi)^3}\int_{[-\pi,\pi]^3}
   \frac{\sum_i|1-\e^{it_i}|^2}{a(t)^2}\,dt
 =\frac1{(2\pi)^3}\int\frac{dt}{a(t)}=V(0)<\infty.
\]

For the finite box $\Lda$, retain only the edges having at least one
endpoint in $\Lda$:
\begin{equation*}\label{eq:truncated-flow}
  J^{(L)}_j(x)
  :=J_j(x)\,
  \one_{\{x\in\Lda\text{ or }x+e_j\in\Lda\}}.
\end{equation*}
Then
\begin{equation}\label{eq:finite-flow-div}
  \diver J^{(L)}(x)=\delta_0(x),
  \qquad x\in\Lda.
\end{equation}
Also
\begin{equation}\label{eq:finite-flow-energy}
  \sum_{x\in\Z^3}\sum_{j=1}^3\bigl(J^{(L)}_j(x)\bigr)^2
  \le R_3,
  \qquad
  R_3:=\frac{G(0)}6<\infty.
\end{equation}

Finally, because $m=0$ outside $\Lda$, discrete summation by parts
and \eqref{eq:finite-flow-div} give
\begin{align}
  \sum_{x\in\Z^3}\sum_{j=1}^3
    J^{(L)}_j(x)D_j^+m(x)
  &=-\sum_{x\in\Z^3}m(x)\diver J^{(L)}(x)\notag\\
  &=-m(0).\label{eq:m0-flow}
\end{align}
Thus $m(0)$ is a linear functional of the slopes $D_j^+m$ with a coefficient vector whose norm is bounded uniformly in $L$ by~\eqref{eq:finite-flow-energy}.

Consider the linear functional $f:\bigoplus_{j=1}^3\Z^{\Lda_{L,j}} \to \R$ defined by
\begin{align*}
    f(y_1,y_2,y_3)= -\sum_{j=1}^3\sum_{x\in\Z^3}J^{(L)}_j(x)y_j(x).
\end{align*}
Now, using \eqref{eq:m0-flow} and applying Lemma~\ref{lem:slope-relaxation}, we get
\begin{align*}
  \var_{\prob^{IMM}_{\b,\Lda}}(m(0))&=  \var_{\prob^{IMM}_{\b,\Lda}}\left( f(D_1^+m, D_2^+m, D_3^+m)\right)\\
  &\le\var_{\bigotimes_{j=1}^3\prob^{IGFF}_{\b,\Lda_{L,j}}}
       \left(f(h_1,h_2,h_3)\right)\\ 
  &=
  \sum_{j=1}^3
 \var_{\prob^{IGFF}_{\b,\Lda_{L,j}}}\left(
     \sum_{x\in\Z^3}J^{(L)}_j(x)h(x)\right).
\end{align*}
By Proposition~\ref{prop:susceptibility}, 
\[
 \var_{\prob^{IGFF}_{\b,\Lda_{L,j}}}\left(
     \sum_{x\in\Z^3}J^{(L)}_j(x)h(x)\right) \le \chi_\b\sum_{x\in\Z^3}\bigl(J^{(L)}_j(x)\bigr)^2.
\]
Summing over $j$ and using \eqref{eq:finite-flow-energy} gives
\[
  \var_{\prob^{IMM}_{\b,\Lda}}(m(0))
  \le
  \chi_\b\sum_{j=1}^3\sum_{x\in\Z^3}\bigl(J^{(L)}_j(x)\bigr)^2
  \le\chi_\b R_3.
\]
\end{proof}

\subsection{The \texorpdfstring{$d=4$}{d=4} case}
The four-dimensional model will be reduced to the three-dimensional one by removing non-negative square terms from the energy.  For continuous Gaussians, such a variance comparison is immediate from matrix order. For discrete Gaussians it is not automatic. Lemma~\ref{lem:add-squares} shows that it is nevertheless true.

We complete the proof in $d=4$. 

Let $g:\Z^4 \to \Z$ be a function supported in
$\Lda\subset\Z^4$.  By Lemma~\ref{lem:hessian},
\begin{equation}\label{eq:Q4}
  Q(g):=\|\Delta g\|_2^2
  =\sum_{i,j=1}^4\|D_i^+D_j^+g\|_2^2.
\end{equation}
We write
\begin{equation}\label{eq:Qsplit}
  Q(g)=Q_0(g)+Q_1(g),
\end{equation}
where
\begin{equation}\label{eq:Q0}
  Q_0(g):=\sum_{i,j=1}^3\|D_i^+D_j^+g\|_2^2,
\end{equation}
and

\begin{equation}\label{eq:Q1}
 Q_1(g):=Q(g)-Q_0(g)
 =\sum_{\substack{1\le i,j\le4\\i=4\ \text{or}\ j=4}}
       \|D_i^+D_j^+g\|_2^2.
\end{equation}

Define the integer-linear maps $A$ and $B$ on $\Z^{\Lda}$ by
\begin{align*}
    &A(g)= \left(D_i^+ D_j^+g\right)_{1\le i,j\le 3},\\
    &B(g)= \left(D_i^+ D_j^+g\right)_{1\le i,j\le 4,\, i=4 \text{ or }j=4}.
\end{align*}

We first show that $A$ is injective.  Suppose $Ag=0$.
For each fixed fourth coordinate $t$, let
\begin{align*}
  g^{(t)}(x_1,x_2,x_3):=g(x_1,x_2,x_3,t).
\end{align*}
Then
\begin{align*}
  \sum_{i,j=1}^3\|D_i^+D_j^+g^{(t)}\|_2^2=0.
\end{align*}
By Lemma~\ref{lem:hessian},
\begin{align*}
  \|\Delta g^{(t)}\|_2^2=0,
\end{align*}
so $\Delta g^{(t)}=0$.  Since $g^{(t)}$ has finite support,
\begin{align*}
  0=\la g^{(t)}, -\Delta g^{(t)}\ra
   =\sum_{x\in\Z^3}\sum_{i=1}^3(D_i^+g^{(t)}(x))^2,
\end{align*}
and therefore $g^{(t)}\equiv0$.  This holds for every $t$, so $g=0$.
Thus $A$ is injective.

Also, observe that 
\begin{align*}
    Q(g)= \|Ag\|_2^2 + \|Bg\|_2^2,
\end{align*}
and 
\begin{align*}
    Q_0(g)&= \|Ag\|_2^2= \sum_{i,j=1}^3\|D_i^+D_j^+g\|_{\ell^2(\Z^4)}^2\\
    & = \sum_{t=-L}^L
    \sum_{i,j=1}^3\|D_i^+D_j^+g^{(t)}\|_{\ell^2(\Z^3)}^2\notag\\
  &=\sum_{t=-L}^L
    \|\Delta g^{(t)}\|_{\ell^2(\Z^3)}^2.
\end{align*}
In the last equality, we used Lemma~\ref{lem:hessian} separately for each $t$.
Thus the law $\mu_A$ of the field $m$ defined by  $\mu_A(m=g)
  \propto
  \exp\!\left[-\frac\b2\|Ag\|^2\right]
$ is a product of independent three-dimensional integer-valued membrane measures, one for each $t\in[-L,L]$. Now, by Lemma~\ref{lem:add-squares}, we have
\begin{align*}
  \var_{\prob^{IMM}_{\b,\Lda}}(m(0))
  \le
  \var_{\mu_A}(m(0)) = \var_{\prob^{IMM}_{\b,\Lda^{(3)}}}(m(0)),
\end{align*}
where $\prob^{IMM}_{\b,\Lda^{(3)}}$ is the three-dimensional integer-valued membrane measure on $\Lda^{(3)}$. Now apply Proposition~\ref{prop:d3} to complete the proof.

\section{Proof of the lower bound: Theorem~\ref{thm:lower_bound_res}}
\label{sec:lower}
Throughout this section, $d=2$. Also, throughout the paper the distances are measured in the $\ell^\infty$ norm. 

For any function $u$, we define the Hessian by 
\begin{align}
(\HH u)_{ij}(x)&=D_i^+D_j^-u(x), \quad 1\leq i,\,j\leq d,\, x\in\Z^d.\notag
\end{align}
Also, for a function $u$ with finite support, we denote
\begin{align}
\cE(u)&=\|\Delta u\|_2^2.\notag
\end{align}
Then, for functions $u,v$ of finite support, by an argument similar to that in Lemma~\ref{lem:hessian}, we have
\begin{equation}\label{eq:d2-hessidentity}
 \sum_{x\in \Z^d}|\HH u(x)|^2=\cE(u),\qquad
 \sum_{x\in \Z^d}\sum_{1\leq i,j\leq d}(\HH u)_{ij}(x)(\HH v)_{ij}(x)
       =\la u,\Delta^2v\ra.
\end{equation}
Throughout, the norm of the Hessian $\HH u(x)$ is the Frobenius norm.

\subsection{Proof of Theorem~\ref{thm:lower_bound_res} using Theorem~\ref{thm:1.5_res}}
In this section we prove Theorem~\ref{thm:lower_bound_res} using Theorem~\ref{thm:1.5_res}, which is proved in Section~\ref{sec:thm 1.5}.
For any $N\ge 1$, define the Fej\'er kernel $F_N : \R \to \R$ by
\begin{align*}
F_N(x):= 1+ \sum_{q=1}^N 2\left(1-\frac{q}{N}\right) \cos({qx}).
\end{align*}
We define a measure
\begin{align*}
    \mathrm{d}\prob_{\b,\Lda, F_N}(\vr):=\frac{1}{Z_{\b,\Lda, F_N}(0)}\prod_{x\in\Lda}F_N(\vr(x))\mathrm{d} \prob^{MM}_{\b, \Lda}(\varphi),
\end{align*}
where
\begin{align*}
 Z_{\b,\Lda,F_N}(g)&=
 \E^{MM}_{\b,\Lda}\left[\prod_{x\in\Lda}F_N(\vr(x)+g(x))\right], \quad g:\Lda\to \R.
\end{align*}
The following theorem is analogous to Theorem 1.5 of \cite{K-P_2017} and Theorem 11 of \cite{Wirth_2019}.
\begin{theorem}\label{thm:1.5_res}
For any $\eps>0$, there exists $\b'_{\eps}>0$ such that the following holds. Let $0<\b<\b'_{\eps}$. Then for all $f:\Lda \to \R$ and $N\ge 1$
\begin{align*}
\E_{\b,\Lda, F_N}\left[ \e^{\la \vr, f\ra }\right] \ge \exp\left[\frac{1}{2(1+\eps)\b}\la \sigma_f, f \ra \right], 
\end{align*}
where $\sigma_f$ is the solution of \eqref{eq:bvp_sigma} with input $f$.
\end{theorem}
We now prove Theorem~\ref{thm:lower_bound_res} using the above theorem. 
\begin{proof}[Proof of Theorem~\ref{thm:lower_bound_res}]
We first show that 
\begin{align}\label{eq:lbprf1_res}
   \lim_{N\to \infty} \E_{\b/(2\pi)^2,\Lda, F_N}\left[ \e^{\la \vr, \frac{1}{2\pi}f\ra }\right] = \E_{\b,\Lda}^{IMM}\left[ \e^{\la m, f\ra }\right].
\end{align}
Let $\overline{\Lda}:=\Lda \cup \partial_2\Lda$ and
\begin{align*}
&\R^{\overline{\Lda}}_0:=\{g\in\R^{\overline{\Lda}} :\, g(x)=0, \, x\in\partial_2\Lda\} \\
&\Z^{\overline{\Lda}}_0:= \{g\in\R^{\overline{\Lda}}_0:\, g(x) \in\Z, \, x\in\Lda \}.
\end{align*}
For $g\in\Z^{\overline{\Lda}}_0$, let $\Omega_g\subseteq \R^{\overline{\Lda}}_0$ be the set of functions $\varphi$ satisfying $$\varphi(x) \in [-\pi+2\pi g(x), \pi+2\pi g(x)), \quad x\in\Lda.$$
Since $F_N$ is a positive summability kernel, we have, for $g\in\Z^{\overline{\Lda}}_0$ and any continuous function $\Psi: \R^{\overline{\Lda}}_0 \to \R$,
\begin{align}\label{eq:conv_sum_kernel_FN_res}
\lim_{N\to\infty} (2\pi)^{-|\Lda|} \int_{\Omega_g} \Psi(\varphi) \prod_{x\in\Lda} F_N(\varphi(x)) \mathrm{d}\varphi(x)\prod_{x\in\partial_2\Lda} \delta_0(\mathrm{d} \varphi(x)) = \Psi(2\pi g).
\end{align}
Abbreviating 
$$\mu_N(\De \vr):= \prod_{x\in\Lda} F_N(\vr(x))\prod_{x\in\Lda} \mathrm{d}\varphi(x) \prod_{x\in\partial_2\Lda} \delta_0(\mathrm{d} \varphi(x)),$$
we have
\begin{align*}
    \E_{\b/(2\pi)^2,\Lda, F_N}\left[ \e^{\la \vr, \frac{1}{2\pi}f\ra }\right] &=  \frac{ \int_{\R^{\overline{\Lda}}} \e^{\la \varphi, \frac1{2\pi} f\ra} \e^{-\frac{\b}{2(2\pi)^2}\la \varphi, \Delta^2\varphi\ra }\mu_N(\De \vr)} { \int_{\R^{\overline{\Lda}}} \e^{-\frac{\b}{2(2\pi)^2}\la \varphi, \Delta^2\varphi\ra }\mu_N(\De \vr)} \\
 &= \frac{ \sum_{g\in\Z^{\overline{\Lda}}_0}\int_{\Omega_g} \e^{\la \varphi, \frac1{2\pi} f\ra} \e^{-\frac{\b}{2(2\pi)^2}\la \varphi, \Delta^2\varphi\ra }\mu_N(\De \vr)}{ \sum_{g\in\Z^{\overline{\Lda}}_0}\int_{\Omega_g} \e^{-\frac{\b}{2(2\pi)^2}\la \varphi, \Delta^2\varphi\ra }\mu_N(\De \vr)}.
\end{align*}
For fixed $\beta,\Lambda,f$, one can use the positive definiteness of $\Delta^2_\Lda$ to obtain constants $c,C>0$ such that
\[
 \sup_{\varphi\in\Omega_g}
 \exp\left\{\ip{\varphi}{f/(2\pi)}
       -\frac{\beta}{2(2\pi)^2}\ip{\varphi}{\Delta^2\varphi}\right\}
 \le C \mathrm{e}^{-c\|g\|_2^2}.
\]
Now, by dominated convergence and~\eqref{eq:conv_sum_kernel_FN_res}, we conclude that the above expression converges to 
\begin{align*}
&\frac{ \sum_{g\in\Z^{\overline{\Lda}}_0} \e^{\la 2\pi g, \frac1{2\pi} f\ra} \e^{-\frac{\b}{2(2\pi)^2}\la 2\pi g, \Delta^2(2\pi g)\ra }}{ \sum_{g\in\Z^{\overline{\Lda}}_0} \e^{-\frac{\b}{2(2\pi)^2}\la 2\pi g, \Delta^2(2\pi g)\ra }}  = \E^{IMM}_{\b, \Lda} \left[ \e^{\la m, f\ra}\right]
\end{align*}
as $N$ goes to infinity. This completes the proof of~\eqref{eq:lbprf1_res}.

Consequently, by Theorem~\ref{thm:1.5_res} we have that
\begin{align*}
\E_{\b/(2\pi)^2, \Lda, F_N}\left[ \e^{\la \vr, \frac{1}{2\pi}f\ra }\right] \ge \exp\left[\frac{1}{2(1+\eps)\b} \la \sigma_f, f \ra \right] ,\quad \b<4\pi^2\b'_{\eps}.
\end{align*}
This, together with~\eqref{eq:lbprf1_res}, implies that
\begin{align*}
\E^{IMM}_{\b, \Lda} \left[ \e^{\la m, f\ra}\right] \ge \exp\left[\frac{1}{2(1+\eps)\b} \la \sigma_f, f \ra \right] , \quad \b<\b_\eps:=4\pi^2\b'_{\eps}.
\end{align*}

To prove the second part of the statement concerning the expectation of $\la m, f\ra^2$, let $\delta>0$ and $g=\delta f$. Then
\begin{align*}
\E^{IMM}_{\b, \Lda} \left[ \e^{\delta\la m, f\ra}\right] &\ge \exp\left[\frac{\delta^2}{2(1+\eps)\b} \la \sigma_f, f \ra \right] \\
& \ge 1+\frac{\delta^2}{2(1+\eps)\b} \la \sigma_f, f \ra.
\end{align*}
On the other hand, a Taylor expansion around $\delta=0$ gives
\begin{align*}
\E^{IMM}_{\b, \Lda} \left[ \e^{\delta\la m, f\ra}\right] &\le \E^{IMM}_{\b, \Lda} \left[ 1+ \delta \la m, f\ra + \frac{\delta^2}{2}\la m, f\ra^2 + \frac{\delta^3}{6}|\la m, f\ra|^3 \e^{|\delta\la m, f\ra|} \right]\\
& \le 1+ \frac{\delta^2}{2}\E^{IMM}_{\b, \Lda} \left[\la m, f\ra^2\right] + D(f)\delta^3.
\end{align*}
Thus
\begin{align*}
\E^{IMM}_{\b, \Lda} \left[ \la m, f\ra^2\right] + 2D(f)\delta \ge \frac{1}{(1+\eps)\b} \la \sigma_f, f\ra.
\end{align*}
Now, by letting $\delta\downarrow 0$, we get
\begin{align*}
\E^{IMM}_{\b, \Lda} \left[ \la m, f\ra^2\right]\ge\frac{1}{(1+\eps)\b}\la \sigma_f, f\ra.
\end{align*}
\end{proof}

\subsection{Proof of Theorem~\ref{thm:1.5_res} using Theorem~\ref{thm:1.6}}\label{sec:thm 1.5}
In this section we prove Theorem~\ref{thm:1.5_res} using Theorem~\ref{thm:1.6}. 

We start with a few definitions. A \textit{density} is a function $\r:\Z^2 \to \Z$ such that $\supp(\r) \neq \emptyset$ and $\supp(\r) \subseteq \Lda$. An \textit{ensemble} is a finite (possibly empty) collection of densities whose supports are mutually disjoint. For a density $\r$, define
\begin{align}
 Q(\r)&=\sum_{x\in\Z^2}\r(x),&
 P_i(\r)&=\sum_{x\in\Z^2} x_i\r(x),\quad i=1,2,\notag\\
 P(\r)&=(P_1(\r),P_2(\r)),&
 \mom(\r)&=(Q(\r),P_1(\r),P_2(\r)).\label{d2:eq:moments}
\end{align}
Observe that for an affine function $f(x)=a+b_1x_1+b_2x_2$,
\begin{equation}\label{d2:eq:affinecancellation}
 \la f,\r\ra=aQ(\r)+b_1P_1(\r)+b_2P_2(\r).
\end{equation}
Hence $\la f,\r\ra=0$ for every affine $f$ if and only if $\mom(\r)=0$.
We call $\r$ \emph{moment-neutral} if $\mom(\r)=0$, and
\emph{moment-active} otherwise. Write
\begin{equation}\label{d2:eq:momentclasses}
 \cN_0=\{\r\in\cN:\mom(\r)=0\},\qquad \cN_c=\cN\setminus\cN_0.
\end{equation}
In particular, $\cN_c$ includes densities with $Q(\r)=0$ and $P(\r)\ne0$. For any $A, B\subset \Z^2$, define 
\begin{align*}
    \diam(A)&:=\sup\{|x-y|_\infty: x,y\in A\},\\
    \dist(A, B)&:= \inf\{|x-y|_\infty: x\in A,\,y\in B\}.
\end{align*}
Set
\begin{equation}\label{d2:eq:radius}
\begin{aligned}
 d(\r)&=\max\{1,\diam(\supp(\r))\},&
 \delta(\r)&=\dist(\r,\Lda^c),\\
 d_\Lda(\r)&=
 \begin{cases}
 d(\r),&\mom(\r)=0,\\
 \max\{d(\r),\delta(\r)\},&\mom(\r)\ne0.
 \end{cases}
\end{aligned}
\end{equation}
Above and throughout the rest of the paper, $\dist(\r,A) = \dist(\supp(\r),A)$ for $A\subseteq \Z^2.$ Also, we write $\dist(\r_1,\r_2)$ for $\dist(\supp(\r_1),\supp(\r_2))$.

Note that $d_\Lda(\r)\ge 1$ for any $\r$.
For each density $\r$, we choose an arbitrary fixed $x_\r \in \supp(\r)$. For any density $\r$, we write $$\|\r \|_2:=\sqrt{\sum_{x\in\Z^2}\r(x)^2}.$$
The constants
\begin{align}\label{eq:def_alpha_M}
\frac32 < \a <2,\, M=10^{12}
\end{align}
are fixed throughout.

The following theorem is analogous to Theorem 1.6 of \cite{K-P_2017} or Theorem 13 of \cite{Wirth_2019}.

\begin{theorem}\label{thm:1.6}
There exist constants $\b_1, c_1>0$ such that the following holds. There exist:
\begin{enumerate}
    \item a finite collection of ensembles $\cF$,
    \item positive coefficients $(c_{\cN})_{\cN\in\cF}$ summing to $1$,
    \item real coefficients $(z(\b,\r,\cN))_{\r\in\cN,\cN\in\cF}$,
    \item functions $a_\r:\Lda\to\R$ for each $\r\in\cN,\cN\in\cF$,
\end{enumerate}
such that for any $g:\Lda \to\R$
\begin{align*}
\E_{\b,\Lda}^{MM}\left[ \prod_{x\in\Lda} F_N(\vr(x)+g(x))\right] = \sum_{\cN \in \cF} c_{\cN} \E_{\b,\Lda}^{MM}\left[ \prod_{\r\in\cN} \left[ 1+ z(\b,\r,\cN)\cos\left(\la \vr, \r-\b\Delta_\Lda^2a_\r \ra + \la g, \r \ra \right)\right]\right].
\end{align*}
Moreover, the following properties are satisfied for every $\cN \in \cF$:
\begin{enumerate}
    \item[(i)]\label{item:thm1.6.1} For distinct $\r_1,\r_2\in \cN$,  $$\dist(\r_1, \r_2) \ge M\min\{d_\Lda(\r_1),d_\Lda(\r_2)\}^{\a}. $$
    \item[(ii)] For every $\r\in\cN$
    \begin{align*}
        |z(\b,\r,\cN)| \le \exp \left[-\frac{c_1}{\b}\left( \|\r\|_2^2 + \log_2(d_\Lda(\r)+1)\right) \right], \qquad 0<\b<\b_1.
    \end{align*}

\end{enumerate}
\end{theorem}

\begin{proof}[Proof of Theorem~\ref{thm:1.5_res}]
%Recall that $\sigma:\Lda\to\R$ be the solution of 
%\begin{align*}
%    \Delta^2\sigma=\frac{1}{\b}f \qquad \mbox{on } \Lda^{\circ}\\
%    \sigma=0 \qquad \mbox{on } \partial_2\Lda
%\end{align*}
Denote $\sigma=\sigma_f$. We have 
\begin{align*}
 \E_{\b,\Lda, F_N}\left[ \e^{\la \vr, f\ra }\right]&= \frac{1}{Z_{\b,\Lda, F_N}(0)}\int \e^{\la \vr, f\ra}\prod_{x\in\Lda}F_N(\vr(x))\mathrm{d} \prob^{MM}_{\b, \Lda}(\varphi)\\
 &= \frac{1}{Z_{\b,\Lda, F_N}(0) Z_{\b,\Lda}^{MM}}\int \prod_{x\in\Lda}F_N(\vr(x)) \e^{\la \vr, f\ra -\frac{\b}{2}\la \vr, \Delta^2\vr\ra} \prod_{x\in\Lda}\mathrm{d}\vr(x) \prod_{x\in\Z^2\setminus\Lda} \delta_0(\mathrm{d}\vr(x)).
\end{align*}
By the change of variables $\vr\to\vr+\sigma/\b$, we have
\begin{align}\label{eq:t1.5.1_res}
 \E_{\b,\Lda, F_N}\left[ \e^{\la \vr, f\ra }\right]&= \frac{1}{Z_{\b,\Lda, F_N}(0)Z_{\b,\Lda}^{MM}}\int \prod_{x\in\Lda}F_N\left(\vr(x) + \frac{\sigma(x)}{\b}\right) \exp\left[\la \vr, f\ra +\frac{1}{\b}\la\sigma, f\ra -\frac{\b}{2}\la \vr, \Delta^2\vr\ra \right.\nonumber\\
 &-\left.\frac{1}{2}\la \sigma, \Delta^2\vr\ra - \frac{1}{2}\la \vr, \Delta^2\sigma\ra -\frac{1}{2\b}\la\sigma, f\ra\right] \prod_{x\in\Lda}\mathrm{d}\vr(x) \prod_{x\in\Z^2\setminus \Lda} \delta_0(\mathrm{d}\vr(x))\nonumber\\
 &=\frac{1}{Z_{\b,\Lda, F_N}(0)} \e^{\frac{1}{2\b}\la\sigma, f\ra} Z_{\b,\Lda, F_N}(\sigma/\b).
\end{align}

Our goal is to estimate $Z_{\b,\Lda, F_N}(\sigma/\b)/Z_{\b,\Lda, F_N}(0)$. By Theorem~\ref{thm:1.6}, we have
\begin{align*}
  Z_{\b,\Lda, F_N}(\sigma/\b)=\sum_{\cN\in\cF} c_{\cN}Z_{\cN}(\sigma/\b), 
\end{align*}
where
\begin{align*}
 Z_{\cN}(\sigma/\b):=\int \prod_{\r\in\cN}\left[ 1+ z(\b,\r,\cN)\cos\left(\la \vr, \r-\b\Delta^2a_\r \ra + \la \sigma/\b, \r \ra \right) \right] \mathrm{d} \prob^{MM}_{\b, \Lda}(\varphi)  .
\end{align*}
First we use the following trigonometric claim:
\begin{claim}[Claim 3.1 of \cite{K-P_2017} or Claim 14 of \cite{Wirth_2019}]\label{claim:3.1_res}There exists an absolute constant $D>0$ such that for all $x,y\in\R$ and $|z|<1/8$,
\begin{align*}
    1+z\cos(x+y) \ge \exp\left(-\frac{z\sin x \sin y}{1+z\cos x}-D|z|y^2 \right)\left(1+z\cos x\right).
\end{align*}
\end{claim}
By Theorem~\ref{thm:1.6}, we can choose $\b$ small enough that $|z(\b,\r,\cN)|<1/8$.

Now, applying Claim~\ref{claim:3.1_res} and using Jensen's inequality, we have, for each $\cN\in\cF$,
\begin{align*}
    \frac{Z_{\cN}(\sigma/\b)}{Z_{\cN}(0)}&=\frac{1}{Z_{\cN}(0)}\int \prod_{\r\in\cN}\left[ 1+ z(\b,\r,\cN)\cos\left(\la \vr, \r-\b\Delta^2a_\r \ra + \la \sigma/\b, \r \ra \right) \right] \mathrm{d} \prob^{MM}_{\b, \Lda}(\varphi)\\
    & \ge \exp\left[ -D\sum_{\r\in\cN}|z(\b,\r,\cN)| \la \sigma/\b, \r \ra^2 \right]
    Z_{\cN}(0)^{-1} \int \prod_{\r\in\cN}\left[1+z(\b,\r,\cN)\cos\la \vr, \r-\b\Delta^2a_\r \ra \right]\\
    & \quad \quad \exp\left[-\sum_{\r\in\cN} \frac{z(\b,\r,\cN)\sin\la \vr, \r-\b\Delta^2a_\r \ra\sin\la \sigma/\b, \r \ra}{1+z(\b,\r,\cN)\cos\la \vr, \r-\b\Delta^2a_\r \ra}\right]\mathrm{d} \prob^{MM}_{\b, \Lda}(\varphi)\\
    & \ge \exp\left[-D\sum_{\r\in\cN}|z(\b,\r,\cN)| \la \sigma/\b, \r \ra^2\right].
\end{align*}
The proof will be complete if we can prove the following claim, which is analogous to Claim 3.2 of \cite{K-P_2017} and Claims 15--16 of \cite{Wirth_2019}.
\begin{claim}\label{claim:3.2_res}
Let $D_1>0$. Then there exists $0<\b_2\le \b_1$ such that
\begin{align*}
  \sum_{\r\in\cN}|z(\b,\r,\cN)| \la \sigma, \r \ra^2 \le \frac{\b}{D_1} \la \sigma, \Delta^2\sigma \ra, \quad \b<\b_2,\, \cN\in\cF.
\end{align*}
\end{claim}
We now show how Theorem~\ref{thm:1.5_res} follows from Claim~\ref{claim:3.2_res}. Indeed, for any $\eps>0$, applying the claim with 
$$D_1:=\frac{2(1+\eps)}{\eps}D,$$
we have that there exists $\b'_{\eps}>0$ such that for $0<\b<\b'_{\eps}$ and for each $\cN\in\cF$
\begin{align*}
    \frac{Z_{\cN}(\sigma/\b)}{Z_{\cN}(0)} & \ge \exp\left[-D\sum_{\r\in\cN}|z(\b,\r,\cN)| \la \sigma/\b, \r \ra^2\right]\\
    &\ge \exp\left[-\frac{\eps}{2(1+\eps)\b} \la \sigma, \Delta^2\sigma \ra \right].
\end{align*}
Hence
\begin{align*}
Z_{\b,\Lda, F_N}(\sigma/\b)/Z_{\b,\Lda, F_N}(0) & =\frac{\sum_{\cN\in\cF} c_{\cN}Z_{\cN}(\sigma/\b)}{\sum_{\cN\in\cF} c_{\cN}Z_{\cN}(0)}\\
& \ge \exp\left[-\frac{\eps}{2(1+\eps)\b} \la \sigma, \Delta^2\sigma \ra \right].
\end{align*}

Therefore, from \eqref{eq:t1.5.1_res}, we conclude that for any $\eps>0$, there exists $\b'_{\eps}>0$ such that for $0<\b<\b'_{\eps}$
\begin{align*}
 \E_{\b,\Lda, F_N}\left[ \e^{\la \vr, f\ra }\right] \ge \exp\left[\frac{1}{2(1+\eps)\b} \la \sigma, f \ra \right], 
\end{align*}
concluding the proof of Theorem~\ref{thm:1.5_res}.\end{proof}
We now prove Claim~\ref{claim:3.2_res}. We will show two separate bounds on
\[
\sum_{\r\in\cN}|z(\b,\r,\cN)| \la \sigma, \r \ra^2
\]
by splitting the sum into contributions from moment-neutral (Claim~\ref{claim:3.21_res}) and moment-active (Claim~\ref{claim:3.22_res}) charges. Together, they imply Claim~\ref{claim:3.2_res}.

We first give the local estimates for $\cN_0$ and $\cN_c$. For any $x\in\Z^2$ and $R\ge 1$, we define 
\begin{align*}
 \bx_R(x):=\{y\in\Z^2:|y-x|_\infty\le R\}.
\end{align*}
\begin{lemma}\label{d2:lem:pairing}
There exists an absolute constant $C$ such that, for every density
$\r$,
\begin{equation}\label{d2:eq:localpair}
 |\la\sigma,\r\ra|^2\le
 C d_\Lda(\r)^{C_1}\|\r\|_1^2
 \sum_{x\in\bx_{C_0\,d_\Lda(\r)}(x_\r)}|\HH\sigma(x)|^2,
\end{equation}
where $C_0=3,\, C_1=4$ if $\mom(\r)=0$, and $C_0=9,\, C_1=6$ if $\mom(\r)\neq 0$.
\end{lemma}

\begin{proof}
For $z\in\Z^2$ and $R\ge1$, define
\[
 A_{\sigma,z}(x)=\sigma(z)+\sum_{i=1}^2(x_i-z_i)D_i^+\sigma(z),
 \qquad
 S_R=\max_{x\in\bx_{3R}(z)}|\HH\sigma(x)|.
\]
Observe that
\[
 D_i^+\sigma(y+e_j)-D_i^+\sigma(y)
       =(\HH\sigma)_{ij}(y+e_j).
\]
Therefore, for any $y\in\bx_R(z)$
\[
 |D_i^+\sigma(y)-D_i^+\sigma(z)|\le4RS_R.
\]
Let $x\in\bx_R(z)$ and $z=y_0,\ldots,y_m=x$ be a path from $z$ to $x$ with $m\le2R$. Write
$y_{\ell+1}-y_\ell=\varepsilon_\ell e_{i_\ell}$, where $e_{i_\ell}\in\{e_1, e_2\}$ and put
$w_\ell=y_\ell$ if $\varepsilon_\ell=1$, and
$w_\ell=y_{\ell+1}$ otherwise. Then
\begin{align*}
 \sigma(x)-A_{\sigma,z}(x)
 &=\sum_{\ell=0}^{m-1}\varepsilon_\ell
   \left(D_{i_\ell}^+\sigma(w_\ell)-D_{i_\ell}^+\sigma(z)\right).
\end{align*}   
Hence
\begin{align}
 \sup_{x\in\bx_R(z)}|\sigma(x)-A_{\sigma,z}(x)|
 &\le8R^2S_R
 \le8R^2\left(\sum_{x\in\bx_{3R}(z)}|\HH\sigma(x)|^2\right)^{1/2}.
       \label{d2:eq:taylorlocal}
\end{align}

Now suppose $\mom(\r)=0$. Since $ A_{\sigma,z}$ is affine, we have 
\begin{align}\label{eq:affinecancellation}
    \la  A_{\sigma,z}, \r\ra=0.
\end{align}

Take $z=x_\r$ and $R=d_\Lda(\r)=d(\r)$.
By \eqref{eq:affinecancellation},
\begin{align}
 |\la\sigma,\r\ra|^2
 &=|\la\sigma-A_{\sigma,x_\r},\r\ra|^2\notag\\
 &\le64d_\Lda(\r)^4\|\r\|_1^2
       \sum_{x\in\bx_{3d_\Lda(\r)}(x_\r)}|\HH\sigma(x)|^2.
       \label{d2:eq:neutralpair}
\end{align}

Next, suppose $\mom(\r)\ne0$. In this case, $d_\Lda(\r)=\max\{d(\r), \delta(\r)\}$. Choose $y\in\supp(\r)$ with
$\dist(y,\Lda^c)=\delta(\r)$ and a nearest exterior point $a$
across a face of $\Lda$. Let $u$ be the outward coordinate unit
vector and $v$ a tangential coordinate unit vector. Then
$a,a+u,a+v\in\Lda^c$. Write
\[
 P(\r)-Q(\r)a=p_u u+p_v v,
\]
and define
\begin{equation}\label{d2:eq:boundarycompletion}
 \widehat\r=\r+(-Q(\r)+p_u+p_v)\delta_a
                  -p_u\delta_{a+u}-p_v\delta_{a+v}.
\end{equation}
Since $p_u,p_v\in\Z$,
\begin{align*}
 Q(\widehat\r)&=Q(\r)-Q(\r)+p_u+p_v-p_u-p_v=0,\\
 P(\widehat\r)&=P(\r)-aQ(\r)-p_u u-p_v v=0,\\
 \la\sigma,\widehat\r\ra&=\la\sigma,\r\ra.
\end{align*}
Note that for any $x\in\supp(\r)$, we have $|a-x|_\infty\le d(\r)+\delta(\r)\le2 d_\Lda(\r)$. Hence $\supp(\widehat\r)\subset\bx_{3d_\Lda(\r)}(x_\r)$. Also,
\begin{align*}
 |p_u|+|p_v|&\le\sum_x |x-a|_1|\r(x)|\le 4d_\Lda(\r)\|\r\|_1,
\end{align*}
and hence $ \|\widehat\r\|_1\le Cd_\Lda(\r)\|\r\|_1$ for some constant $C$. We now apply \eqref{d2:eq:taylorlocal} to $\widehat\r$, with $R=3d_\Lda(\r)$,
and use $\la A_{\sigma,x_\r},\widehat\r\ra=0$ from \eqref{eq:affinecancellation} to obtain \eqref{d2:eq:localpair}, with $C_0=9$ and $C_1=6$.
\end{proof}

 Now we prove the following two claims to complete the proof of Claim~\ref{claim:3.2_res}.
\begin{claim}[Control of moment-neutral densities]\label{claim:3.21_res}
Let $D_1>0$. Then there exists $\b_3>0$ such that
\begin{align*}
  \sum_{\r\in\cN_0}|z(\b,\r,\cN)| \la \sigma, \r \ra^2 \le \frac{\b}{2D_1} \la \sigma, \Delta^2\sigma \ra, \quad \b<\b_3,\, \cN\in\cF.
\end{align*}
\end{claim}
\begin{proof}[Proof of Claim~\ref{claim:3.21_res}]
For any $\r\in\cN$, let $B_\r:=\bx_{3d_\Lda(\r)}(x_\r)$. We know from Theorem~\ref{thm:1.6} that, for distinct $\r_1,\r_2\in \cN$,  $$\dist(\r_1, \r_2) \ge M\min\{d_\Lda(\r_1),d_\Lda(\r_2)\}^{\a}. $$
Therefore, if $\r_1\neq \r_2$ and $2^k-1 \le d_\Lda(\r_1), d_\Lda(\r_2)<2^{k+1}-1$, then $B_{\r_1}\cap B_{\r_2}=\emptyset$.
Now, using \eqref{d2:eq:localpair}, we obtain for any $\r\in\cN_0$ and $\b<\b_1$ 
\begin{align*}
    |z(\b,\r,\cN)| \la \sigma, \r \ra^2 &\le C \exp\left[-\frac{c_1}\b \left( \|\r\|_2^2 + \log_2(d_\Lda(\r)+1)\right)\right]d_\Lda(\r)^{4}\|\r\|_1^2
 \sum_{x\in B_\r}|\HH\sigma(x)|^2\\
    & \le C \|\r\|_2^4 \exp\left[-\frac{c_1}\b \|\r\|_2^2\right] d_\Lda(\r)^4\exp\left[-\frac{c_1}\b\log_2(d_\Lda(\r)+1)\right]\sum_{x\in B_\r}|\HH\sigma(x)|^2.
\end{align*}
Hence, by choosing $0<\b_3<\b_1$, we obtain for each $\b<\b_3$ that
\begin{align*}
    |z(\b,\r,\cN)| \la \sigma, \r \ra^2 &\le \frac{\b}{2D_1(d_\Lda(\r)+1)}\sum_{x\in B_\r}|\HH\sigma(x)|^2.
\end{align*}
Now we conclude 
\begin{align*}
\sum_{\r\in\cN_0}|z(\b,\r,\cN)| \la \sigma, \r \ra^2 &
%\le \sum_{\r\in\cN_0} \frac{\b}{DC(d_\Lda(\r)+1)}\sum_{j,\ell\in D(\r), j\sim \ell}(\sigma_j-\sigma_\ell)^2\\
\le \frac{\b}{2D_1} \sum_{k\ge 1} 2^{-k} \sum_{\r\in\cN_0, 2^k\le d_\Lda(\r)+1 <2^{k+1}} \sum_{x\in B_\r}|\HH\sigma(x)|^2\\
& \le \frac{\b}{2D_1} \sum_{k\ge 1} 2^{-k}\sum_{x\in \Z^2}|\HH\sigma(x)|^2 \\
&\le \frac{\b}{2D_1} \sum_{x\in \Z^2}|\HH\sigma(x)|^2. 
\end{align*}
Thus
\begin{align*}
\sum_{\r\in\cN_0}|z(\b,\r,\cN)| \la \sigma, \r \ra^2 \le \frac{\b}{2D_1} \la \sigma, \Delta^2\sigma \ra.\qedhere
\end{align*}
\end{proof}

\begin{claim}[Control of moment-active densities]\label{claim:3.22_res}
Let $D_1>0$. Then there exists $\b_4>0$ such that
\begin{align*}
  \sum_{\r\in\cN_c}|z(\b,\r,\cN)| \la \sigma, \r \ra^2 \le \frac{\b}{2D_1} \la \sigma, \Delta^2\sigma \ra, \quad \b<\b_4,\, \cN\in\cF.
\end{align*}
\end{claim}
\begin{proof}
The proof is similar to that of Claim~\ref{claim:3.21_res} and is therefore omitted. 
\end{proof}

\subsection{Proof of Theorem~\ref{thm:1.6}}
\label{d2:sec:gaussian}
The proof of Theorem~\ref{thm:1.6} is completed by proving Theorem~\ref{thm:2.2} and Theorem~\ref{thm:2.3}. Theorem~\ref{thm:1.6} is an immediate consequence of these two theorems. We use a convex expansion and a measure modification, as in
Theorems~2.2--2.3 of \cite{K-P_2017}.

\subsubsection{Covering densities with squares}
Let $k\ge 0$ be an integer. We call $s$ a square of side length $2^k$ if there exists a vertex $(a_1,a_2)\in \Z^2$ such that
\begin{align*}
s=\{(b_1, b_2):a_i-b_i \in \{0,\ldots, 2^k-1\} \text{ for } i\in \{1,2\}\}.
\end{align*}
For $k\ge 0$, let $\cS_k(\r)$ be a minimal collection of squares of side length $2^k$ covering the support of $\r$. 
%The precise choice of $\cS_k(\r)$, when more than one minimal collection exists will be described later in Section~\ref{sec:prop 2.4}.

Next we define $\cS^{sep}_k(\r)$ for $k\ge 1$ as follows.
\begin{align*}
\cS^{sep}_k(\r)\coloneqq\begin{cases}\{s\in \cS_k(\r):  & \dist(s,s^{'}) \ge 2M\left(1000\cdot 2^k\right)^{\a}, \\ &\,\forall \,s^{'}\in \cS_k(\r)\setminus \{s\}\} \,\,\text{ if } |\cS_k(\r)|>1\\
\emptyset, & \qquad\qquad\qquad\qquad\quad\text{if } |\cS_k(\r)|=1 \text{ and } \mom(\r)=0\\
\cS_k(\r), & \qquad\qquad\qquad\qquad\quad \text{if } |\cS_k(\r)|=1, \mom(\r)\neq 0,\\
&\qquad\qquad\qquad\qquad\quad\text{ and }\delta(\r) \ge 1000 \cdot 2^{k}\\
\emptyset &  \qquad\qquad\qquad\qquad\quad\text{otherwise. }
\end{cases}
\end{align*}

Define $\cS^{"}_{k}(\r):=\cS_{k}(\r)\setminus \cS^{sep}_{k}(\r)$ and set $\cS^{sep}_{0}(\r)=\cS_{0}(\r)$. Also define
\begin{align}\label{eq:def_nrho}
&n(\r) :=\lceil \log_2(M. d_\Lda(\r)^\a)\rceil,\\
& A(\r):=\sum_{k=0}^{n(\r)} |\cS_k(\r)|.\nonumber\\
&I(\r):=\sum_{k=8}^{n(\r)} |\cS^{sep}_k(\r)|.
\end{align}
The following proposition allows us to control the size of $A(\r)$. This is analogous to Proposition 20 of \cite{Wirth_2019} and Proposition 2.1 of \cite{K-P_2017}. We prove this in Appendix~\ref{d2:sec:regrouping}.
\begin{proposition}[]\label{prop:2.1}
There exists a positive absolute constant $D_1$ such that for any density $\r$
\begin{align*}
\log_2(d_\Lda(\r) + 1) \le A(\r) \le D_1 \left( |\cS_0(\r)| + I(\r) \right).
\end{align*}
\end{proposition}

\subsubsection{Expansion as a convex combination}

Below and throughout the rest of the paper, we write $\r_1\subset \r$ if $\supp(\r_1)\subset \supp(\r)$ and $\r_1(x)=\r(x)$ for all $x\in\supp(\r_1)$. The following theorem is analogous to Theorem 21 of \cite{Wirth_2019} and Theorem 2.2 of \cite{K-P_2017}. The proof is also similar. We give the proof in Appendix~\ref{d2:sec:regrouping}.

\begin{theorem}[]\label{thm:2.2}
There exist a positive constant $D_2$, a finite collection of ensembles $\cF$, positive coefficients $(C_\cN)_{\cN\in\cF}$ summing to $1$ and real coefficients $(K(\r))_{\r\in\cN\in\cF}$ such that, for every $\psi:\Lda \to \R$,
\begin{align*}
    \prod_{x\in\Lda} F_N(\psi(x)) = \sum_{\cN\in\cF} C_{\cN} \prod_{\r\in\cN}\left[1+ K(\r)\cos(\la \psi, \r \ra)\right].
\end{align*}
Additionally, the following properties are satisfied for each $\cN\in\cF$:
\begin{enumerate}[label=(\arabic*),ref=(\arabic*)]
\item\label{item:one_6.4} If $\r, \r^{'}\in\cN$ are distinct, then
\begin{align*}
    \dist(\r, \r^{'}) \ge M [\min\{d_\Lda(\r), d_\Lda(\r^{'})\}]^\a.
\end{align*}
\item\label{item:two_6.4} If $0\neq \r_1\subseteq \r\in\cN$ with $\r_1\neq \r$ satisfies $\dist(\r_1, \r-\r_1) \ge 2Md(\r_1)^\a$, then 
\begin{align*}
   \mom(\r_1)\neq 0, \qquad 2M \delta(\r_1)^\a>\dist(\r_1, \r-\r_1).
\end{align*}
\item\label{item:three_6.4} The coefficients $K(\r)$ satisfy
\begin{align*}
    |K(\r)| \le e^{D_2A(\r)} \prod_{x\in \supp(\r)}e^{\r(x)^2}\left|\left(1-\frac{|\r(x)|}{N}\right)\right|.
\end{align*}

\end{enumerate}
\end{theorem}

\subsubsection{Measure modification}
The following theorem is in line with Theorem~2.3 of \cite{K-P_2017} and Theorem~22 of \cite{Wirth_2019}.

\begin{theorem}\label{thm:2.3}
There exists an absolute constant $D_3>0$ such that the following holds. Let $\cN$ be an ensemble satisfying properties~\ref{item:one_6.4}--\ref{item:two_6.4} of Theorem~\ref{thm:2.2} and let $(K(\r))_{\r\in\cN}$ be real coefficients. Then there exist $z(\b,\r,\cN)$ and functions $a_{\r, \cN}:\Lda\to\R$ such that for every $g:\Lda\to\R$
\begin{align*}
    &\int\prod_{\r\in\cN}\left[1+K(\r)\cos\left( \la \vr, \r\ra + \la g, \r\ra \right)\right] \mathrm{d} \prob^{MM}_{\b, \Lda}(\varphi)\\
    & = \int\prod_{\r\in\cN}\left[1+z(\b, \r,\cN)\cos\left( \la \vr, \r - \b\Delta^2a_{\r,\cN}\ra + \la g, \r\ra \right)\right] \mathrm{d} \prob^{MM}_{\b, \Lda}(\varphi).
\end{align*}
Additionally,
\begin{align}\label{eq:t2.3.2}
    |z(\b, \r,\cN)| \le |K(\r)| \exp\left[-\frac1\b \left( \frac1{144}\|\r\|_2^2 + D_3 \sum_{k=8}^{n(\r)} |\cS^{sep}_k(\r)|\right)\right].
\end{align}

If the coefficients $(K(\r))_{\r\in\cN}$ satisfy condition~\ref{item:three_6.4} of Theorem~\ref{thm:2.2}, then there exist $\b_0, C_2>0$ such that
\begin{align}\label{eq:t2.3.1}
    |z(\b, \r,\cN)| \le \exp\left[-\frac{C_2}\b \left( \|\r\|_2^2 + \log_2(d_\Lda(\r)+1)\right)\right], \quad 0<\b<\b_0.
\end{align}
%\textcolor{red}{\begin{align*}
%    |z(\b, \r,\cN)| \le \exp\left[-\frac{C_2}\b \left( \|\r\|_2^2 + \log_2(d_\Lda(\r)+1)\right)\right], \quad 0<\b<\b_0(L).
%\end{align*}}
\end{theorem}

The proof of Theorem~\ref{thm:2.3} is completed in the following three steps. We fix an ensemble $\cN$ satisfying properties~\ref{item:one_6.4}--\ref{item:two_6.4} of Theorem~\ref{thm:2.2}, a sequence $(K(\r))_{\r\in\cN}$ and $g:\Lda\to\R$. For a density $\r$ and a function $a:\Lda\to\R$, we also define the energy
\begin{align*}
    E_{\b}(a,\r):= \la a, \r\ra -\frac{\b}2\la a, \Delta^2a\ra.
\end{align*}

\paragraph{Step 1: constructing a spin wave.}

For $\r\in\cN$, let $\cN(\r)$ be the following subensemble:
\begin{align*}
    \cN(\r):=\{\eta\in\cN: d_{\Lda}(\eta) \le 2 d_{\Lda}(\r), \eta\neq \r\}.
\end{align*}

The following proposition is analogous to Proposition 23 of \cite{Wirth_2019} and Proposition 2.4 of \cite{K-P_2017}. 
\begin{proposition}\label{prop:2.4}
    There exists an absolute constant $D_3>0$ such that the following holds. For each $\r\in\cN$, there exists a function $a_{\r,\cN}$, denoted by $a_\r$ for clarity of notation, such that the following hold:
    \begin{enumerate}[label=(\arabic*),ref=(\arabic*)]
    \item \label{item:prop1_2.4}$\supp(a_\r) \subseteq \Lda \cap \bx_{10d_\Lda(\r)}(x_\r) $. 
        \item\label{item:prop2_2.4} $a_\r$ is affine on $\bx_{100d_\Lda(\eta)}(x_\eta)$ for
$\eta\in\cN(\r)$.
        \item\label{item:prop3_2.4} $\supp(a_\r)\cap\supp(\eta)=\varnothing$ for
$\eta\ne\r$ with $\mom(\eta)\ne0$.
\item\label{item:prop4_2.4} For $\eta\ne\r$,
\begin{equation*}
 \la a_\r,\eta\ra=0,\qquad
 \la a_\r,\Delta_\Lda^2a_\eta\ra=0.
\end{equation*}
        \item\label{item:prop5_2.4} $E_{\b}(a_\r,\r) \ge \frac1\b \left( \frac1{144}\|\r\|_2^2 + D_3 \sum_{k=8}^{n(\r)} |\cS^{sep}_k(\r)|\right)$.
        \end{enumerate}
\end{proposition}

We call $a_\r$ a spin wave associated with $\r$. We prove this proposition in Section~\ref{sec:prop 2.4}.

\paragraph{Step 2: complex translation by spin waves.}
In this step, we will use the following Gaussian integral identity: for $a,\tau:\Z^2\to\R$ with $\supp(a)\subseteq\Lda$
\begin{align}\label{eq:gaussian_integral}
 \int \e^{\iota\la \vr,\tau\ra}  \mathrm{d} \prob^{MM}_{\b, \Lda}(\varphi) = \e^{-E_{\b}(a,\tau)} \int \e^{\iota\la \vr,\tau-\b\Delta^2a\ra}  \mathrm{d} \prob^{MM}_{\b, \Lda}(\varphi).
\end{align}
To prove the above identity, we observe that
\begin{align*}
     \e^{-E_{\b}(a,\tau)}& \int \e^{\iota\la \vr,\tau-\b\Delta^2a\ra}  \mathrm{d} \prob^{MM}_{\b, \Lda}(\varphi)\\
     &= \exp\left(-\left(\la a, \tau\ra -\frac{\b}2\la a, \Delta^2a\ra\right)-\frac12\sum_{x\in\Lda} (\tau-\b\Delta^2a)(x)\left(\frac1{\b}\sum_{y\in\Lda}(\Delta^2_{\Lda})^{-1}(y,x)(\tau-\b\Delta^2a)(y)\right)\right)\\
    &= \exp\left(-\left(\la a, \tau\ra -\frac{\b}2\la a, \Delta^2a\ra\right) - \frac12\sum_{x\in\Lda} (\tau-\b\Delta^2a)(x) \left(\frac1{\b}\sum_{y\in\Lda}(\Delta^2_{\Lda})^{-1}(y,x)\tau(y)- a(x)\right)\right)\\
%    &= \exp\left(-\frac12\sum_{x\in\Lda} \tau(x)\cdot \left(\frac1{\b}\sum_{y\in\Lda}(\Delta^2_{\Lda})^{-1}(y,x)\tau(y)\right)-\la a, \tau \ra+\frac{\b}2 \la a, \Delta^2a\ra +\frac12\la a, \tau \ra\right.\\
     &\left.  +\frac12\sum_{x\in\Lda} \Delta^2a(x) \sum_{y\in\Lda}(\Delta^2_{\Lda})^{-1}(y,x)\tau(y) -\frac{\b}2 \la\Delta^2a,\, a\ra\right)= \int \e^{\iota\la \vr,\tau\ra}  \mathrm{d} \prob^{MM}_{\b, \Lda}(\varphi).
\end{align*}

For each $\r\in\cN$, let $a_\r$ be the function given by Proposition~\ref{prop:2.4}. We write
\begin{align*}
    &\int \prod_{\r\in\cN} \left[ 1+ K(\r)\cos\left(\la\vr,\r\ra + \la g,\r\ra\right)\right]\mathrm{d} \prob^{MM}_{\b, \Lda}(\varphi)\\
    &=\int \prod_{\r\in\cN} \left[ 1+ \frac{K(\r)}2\left[\exp\left(\iota\left(\la\vr,\r\ra + \la g,\r\ra\right)\right) + \exp\left(-\iota\left(\la\vr,\r\ra + \la g,\r\ra\right)\right)\right]\right]\mathrm{d} \prob^{MM}_{\b, \Lda}(\varphi)\\
    &=\int\sum_{\zeta:\cN\to\{-1,0,1\}}\prod_{\r\in\cN}\left(\frac{K(\r)}2\right)^{|\zeta(\r)|} \e^{\iota \zeta(\r)\left(\la\vr,\r\ra + \la g,\r\ra\right)}\mathrm{d} \prob^{MM}_{\b, \Lda}(\varphi)\\
    &=\sum_{\zeta:\cN\to\{-1,0,1\}} \left[\prod_{\r\in\cN}\left(\frac{K(\r)}2\right)^{|\zeta(\r)|}\right]\int\exp\left[\iota \left(\la \vr, \sum_{\r\in\cN}\zeta(\r)\r\ra + \la g, \sum_{\r\in\cN}\zeta(\r)\r\ra\right)\right]\mathrm{d} \prob^{MM}_{\b, \Lda}(\varphi).
\end{align*}
Fix $\zeta:\cN\to\{-1,0,1\}$. Let
\begin{align*}
    a_\zeta:=\sum_{\r\in\cN}\zeta(\r)a_\r, \qquad \tau_\zeta:=\sum_{\r\in\cN}\zeta(\r)\r.
\end{align*}
Note that $\supp(a_\zeta)\subseteq\Lda$ by property~\ref{item:prop1_2.4} of Proposition~\ref{prop:2.4}. Now, applying \eqref{eq:gaussian_integral}, we have
\begin{align*}
    &\int \exp\left[\iota\left( \la \vr, \tau_\zeta\ra + \la g, \tau_\zeta\ra\right)\right]\mathrm{d} \prob^{MM}_{\b, \Lda}(\varphi)=\e^{\iota\la g, \tau_\zeta\ra} \int \exp\left[\iota \la \vr, \tau_\zeta\ra\right]\mathrm{d} \prob^{MM}_{\b, \Lda}(\varphi)\\
    &\qquad\qquad =\e^{\iota\la g, \tau_\zeta\ra} \e^{-E_{\b}(a_\zeta, \tau_\zeta)} \int \exp\left[\iota \la \vr, \tau_\zeta - \b\Delta^2a_\zeta\ra\right]\mathrm{d} \prob^{MM}_{\b, \Lda}(\varphi).
\end{align*}
But 
\begin{align*}
  E_{\b}(a_\zeta, \tau_\zeta)& = \la a_\zeta, \tau_\zeta\ra -\frac{\b}2 \la a_\zeta, \Delta^2a_\zeta\ra \\
  &= \sum_{\r\in\cN}\sum_{\r^{'}\in\cN} \zeta(\r)\zeta(\r^{'})\la a_\r, \r^{'}\ra - \frac{\b}2\sum_{\r\in\cN}\sum_{\r^{'}\in\cN} \zeta(\r)\zeta(\r^{'})\la a_\r, \Delta^2a_{\r^{'}}\ra.
\end{align*}

By property~\ref{item:prop4_2.4} of Proposition~\ref{prop:2.4} 
\begin{align*}
  \la a_\r, \r^{'}\ra=0,\quad  \la a_\r, \Delta^2a_{\r^{'}}\ra=0, \qquad \r,\r^{'}\in\cN, \r\neq\r^{'}. 
\end{align*}
Therefore
\begin{align*}
   E_{\b}(a_\zeta, \tau_\zeta)=\sum_{\r\in\cN} |\zeta(\r)|E_{\b}(a_\r, \r).
\end{align*}

Thus we have 
\begin{align*}
    &\int \prod_{\r\in\cN} \left[ 1+ K(\r)\cos\left(\la\vr,\r\ra + \la g,\r\ra\right)\right]\mathrm{d} \prob^{MM}_{\b, \Lda}(\varphi)\\
    &= \sum_{\zeta:\cN\to\{-1,0,1\}} \left[\prod_{\r\in\cN}\left(\frac{K(\r)}2 \e^{-E_{\b}(a_\r, \r)}\right)^{|\zeta(\r)|}\right] \int \e^{\iota \la \vr, \sum_{\r\in\cN}\zeta(\r)\r - \sum_{\r\in\cN}\zeta(\r)\b\Delta^2a_\r\ra + \iota\la g, \sum_{\r\in\cN}\zeta(\r)\r\ra}\mathrm{d} \prob^{MM}_{\b, \Lda}(\varphi)\\
    &=\sum_{\zeta:\cN\to\{-1,0,1\}} \int \prod_{\r\in\cN}\left(\frac{K(\r)}2 \e^{-E_{\b}(a_\r, \r)}\right)^{|\zeta(\r)|} \e^{\iota \zeta(\r)\left(\la \vr, \r -\b\Delta^2a_\r\ra + \la g, \r\ra\right)}\mathrm{d} \prob^{MM}_{\b, \Lda}(\varphi)\\
    &= \int\prod_{\r\in\cN}\left[1+z(\b, \r,\cN)\cos\left( \la \vr, \r - \b\Delta^2a_{\r,\cN}\ra + \la g, \r\ra \right)\right] \mathrm{d} \prob^{MM}_{\b, \Lda}(\varphi),
\end{align*}
where $z(\b, \r,\cN)= K(\r)\e^{-E_{\b}(a_\r, \r)}$. By property~\ref{item:prop5_2.4} of Proposition~\ref{prop:2.4} we obtain
\begin{align*}
    |z(\b, \r,\cN)| \le |K(\r)| \exp\left[-\frac1\b \left( \frac1{144}\|\r\|_2^2 + D_3 \sum_{k=8}^{n(\r)} |\cS^{sep}_k(\r)|\right)\right].
\end{align*}

\paragraph{Step 3: bounding the coefficients $z(\b,\r,\cN)$.}
In this step we finish the proof by showing \eqref{eq:t2.3.1}. Suppose the coefficients $(K(\r))_{\r\in\cN}$ satisfy condition~\ref{item:three_6.4} of Theorem~\ref{thm:2.2}, that is,
\begin{align*}
    |K(\r)| \le e^{D_2A(\r)} \prod_{x\in \supp(\r)}e^{\r(x)^2} \left|1-\frac{|\r(x)|}{N}\right|.
\end{align*}

 We have
 \begin{align*}
     \prod_{x\in \supp(\r)}e^{\r(x)^2}\big|1-\frac{|\r(x)|}{N}\big|
     &\le \prod_{x\in \supp(\r)}e^{\r(x)^2} e^{\frac{|\r(x)|}{N}} \le \prod_{x\in \supp(\r)}e^{2\r(x)^2}  \le \e^{2\|\r\|_2^2}.
 \end{align*}
 Put $\kappa:=\min\{1/288, D_3\}$. Then
 $$\frac{1}{144}\|\rho\|_2^2+ D_3I(\rho)\ge \frac{1}{288}\|\rho\|_2^2+\kappa( \|\rho\|_2^2+ I(\rho)).$$
Since $\rho$ is integer-valued, $|S_0(\rho)|=|\mathrm{supp}\, \rho|\le \|\rho\|_2^2$. Proposition~\ref{prop:2.1} therefore gives
$$\log_2( d_{\Lambda}(\rho)+1)\le A(\rho)\le D_1(\|\rho\|_2^2+ I(\rho)).$$
Inserting both bounds into \eqref{eq:t2.3.2} and using $|K(\rho)|\le \exp\left({D_2 A(\rho)+ 2\|\rho\|_2^2}\right)$,
\[
|z(\beta,\rho,\cN)|\le\exp\Big(-\frac1\beta\Big[\Big(\frac{1}{288}-2\beta\Big)\|\rho\|_2^2
+\big(\kappa-\beta D_1D_2\big)\big(\|\rho\|_2^2+I(\rho)\big)\Big]\Big).
\]
Let $\beta_0:=\min\{1/1152,\ \kappa/(2D_1D_2)\}$. For $0<\beta<\beta_0$, both coefficients are positive: the first is at least $1/576$ and the second at least $\kappa/2$. Using once more $\|\rho\|_2^2+I(\rho)\ge \log_2(d_{\Lambda}(\rho)+1)/D_1$,
\[
|z(\beta,\rho,\cN)|\le\exp\Big(-\frac{C_2}{\beta}\Big[\|\rho\|_2^2+\log_2(d_{\Lambda}(\rho)+1)\Big]\Big),
\qquad C_2:=\min\Big\{\frac1{576},\frac{\kappa}{2D_1}\Big\}.
\]
Hence \eqref{eq:t2.3.1} follows. Note that the numerical coefficients are simply safety margins to ensure that the negative exponential coming from the spin-wave energy dominates the positive combinatorial factors in \(K(\rho)\) in~\eqref{eq:t2.3.2}.

\subsection{Spin-wave construction: proof of Proposition~\ref{prop:2.4}}\label{sec:prop 2.4}
We fix an ensemble $\cN$ satisfying properties~\ref{item:one_6.4}--\ref{item:two_6.4} of Theorem~\ref{thm:2.2}. Also, we fix $\r\in\cN$ for the rest of this proof. The spin wave $a_{\r}$ will be constructed as a sum of spin waves at different scales, which we define below.

\subsubsection{Initial spin wave}

The following lemma is analogous to Lemma 4.5 of \cite{K-P_2017} and Lemma 25 of \cite{Wirth_2019}.
\begin{lemma}\label{lem:4.5}
There exists $a_{0,\r}:\Z^2\to\R$ such that
\begin{enumerate}[label=(\roman*),ref=(\roman*)]
    \item \label{item:roman1}$\supp(a_{0,\r})\subseteq\supp(\r)$.
   \item\label{item:roman4} $E_{\b}(a_{0,\r},\r) \ge \|\r\|_2^2/({144\b})$.
\end{enumerate}
\end{lemma}
\begin{proof}[Proof of Lemma~\ref{lem:4.5}]
  We first define $a_{0,\r}$. Note that $\Z^2$ is bipartite. Let $(\Omega_1, \Omega_2)$ be a bipartition of $\Z^2$, chosen without loss of generality so that
  \begin{align}\label{eq:l4.5.1}
      \sum_{x\in\Omega_1} \r(x)^2 \ge \frac12 \sum_{x\in \Z^2} \r(x)^2.
  \end{align}

  Let 
  \begin{align*}
      \r_1(x)=\begin{cases}
          \r(x) & \text{ if }x\in\Omega_1\\
          0&  \text{otherwise}.
      \end{cases}
  \end{align*}
  Note that \eqref{eq:l4.5.1} implies $\|\r_1\|_2^2 \ge 1/2 \|\r\|_2^2$. Define $a_{0,\r}:\Z^2\to\R$ by 
  \begin{align}\label{eq:def_a0rho}
      a_{0,\r}(x):=\frac{\r_1(x)}{24 \b}.
  \end{align}
Next, we verify the properties of the lemma. Property~\ref{item:roman1} holds by definition. 
We now verify property~\ref{item:roman4}. We have
\begin{align*}
    E_{\b}(a_{0,\r},\r) &= \la a_{0,\r},\r\ra  -\frac{\b}2 \la a_{0,\r}, \Delta^2a_{0,\r} \ra= \sum_{x\in\Omega_1}\frac{(\r(x))^2}{24\b} - \frac{\b}2\sum_{x\in\Z^2} \left(\Delta a_{0,\r}(x)\right)^2 .
\end{align*}
Now observe that for any $x\in\Omega_1$
\begin{align*}
 \left(\Delta a_{0,\r}(x)\right)^2 = \left(\sum_{y\sim x} \left(a_{0,\r}(y)-a_{0,\r}(x)\right)\right)^2 = \left(4 \frac{\r(x)}{24\b}\right)^2= \frac{\r(x)^2}{36\b^2}  .
\end{align*}
Also, for $x\notin \Omega_1$,
\begin{align*}
 \left(\Delta a_{0,\r}(x)\right)^2 = \left(\sum_{y\sim x} \left(a_{0,\r}(y)-a_{0,\r}(x)\right)\right)^2 = \left(\sum_{y\sim x} \frac{\r(y)}{24\b}\right)^2\le 4 \sum_{y\sim x, y\in\Omega_1}\frac{\r(y)^2}{576\b^2} . 
\end{align*}
From the above two observations we conclude that
\begin{align*}
   \sum_{x\in\Z^2} \left(\Delta a_{0,\r}(x)\right)^2 \le \sum_{x\in\Omega_1}\frac{\r(x)^2}{18\b^2}.  
\end{align*}
Thus 
\begin{align*}
   E_{\b}(a_{0,\r},\r) &\ge \sum_{x\in\Omega_1}\frac{(\r(x))^2}{24\b} - \frac{\b}{2} \sum_{x\in\Omega_1}\frac{\r(x)^2}{18\b^2}\\ 
   &= \left(\frac{1}{24\b}-\frac{1}{36\b}\right)\sum_{x\in\Omega_1}(\r(x))^2\ge \frac{1}{144\b}\|\r\|_2^2.\qedhere
\end{align*}
\end{proof}

\subsubsection{Affine modification on separated squares}
\label{d2:sec:affine}
For any function $u$ and a finite subset $A\subset \Z^2$, we define 
\begin{align*}
&\|\HH u\|_{\infty, A}:= \sup\{|\HH u(x)|: x\in A\},\\
&\|\HH u\|_{\infty}:= \sup\{|\HH u(x)|: x\in \Z^2\}.
\end{align*}

The following local modification is the two-dimensional version of the
construction in~\citet[Lemma 4.3]{schweiger2022pinning}. Schweiger's statement produces, for each function $v$, a local correction $w$. The proof constructs this correction in the form
\[
    w=\xi_{x,r,R}(A_xv-v),
\]
where the cut-off $\xi_{x,r,R}$ depends only on $x,r,R$ and $A_xv$
is the affine Taylor approximation of $v$ at $x$.  Consequently, the
map $v\mapsto v+w$ is linear.  We record this operator form, together
with the fixed-point property \eqref{d2:eq:affinekernel}, because both features
are needed when the modifications are composed in Lemma~\ref{d2:lem:sparse}. 

\begin{lemma}\label{d2:lem:localmodification}
There is $c<\infty$ such that, for integers $R\ge16r\ge16$
and $x\in\Z^2$, there is a linear operator
$\mathcal T_{x,r,R}$ satisfying
\begin{align}
 \supp(\mathcal T_{x,r,R}v-v)&\subset\bx_{R-1}(x),\label{d2:eq:localmodsupport.1}\\
 \supp(\HH(\mathcal T_{x,r,R}v-v))&\subset\bx_R(x),
       \label{d2:eq:localmodsupport}\\
 \mathcal T_{x,r,R}v&\text{ is affine on }\bx_r(x),
       \label{d2:eq:localmodaffine}\\
 \|\HH\mathcal T_{x,r,R}v\|_{\infty,\bx_R(x)}
 &\le\left(1+\frac{c}{\log(R/r)}\right)
             \|\HH v\|_{\infty,\bx_R(x)}.
       \label{d2:eq:localmodbound}
\end{align}
It can be chosen in the form
\begin{equation}\label{d2:eq:localoperatorexplicit}
 \mathcal T_{x,r,R}v
 =v+\xi_{x,r,R}(\mathcal A_xv-v),\qquad
 \mathcal A_xv(y)=v(x)+\sum_{j=1}^2(y_j-x_j)D_j^+v(x),
\end{equation}
where $\xi_{x,r,R}=1$ on $\bx_{r+1}(x)$ and
$\supp(\xi_{x,r,R})\subset\bx_{R-1}(x)$. In particular,
\begin{equation}\label{d2:eq:affinekernel}
 \HH v=0\text{ on }\bx_R(x)\quad\Longrightarrow\quad
 \mathcal T_{x,r,R}v=v.
\end{equation}
\end{lemma}
\begin{proof}
We first prove that there exists a function \(\xi:\mathbb{Z}^2\to\mathbb{R}\) such that
\[
\xi\equiv 1 \quad\text{on }\bx_{r+1}(0),\qquad
\xi\equiv 0 \quad\text{on }\mathbb{Z}^2\setminus \bx_{R-1}(0),\qquad
0\le \xi \leq 1,
\]
and
\begin{align}
|D_i^+\xi(y)|&\leq
\frac{C}{(1+\|y\|_2)\log(R/r)},\label{eq:5.13.6}\\
|D_i^+ D_j^+\xi(y)|&\leq
\frac{C}{(1+\|y\|_2)^2\log(R/r)},\label{eq:5.13.7}
\end{align}
for all $i,j\in\{1,2\}$ and all \(y\in \Z^2\).

To prove this, we choose a function $\phi \in C^\infty(\mathbb{R}^2)$ such that $0\le \phi \le 1$, 
$\phi\equiv 1$ on $[-1,1]^2$, $\phi\equiv 0$ on $\R^2\setminus [-2,2]^2$. Then for any $\delta>0$, we define
\[
\phi_\delta(y) = \phi\left(\frac{y}{\delta}\right).
\]
We define $\widetilde{\xi}: \mathbb{R}^2 \to \mathbb{R}$ by
\[
\widetilde{\xi}(y)
= \phi_{2r}(y)
+
\bigl(\phi_{R/4}(y)-\phi_{2r}(y)\bigr)
\frac{\log R-\log\|y\|_2}{\log R-\log r}.
\]
One can check that $0 \leq \widetilde{\xi} \leq 1$,
$\widetilde{\xi}\equiv 1$ on $[-2r,2r]^2 \supset \bx_{r+1}(0)$,
$\widetilde{\xi} \equiv 0$ on
$\mathbb{R}^2 \setminus [-R/2,R/2]^2
\supset \mathbb{Z}^2 \setminus \bx_{R-1}(0)$ and also
\[
|\partial^\g \widetilde{\xi}(y)|
\leq
\frac{C_{\g}}{(1+\|y\|_2)^{|\g|}(\log R-\log r)},
\]
for a multi-index $\g=(\g_1,\g_2)$ with $|\g|=\g_1+\g_2$ and $1\le |\g| \le 2$, where $\partial^\g=\frac{\partial^{\g_1}}{\partial y_1} \frac{\partial^{\g_2}}{\partial y_2}$. We now define $\xi$ to be the restriction of
$\widetilde{\xi}$ to $\mathbb{Z}^2$. The above estimates on
$\widetilde{\xi}$ together with Taylor's theorem then imply \eqref{eq:5.13.6} and \eqref{eq:5.13.7}.

For \(v:\mathbb{Z}^2\to\mathbb{R}\), define
\[
\mathcal A_x v(y)
:=v(x)+\sum_{j=1}^{2}(y_j-x_j)D_j^+v(x)
=v(x)+(y-x)\cdot\nabla v(x),
\]
and 
\[
\mathcal T_{x,r,R}v=v+\xi_{x,r,R}(\mathcal A_x v-v),
\]
where
\[
\xi_{x,r,R}(y)=\xi(y-x).
\]
Then
\[
\mathcal T_{x,r,R}v -v = \xi_{x,r,R}(\mathcal A_x v-v).
\]
Therefore
\[
\supp(\mathcal T_{x,r,R}v -v) \subset \bx_{R-1}(x), 
\qquad \supp(\HH(\mathcal T_{x,r,R}v -v)) \subset \bx_{R}(x).
\]
Since $\xi_{x,r,R} \equiv 1$ on $\bx_{r+1}(x)$, we observe that $\mathcal T_{x,r,R}v$ is affine on $\bx_{r}(x)$.

We now prove \eqref{d2:eq:localmodbound}. Fix $y\in\bx_R(x)$. We have
\begin{align*}
\HH \mathcal T_{x,r,R}v(y)=\HH v(y)+ \HH \left(\xi_{x,r,R}(\mathcal A_x v-v)(y) \right).
\end{align*}
Let us write $\xi_x=\xi_{x,r,R}$. Then we obtain 
\begin{align*}
 \HH \left(\xi_{x}(\mathcal A_x v-v) \right)_{ij}(y)
&=-\xi_x(y+e_i-e_j)D_i^+D_j^- v(y) +D_i^+\xi_x(y-e_j)
D_j^-\left(\mathcal A_x v- v\right)(y)\\
& +D_j^-\xi_x(y+e_i)
D_i^+\left(\mathcal A_x v- v\right)(y) +\left(\mathcal A_x v- v\right)(y)
D_i^+D_j^-\xi_x(y).
\end{align*}
Hence,
\begin{align}\label{eq:5.13.5}
\left|
\left(\HH\mathcal T_{x,r,R}v\right)_{ij}(y)
\right|
&\leq
\left|(\HH v)_{ij}(y)\right|
+\left|D_i^+\xi_x(y-e_j)\right|
\left|D_j^-\left(\mathcal A_x v- v\right)(y)\right|\notag\\
&\quad
+\left|D_j^-\xi_x(y+e_i)\right|
\left|D_i^+\left(\mathcal A_x v- v\right)(y)\right|\notag\\
&\quad
+\left|\left(\mathcal A_x v- v\right)(y)\right|
\left|(\HH\xi_x)_{ij}(y)\right|.
\end{align}

Note that by \eqref{eq:5.13.6} and \eqref{eq:5.13.7} we get for any $i,j$
\begin{align}\label{eq:5.13.3}
\left|D_i^+\xi_x(y-e_j)\right|, \left|D_j^-\xi_x(y+e_i)\right| \le \frac{C}{(1+|y-x|)\log(R/r)},
\end{align}
and 
\begin{align}\label{eq:5.13.4}
\left|(\HH\xi_x(y))_{ij}\right| \le  \frac{C}{(1+|y-x|)^2\log(R/r)}.
\end{align}
We now show that for all $z\in\bx_R(x)$

\begin{equation}
\left|\mathcal A_x  v(z)-  v(z)\right|
\leq C |z-x|^2\,\|\HH v\|_{\infty,\bx_R(x)}.
 \label{eq:5.13.2}
\end{equation}

First, we show that for any $i$ and any $z\in\bx_R(x)$,
\begin{equation}
\left|D_i^+  v(z)-D_i^+  v(x)\right|
\leq C |z-x|\,\|\HH v\|_{\infty,\bx_R(x)}. \label{eq:5.13.1}
\end{equation}

Let $z\in\bx_R(x)$ and \(x=a_0,a_1,\ldots,a_n=z\) be a shortest path in $\bx_R(x)$,
where $a_{k+1}=a_k\pm e_j,\quad j\in\{1,2\}.$
Then
\[
D_i^+ v(z)-D_i^+ v(x)
=
\sum_{k=0}^{n-1}
\left(
D_i^+  v(a_{k+1})
-
D_i^+  v(a_k)
\right).
\]
If \(a_{k+1}=a_k+e_j\), then \[\left|
D_i^+  v(a_{k+1})
-D_i^+  v(a_k)
\right|
=\left|(\HH v)_{ij}(a_{k+1})\right|
\leq \|\HH v\|_{\infty,\bx_R(x)}.
\]
If \(a_{k+1}=a_k-e_j\), then
\begin{align*}
\left|
D_i^+  v(a_{k+1})
-D_i^+  v(a_k)
\right|\
=\left|- (\HH v)_{ij}(a_k)\right|
\leq \|\HH v\|_{\infty,\bx_R(x)}.
\end{align*}
Thus,
\begin{align*}
\left|D_i^+  v(z)-D_i^+  v(x)\right|
\leq |z-x|_1\,\|\HH v\|_{\infty,\bx_R(x)}
\leq C|z-x|\,\|\HH v\|_{\infty,\bx_R(x)}.
\end{align*}

We now prove \eqref{eq:5.13.2}. Let $z\in\bx_R(x)$ and set
\begin{align*}
\widetilde{\mathcal A}_x  v(z)= v(z)-\mathcal A_x  v(z)
= v(z)-v(x)-\sum_{j=1}^{2}(z_j-x_j)D_j^+  v(x).
\end{align*}
Choose a shortest path \(x=a_0,a_1,\ldots,a_n=z\) in $\bx_R(x)$. Let
\[
a_{k+1}=a_k+\varepsilon_k e_{j_k},
\qquad \varepsilon_k\in\{-1,1\}, j_k\in\{1,2\}.
\]
Set $w_k=a_k$ when $\varepsilon_k=1$ and $w_k=a_{k+1}$ otherwise.
Then
\[
  v(z)-  v(x)
=
\sum_{k=0}^{n-1}
\left(  v(a_{k+1})-  v(a_k)\right) = \sum_{k=0}^{n-1}\varepsilon_kD_{j_k}^+  v(w_k).
\]
Therefore,
\begin{align*}
\widetilde{\mathcal A}_x  v(z)
&=
\sum_{k=0}^{n-1}
\varepsilon_kD_{j_k}^+  v(w_k)
-\sum_{j=1}^{2}(z_j-x_j)D_j^+  v(x)\\
&=
\sum_{k=0}^{n-1}\varepsilon_k
\left(
D_{j_k}^+  v(w_k)
-D_{j_k}^+  v(x)
\right).
\end{align*}
Now, using \eqref{eq:5.13.1}, we get
\begin{align*}
\left|\widetilde{\mathcal A}_x  v(z)\right|
&\leq
\sum_{k=0}^{n-1}
\left|
D_{j_k}^+  v(w_k)
-D_{j_k}^+  v(x)
\right|\\
&\leq
\sum_{k=0}^{n-1}
|w_k-x|_1\,\|\HH v\|_{\infty,\bx_R(x)}\\
&\leq n^2\|\HH v\|_{\infty,\bx_R(x)} \le C |z-x|^2\|\HH v\|_{\infty,\bx_R(x)}.
\end{align*}
Using \eqref{eq:5.13.5}, \eqref{eq:5.13.3}, \eqref{eq:5.13.4}, \eqref{eq:5.13.1} and \eqref{eq:5.13.2}, we have the desired bound \eqref{d2:eq:localmodbound}.

Finally, suppose
\[
\HH v = 0 \quad \text{ on } \bx_R(x).
\]
Let \(g_i = D_i^+ v\). Then
\[
g_i(y+e_j)-g_i(y)=(Hv)_{ij}(y+e_j).
\]
Thus \(g_i\) is constant on \(\bx_{R-1}(x)\).
Hence $g_i(y)=D_i^+ v(x)$ for all $y\in\bx_{R-1}(x)$.

Let \(y\in\bx_{R-1}(x)\) and let
\(x=z_0,z_1,\ldots,z_m=y\) be a path in $\bx_{R-1}(x)$. Then
\begin{align*}
v(y)-v(x)
&=\sum_{k=0}^{m-1} [v(z_{k+1})-v(z_k)]\\
&=\sum_{j=1}^{2}(y_j-x_j)D_j^+ v(x),
\end{align*}
and hence $\mathcal T_{x,r,R}v(y)=v(y)$. Since $\supp(\xi_{x,r,R})\subset\bx_{R-1}(x)$, we conclude that $\mathcal T_{x,r,R}v=v$.
\end{proof}

\begin{lemma}\label{lem:4.6}\label{d2:lem:components}
Let $(\bx_{R_i}(x_i))_{i=1}^m$ be a finite family, $R_i\ge1$,
such that, for $i\ne j$,
\[
 |x_i-x_j|_\infty\ge\kappa\min(R_i,R_j)^q,\qquad
 q>1,\qquad\kappa\ge4\cdot10^q.
\]
Every connected component $\mathcal C$ of $\,\bigcup_{i=1}^{m} \bx_{R_i}(x_i)$ satisfies
\[
 \diam\mathcal C\le10\max_{\bx_{R_i}(x_i)\subset\mathcal C}R_i.
\]
Connectedness here refers to the intersection graph of the squares.
\end{lemma}

\begin{proof}
We prove this by induction on $m$. The lemma is immediate for $m=1$.
Suppose that it holds for $1\le m\le k$. Let $\mathcal C$ be a connected component of $\bigcup_{i=1}^{k+1} \bx_{R_i}(x_i)$.
Let
\begin{align*}
R=\max\{R_i: \bx_{R_i}(x_i)\subseteq \mathcal C\},
\end{align*}
and
\begin{align*}
\mathcal{M}=\{i: \bx_{R_i}(x_i)\subseteq \mathcal C \text{ and } R_i\ge R/10\}.
\end{align*}

We first show that $|\mathcal{M}|=1$, which implies that there is a unique $R_{i_*}$ such that $R_{i_*}=R$ and $\bx_{R_{i_*}}(x_{i_*})\subseteq \mathcal C$.
Assume that there are two distinct $i_1,i_2\in\mathcal{M}$. Suppose $\bx_{R_{i_1}}(x_{i_1})\cap \bx_{R_{i_2}}(x_{i_2})\neq\emptyset$.
Then
\begin{align*}
|x_{i_1}-x_{i_2}|_\infty\le R_{i_1}+R_{i_2}\le 2R.
\end{align*}
But
\begin{align*}
|x_{i_1}-x_{i_2}|_\infty
&\ge \kappa\min(R_{i_1},R_{i_2})^q\\
&\ge \kappa\left(\frac{R}{10}\right)^q\\
&\ge 4R^q\ge 4R,
\end{align*}
which is a contradiction. Therefore, for any distinct $i_1,i_2\in\mathcal{M}$, we must have $\bx_{R_{i_1}}(x_{i_1})\cap \bx_{R_{i_2}}(x_{i_2})=\emptyset$.

Since $\mathcal C$ is connected, there exist $j_1,j_2\in\mathcal{M}$ which intersect some connected component $\mathcal C'$ of $\bigcup_{\substack{i=1\\i\notin\mathcal{M}}}^{k+1}\bx_{R_i}(x_i)$. Then by the induction hypothesis,
\begin{align*}
\kappa\left(\frac{R}{10}\right)^q
&\le |x_{j_1}-x_{j_2}|_\infty\\
&\le R_{j_1}+\diam(\mathcal C')+R_{j_2}\\
&\le 2R+10\cdot\max_{\substack{\bx_{R_i}(x_i)\subseteq \mathcal C'}}R_i\\
&\le 2R+10\cdot\frac{R}{10},
\end{align*}
which is not possible as $\kappa\ge 4\cdot 10^q $. Therefore $|\mathcal{M}|=1$.

Let $\mathcal{M}=\{i_*\}$. After removing the unique square $Q_R(x_{i_*})$, every
remaining component of $\mathcal C$ has diameter at most $R$ by induction and intersects $Q_R(x_{i_*})$, since $\mathcal C$ is connected. Hence
$\mathcal C\subset Q_{2R}(x_{i_*})$ and
$\diam\mathcal C\le4R\le10R$, and the proof is complete by induction.
\end{proof}

\begin{lemma}\label{d2:lem:sparse}
Fix $a=100$ and $T = 128a$. Let $(x_i,r_i)_{i=1}^m$ be a finite family
with $x_i\in\Z^2$, $r_i\ge1$ integers, and
\begin{equation}\label{d2:eq:centreseparation}
 |x_i-x_j|_\infty\ge M\min(r_i,r_j)^\a\qquad(i\ne j).
\end{equation}
Order the indices so that $r_1\ge\cdots\ge r_m$, and set
\[
 s_i=4ar_i,\qquad R_i=\lceil Tr_i^{3/2}\rceil,\qquad D_i=\bx_{R_i}(x_i).
\]
For $v$ with bounded Hessian, define
\begin{equation}\label{d2:eq:finiteiteration}
 v^{(0)}=v,\qquad
 v^{(i)}=\mathcal T_{x_i,s_i,R_i}v^{(i-1)},\quad1\le i\le m,
 \qquad \mathcal T v=v^{(m)}.
\end{equation}
Then:
\begin{enumerate}[label=(\roman*),leftmargin=2.2em]
\item \label{item:lem_4.15.1} $\mathcal T v$ is affine on $\bx_{ar_i}(x_i)$ for every $i$.
\item \label{item:lem_4.15.2} For $ C_H:=\exp\left[\frac{c}{\log(T/(4a))(1-3/(2\a))}\right]$
\[
 \|\HH\mathcal T v\|_\infty\le C_H\|\HH v\|_\infty.
\]
\item \label{item:lem_4.15.3} If $K=\supp(\HH v)$ and $\mathcal D(K)$ is the union of
components of $\bigcup_iD_i$ meeting $K$, then
\begin{equation}\label{d2:eq:localisationprecise}
 \supp(\mathcal T v-v)\subset\mathcal D(K),\qquad
 \supp(\HH\mathcal T v)\subset K\cup\mathcal D(K).
\end{equation}
\end{enumerate}
For $m=0$, take $\mathcal T=\mathrm{id}$.
\end{lemma}
\begin{proof}
First, note that $R_i=\lceil T r_i^{3/2}\rceil \geq 32 s_i$. Hence the operator $\mathcal T_{x_i,s_i,R_i}$ is defined.

Since \(Tr_i^{3/2}\leq R_i\leq 2Tr_i^{3/2}\), we have for any \(i\ne j\),
\begin{align*}
 |x_i-x_j|_\infty
 &\geq M\min(r_i,r_j)^\a\\
 &\geq M\min\left(\left(\frac{R_i}{2T}\right)^{2/3},
                   \left(\frac{R_j}{2T}\right)^{2/3}\right)^\a\\
 &=\frac{M}{(2T)^{2\a/3}}\min(R_i,R_j)^{2\a/3}.
\end{align*}
Since $ M\geq 4\cdot 10^{2\a/3}(2T)^{2\a/3},$ we can apply Lemma~\ref{d2:lem:components}. By Lemma~\ref{d2:lem:components}, every connected component \(\mathcal C\) of the union, including connected components of unions of subfamilies, satisfies
\begin{align}\label{eq:4.15.1}
 \diam\mathcal C
 \leq 10\max_{\bx_{R_i}(x_i)\subseteq\mathcal C} R_i
 \leq 20 T r_*^{3/2},
 \qquad
 r_*:=\max_{\bx_{R_i}(x_i)\subseteq\mathcal C}r_i.
\end{align}

We first prove \ref{item:lem_4.15.1}. By Lemma~\ref{d2:lem:localmodification}, the function \(v^{(i)}\) is affine on
\(\bx_{4 a r_i}(x_i)\). Let \(\mathcal C\) be a connected component of
\(\bigcup_{j>i} \bx_{R_j}(x_j)\) such that $\mathcal C\cap \bx_{a r_i}(x_i)\neq\emptyset.$
Let
\[
 r_*:=\max_{\substack{\bx_{R_j}(x_j)\subseteq\mathcal C\\j>i}}r_j.
\]
Then 
\[
 |x_i-x_*|\geq M r_*^\a,
\]
and, by \eqref{eq:4.15.1}
\[
 |x_i-x_*|\leq a r_i+20T r_*^{3/2}.
\]
Consequently,
\begin{align*}
 M r_*^\a\leq  a r_i+20T r_*^{3/2}
  \Longleftrightarrow  Mr_*^\a\left(1-\frac{20T}{M}r_*^{3/2-\a}\right)\leq  a r_i.
\end{align*}
Note that by \eqref{eq:def_alpha_M}, we have \((20T/M)r_*^{3/2-\a}<1/2\). Therefore
\[
 r_*^\a\leq \frac{2 a r_i}{M}.
\]

Using \eqref{eq:4.15.1} again, we obtain
\begin{align*}
 \diam(\mathcal C)
 &\leq 20T r_*^{3/2}
% =20T(r_*^\a)^{p/\a}\\
 \leq 20T\left(\frac{2 a}{M}\right)^{3/(2\a)}r_i^{3/(2\a)}\\
 &\leq 20T\left(\frac{2 a}{M}\right)^{3/(2\a)}
       r_i \le a r_i.
\end{align*}
In the above, we used \eqref{eq:def_alpha_M} to observe that 
$
 M\geq \left(\frac{20T}{a}\right)^{2\a/3}2a.
$
Now, since \(\mathcal C\cap \bx_{a r_i}(x_i)\neq\emptyset\) and
\(\diam(\mathcal C)\leq a r_i\), it follows that
\[
 \mathcal C\subseteq \bx_{2a r_i}(x_i)
 \subseteq \bx_{4a r_i}(x_i).
\]
Since \(v^{(i)}\) is affine on \(\bx_{4a r_i}(x_i)\), it is an affine function on \(\mathcal C\). Every later square in this component lies in a region on which $v^{(i)}$ is affine. By the explicit operator formula \eqref{d2:eq:localoperatorexplicit}, its first correction is zero.
Inductively, every correction in this component is zero. Operations in other later components do not change values in \(\mathcal C\). All other components avoid $\bx_{a r_i}(x_i)$. Thus
\[\mathcal Tv = v^{(m)}=v^{(i)} \text{ on } \bx_{a r_i}(x_i).\]

Next we prove \ref{item:lem_4.15.2}. Define
\[
 h_i:=\log\left(\frac{R_i}{s_i}\right) =\log(\lceil Tr_i^{3/2}\rceil)-\log(4 a r_i)>0.
\]

Let \(r_i\leq r_j\) and suppose \(\bx_{R_i}(x_i)\cap \bx_{R_j}(x_j)\neq\varnothing\). For \(x\in \bx_{R_i}(x_i)\cap \bx_{R_j}(x_j)\),
\begin{align*}
 M r_i^\a &\leq |x_i-x_j|_\infty \leq |x_i-x|_\infty+|x-x_j|_\infty\\
 &\leq R_i+R_j \leq 2Tr_i^{3/2}+2Tr_j^{3/2} \leq 4Tr_j^{3/2}.
\end{align*}
Thus
\begin{align}\label{eq:5.15.2}
%r_i^\a\leq \frac{4T}{M}r_j^{3/2}
 \log r_i\leq \frac{1}{\a}\log\left(\frac{4T}{M}\right)
                +\frac{3}{2\a}\log r_j.
\end{align}

Then
\[
\log\left(\frac{T}{4 a}\right) + \frac12 \log r_i \le  h_i \le \log 2+\log\left(\frac{T}{4 a}\right) + \frac12 \log r_i.
\]
Now, using \eqref{eq:5.15.2}, we obtain
\begin{align}
 h_i-\frac{3}{2\a}h_j
 &\leq \log 2+\left(1-\frac{3}{2\a}\right)
       \log\left(\frac{T}{4a}\right)
 +\frac12\left(\log r_i-\frac{3}{2\a}\log r_j\right)\notag\\
 &\leq \log 2+\left(1-\frac{3}{2\a}\right)
       \log\left(\frac{T}{4a}\right)
       +\frac{1}{2\a}\log\left(\frac{4T}{M}\right)\notag\\
 &\leq 0.\label{eq:5.15.3}
\end{align}
We obtain the last inequality by using \eqref{eq:def_alpha_M} and observing that
\[
 M\ge4T\exp\!\left\{2\alpha
     \left[\log2+(1-\frac3{2\a})\log\frac{T}{4a}\right]\right\}.
\]
%\[
% M\geq 4T 2^{
% 2\a\left(1-\frac{3}{2\a}\right)
% \log\left(\frac{T}{4a}\right)+\log 2}.
%\]
Define
\[
 c_i=1+\frac{c}{h_i},
\]
where \(c\) is the constant in Lemma~\ref{d2:lem:localmodification}. Then \(c_i>1\) for every \(i\).

Define recursively
\[
 B_i=c_i
      \max\bigl(\{1\}\cup
         \{B_j:j<i,\ D_j\cap D_i\ne\varnothing\}\bigr),
\]
and
\[
 F_i(x)=\max\bigl(\{1\}\cup\{B_j:j\le i,\ x\in D_j\}\bigr),
 \qquad F_0(x)=1.
\]
We claim, by induction, that
\[
 |\HH v^{(i)}(x)|\le F_i(x) \|\HH v\|_\infty.
\]
Outside $D_i$, the Hessian is unchanged. On $D_i$, \eqref{d2:eq:localmodbound} and \eqref{d2:eq:localmodsupport} give
\begin{align*}
 \|\HH v^{(i)}\|_{\infty,D_i}
 &\le c_i
        \|\HH v^{(i-1)}\|_{\infty,D_i}\\
 &\le B_i\|\HH v\|_\infty.
\end{align*}
The second inequality uses only earlier squares meeting $D_i$, which proves
the claim.

Expanding the finite recursion for $B_i$ gives a maximum of products
\[
 \prod_{\ell=1}^{t}c_{i_\ell},
 \quad i_1<\cdots<i_t=i,\quad
 D_{i_\ell}\cap D_{i_{\ell+1}}\ne\varnothing.
\]
Along such a chain, we have by \eqref{eq:5.15.3}
$h_{i_{\ell+1}}\le (3/(2\a)) h_{i_\ell}$. Therefore
\[
 \sum_{\ell=1}^t\frac1{h_{i_\ell}}
 \le\frac1{h_{i_t}}\sum_{j=0}^{t-1}\left(\frac3{2\a}\right)^j
 \le\frac{1}{\log(T/(4a))(1-3/(2\a))}.
\]
Using $\log(1+u)\le u$ yields
\[
 B_i\le\exp\!\left(\frac c{\log(T/(4a))(1-3/(2\a))}\right)=C_H,
 \]
and hence
\[
\|\HH\mathcal T v\|_\infty\le C_H\|\HH v\|_\infty.
\]
This completes the proof of \ref{item:lem_4.15.2}.

Finally, we prove \ref{item:lem_4.15.3}. Let $K=\supp(\HH v)$
and let \(\mathcal D(K)\) be the union of the connected components of \(\bigcup \bx_{R_i}(x_i)\) that meet \(K\). Let \(\mathcal C\) be a component of \(\bigcup \bx_{R_i}(x_i)\) such that
\(\mathcal C\cap K=\varnothing\). Then \(\HH v=0\) on each \(\bx_{R_i}(x_i)\subseteq\mathcal C\).
Since \(\mathcal T_{x_i,s_i,R_i}v=v\) on \(\bx_{R_i-1}(x_i)^c\), an operation in another component changes neither the function nor its Hessian on \(\bx_{R_i}(x_i)\subseteq\mathcal C\).

Also, since \(\HH v=0\) on \(\bx_{R_i}(x_i)\), the first operation in \(\mathcal C\) changes neither the function nor its Hessian on \(\bx_{R_i}(x_i)\subseteq\mathcal C\), by \eqref{d2:eq:affinekernel} and \eqref{d2:eq:localmodsupport.1}. Induction shows every operation in \(\mathcal C\) is zero. By \eqref{d2:eq:localmodsupport}, we conclude $\supp\bigl(\HH(\mathcal Tv-v)\bigr)
 \subseteq \mathcal D(K)$. Hence $\supp( \HH(\mathcal Tv))
 \subseteq K\cup \mathcal D(K)$.
\end{proof}

\subsubsection{A spin wave for each square in $\cS_k^{sep}(\r)$}
\begin{proposition}\label{prop:4.7}
There is $C_E<\infty$ such that, for $s\in\cS_k^{sep}(\r)$,
$8\le k\le n(\r)$ and $\ell=2^k$, there is a function $a_{s,\r}$ with
\begin{align}
 \la a_{s,\r},\r\ra&\ge1,\label{d2:eq:wavepair}\\
 \supp( a_{s,\r})&\subset\Lda\cap\bx_{3\ell}(c_s),\label{d2:eq:wavesupport}\\
 \supp(\HH a_{s,\r})
 &\subset\{x:\ell/8\le\dist(x,\supp(\r))\le5\ell\},
       \label{d2:eq:waveannulus}\\
 \cE(a_{s,\r})&\le C_E,\label{d2:eq:waveenergy}\\
 a_{s,\r}&\text{ is affine on }\bx_{ad_\Lda(\eta)}(x_\eta)
       \quad(\eta\in\cN(\r)),\label{d2:eq:waveaffine}\\
 \supp(a_{s,\r})\cap\supp(\eta)&=\varnothing
       \quad(\eta\ne\r,\ \mom(\eta)\ne0).\label{d2:eq:activeavoid}
\end{align}
Here $s\subset\bx_{\ell/2}(c_s)$ and $c_s\in\Z^2$.
\end{proposition}

\begin{proof}
Fix $s\in\cS_k^{sep}(\r)$, $k\ge 8$, and define $\tau=\r\one_s$. Then
\begin{equation}\label{d2:eq:labelmoments}
 \tau\ne0,\qquad
 \mom(\tau)\ne0,\qquad
 \delta(\tau)\ge 1000 \ell,\qquad \ell\le d_\Lda(\r).
\end{equation}
If $|\cS_k(\r)|=1$, then these assertions follow from the definition of $\cS_k^{sep}(\r)$. On the other hand, if $|\cS_k(\r)|>1$, then these assertions follow from the definitions of $\cS_k(\r)$ and $\cS_k^{sep}(\r)$ and property~\ref{item:two_6.4} of Theorem~\ref{thm:2.2}.

Fix $\chi\in C_c^\infty((-2,2)^2)$ with $0\le\chi\le1$ and
$\chi=1$ on $[-1,1]^2$. Set $\chi_s(x)=\chi((x-c_s)/\ell)$. When $Q(\tau)=0$, let
$b_s=P(\tau)/|P(\tau)|_2$. Note that by \eqref{d2:eq:labelmoments}, $b_s$ is defined.

Define
\[
 a_{s,\r}^0(x)=
 \begin{cases}
 \sgn(Q(\tau))\chi_s(x),&Q(\tau)\ne0,\\
 b_s\cdot(x-c_s)\chi_s(x),&Q(\tau)=0.
 \end{cases}
\]
We have
\begin{equation*}\label{d2:eq:basepair}
 \la a_{s,\r}^0,\r\ra= \sum_{x\in s}  a_{s,\r}^0(x)\r(x)=
 \begin{cases}
 |Q(\tau)|,&Q(\tau)\ne0,\\
 b_s\cdot(P(\tau)-c_sQ(\tau))=|P(\tau)|_2,&Q(\tau)=0.
 \end{cases}
\end{equation*}
Since $Q(\tau)$ is integer-valued and $P(\tau)\in \Z^2\setminus \{0\}$ when $Q(\tau)=0$, the above implies $ \la a_{s,\r}^0,\r\ra \ge 1$.
For $i,j\in\{1,2\}$,
\begin{align*}
 D_i^+\chi_s(x)
 &=\ell^{-1}\int_0^1
       \partial_i\chi((x-c_s+t e_i)/\ell)\,\mathrm dt,\\
 (\HH\chi_s)_{ij}(x)
 &=\ell^{-2}\int_0^1\int_0^1
       \partial_{ij}\chi((x-c_s+t e_i-u e_j)/\ell)\,
       \mathrm dt\,\mathrm du.
\end{align*}
If $Q(\tau)= 0$, put $F(y)=b_s\cdot(y-c_s)\chi((y-c_s)/\ell)$.
Then, by the chain rule, we get
\[
 \partial_{ij}F(y)
 =\ell^{-1}(b_s)_i\partial_j\chi
  +\ell^{-1}(b_s)_j\partial_i\chi
  +\ell^{-2}b_s\cdot(y-c_s)\partial_{ij}\chi,
\]
where the derivatives of $\chi$ are evaluated at $(y-c_s)/\ell$.
Thus, for a constant $C_\chi$,
\begin{equation}\label{d2:eq:basehessian}
 \|\HH a_{s,\r}^0\|_\infty
 \le\begin{cases}
 C_\chi\ell^{-2},&Q(\tau)\ne0,\\
 C_\chi\ell^{-1},&Q(\tau)=0.
 \end{cases}
\end{equation}
Observe that, by definition, $\chi_s\equiv 1$ on $\bx_\ell(c_s)$ and $\chi_s\equiv 0$ on $\Z^2\setminus\bx_{2\ell}(c_s)$. Therefore, $a_{s,\r}^0$ is affine on $\bx_\ell(c_s)$ and
\[
 K_s:=\supp(\HH a_{s,\r}^0)
 \subset\bx_{2\ell+2}(c_s)\setminus\bx_{\ell-2}(c_s).
\]
Using $\supp(\tau)\subset\bx_{\ell/2}(c_s)$, separation from the
other covering squares, and $\ell\ge256$, we obtain
\begin{equation}\label{d2:eq:baseannulus}
 K_s\subset\{x:\ell/4\le\dist(x,\supp(\r))\le4\ell\}.
\end{equation}

Now consider the family $(x_\eta,d_\Lda(\eta))_{\eta\in\cN(\r)}$. By property~\ref{item:one_6.4} of Theorem~\ref{thm:2.2}, we have 
\[
|x_{\eta_1}-x_{\eta_2}|_\infty \ge \dist(\eta_1, \eta_2) \ge M [\min\{d_\Lda(\eta_1), d_\Lda(\eta_2)\}]^\a. 
\]
Now apply Lemma~\ref{d2:lem:sparse} to this family and set
$$a_{s,\r}=\mathcal T a_{s,\r}^0.$$ 

We now show that $a_{s,\r}$ satisfies the properties of the proposition. For any $\eta\in\cN(\r)$, we have by property~\ref{item:one_6.4} of Theorem~\ref{thm:2.2}
\begin{equation}\label{d2:eq:obstacleavoid}
 \dist(x_\eta,\supp(\r))
 \ge M2^{-\a}d_\Lda(\eta)^\a
 >\lceil Td_\Lda(\eta)^{\frac32}\rceil.
\end{equation}
Here we use \eqref{eq:def_alpha_M}. Therefore, by the definition of $a_{s,\r}$ and Lemma~\ref{d2:lem:localmodification}, it follows that
$\la a_{s,\r},\r\ra=\la a_{s,\r}^0,\r\ra\ge1$.

For any $\eta\in\cN(\r)$, let $D_\eta:=\bx_{\lceil Td_\Lda(\eta)^\frac32\rceil}(x_\eta)$. Let $\mathcal C$ be a component of $\bigcup_{\eta\in\cN(\r)} D_\eta$ such that $\mathcal C\cap K_s\ne \emptyset$. Let $\eta_*\in\cN(\r)$ such that $d_\Lda(\eta_*)=\max\{d_\Lda(\eta): \eta\in\cN(\r), D_\eta\subset \mathcal C\}$. Then equations \eqref{eq:4.15.1},
\eqref{d2:eq:obstacleavoid} and \eqref{d2:eq:baseannulus} imply
\[
 M2^{-\a}d_\Lda(\eta_*)^\a
 \le\dist(x_{\eta_*},\supp(\r))\le20Td_\Lda(\eta_*)^{3/2}+4\ell.
\]
The bounds in \eqref{eq:def_alpha_M} give
\[
d_\Lda(\eta_*)^\a\le \frac{8\ell\cdot2^\a}{M},\]
and hence\[
 \diam\mathcal C
 \le20T(8\cdot2^\a/M)^{3/(2\a)}\ell^{3/(2\a)}\le\ell/32.
\]

Writing $K_s^{+\ell/32}=\{x:\dist(x,K_s)\le\ell/32\}$,
Lemma~\ref{d2:lem:sparse}(iii) yields
\[
 \supp(a_{s,\r}-a_{s,\r}^0)\subset K_s^{+\ell/32},\qquad
 \supp(\HH a_{s,\r})\subset K_s^{+\ell/32}.
\]
By definition, $\supp(a_{s,\r}^0) \subset \bx_{2\ell}(c_s)$. Also, $K_s\subset \bx_{2\ell+2}(c_s)$. These, together with \eqref{d2:eq:baseannulus}, give
\begin{align}
 &\supp(a_{s,\r})\cup\supp(\HH a_{s,\r}) \subset \bx_{3\ell}(c_s),\label{eq:a_srho_1}\end{align}
and \begin{align}\supp(\HH a_{s,\r}) \subset \{x:\ell/8\le\dist(x,\supp(\r)) \le5\ell\}\label{eq:a_srho_2}.
\end{align}
Every point of $\bx_{3\ell}(c_s)$ lies within $4\ell$ of a point
of $\supp(\tau)$. By \eqref{d2:eq:labelmoments},
\begin{equation*}\label{d2:eq:wavedepth}
 \dist(\supp(a_{s,\r}),\Lda^c)\ge 996\ell.
\end{equation*}
This proves \eqref{d2:eq:wavesupport}.

From Lemma~\ref{d2:lem:sparse}(ii) and \eqref{d2:eq:basehessian} we get
\[
 \|\HH a_{s,\r}\|_\infty\le C_HC_\chi/\ell.
\]
As $ |\supp(\HH a_{s,\r})|\le(6\ell+1)^2\le49\ell^2,$ we have
\begin{equation}\label{d2:eq:correctedenergy}
 \cE(a_{s,\r})=\sum_x|\HH a_{s,\r}(x)|^2
 \le49C_H^2C_\chi^2\le C_E,
 \qquad C_E:=\max\{1,49C_H^2C_\chi^2\}.
\end{equation}
The affine assertion \eqref{d2:eq:waveaffine} follows from
Lemma~\ref{d2:lem:sparse}(i).

Finally, suppose $x\in\supp( a_{s,\r})\cap\supp(\eta)$ and
$\mom(\eta)\ne0$. Then
\[
 996\ell\le\dist(x,\Lda^c)
       \le\delta(\eta)+d(\eta)\le2d_\Lda(\eta),
 \qquad \dist(x,\supp(\r))\le4\ell.
\]
Since $d_\Lda(\r)\ge\ell$, by property~\ref{item:one_6.4} of Theorem~\ref{thm:2.2} we have
\[
 4\ell\ge\dist(\supp(\eta),\supp(\r))
       \ge M\min\{d_\Lda(\eta),d_\Lda(\r)\}^\a
       \ge M\ell^\a>4\ell.
\]
This proves \eqref{d2:eq:activeavoid}.
\end{proof}

\subsubsection{Proof of Proposition~\ref{prop:2.4}}
\begin{proof}
Fix $\r\in\cN$. In this proof, $\sum_{k,s}$ denotes
$\sum_{k=8}^{n(\r)}\sum_{s\in\cS_k^{sep}(\r)}$. For any $x\in\Z^2$, let $\delta_x=\dist(x,\supp(\r))$. If $\delta_x=0$, then by \eqref{eq:a_srho_2}, $\HH a_{s,\r}(x)=0$ for any $8\le k\le n(\r)$ and $s\in\cS_k^{sep}(\r)$. On the other hand, if $\delta_x>0$, then by \eqref{eq:a_srho_2}, for $s\in\cS_k^{sep}(\r)$, $\HH a_{s,\r}(x)\ne0$ implies $\delta_x/5\le 2^k\le8\delta_x.$
Also, for any $8\le k\le n(\r)$, if $s_1, s_2\in\cS_k^{sep}(\r)$, then by the definition of $\cS_k^{sep}(\r)$ we have $\bx_{3\cdot2^k}(c_{s_1})\cap \bx_{3\cdot2^k}(c_{s_2})=\emptyset$. These, together with~\eqref{eq:a_srho_1}, imply that for any $x$
\[
 \sum_{k,s}\one_{\{\HH a_{s,\r}(x)\ne0\}}\le 8.
\]
Now, by the Cauchy--Schwarz inequality, we obtain
\begin{align}\label{d2:eq:overlapenergy}
\cE\left(\sum_{k,s} a_{s,\r}\right) &=\sum_{x}\left|\HH\left(\sum_{k,s} a_{s,\r}(x)\right)\right|^2 \notag\\
&=\sum_{x}\left|\sum_{k,s}\HH a_{s,\r}(x)\right|^2 \notag\\
 &\le 8\sum_s\cE(a_{s,\r})\le 8 C_E I(\r).
\end{align}
Define
\begin{equation}\label{d2:eq:fullwave}
 a_\r=a_{0,\r}+t_\b\sum_{k,s} a_{s,\r},\qquad
 t_\b=\frac1{\b 8C_E}.
\end{equation}
By property~\ref{item:roman1} of Lemma~\ref{lem:4.5} and \eqref{d2:eq:waveannulus}, for any $8\le k\le n(\r)$ and $s\in\cS_k^{sep}(\r)$
\[
 \supp(\HH a_{0,\r})\cap\supp(\HH a_{s,\r})=\varnothing.
\]
Using \eqref{eq:d2-hessidentity}, property~\ref{item:roman4} of Lemma~\ref{lem:4.5}, \eqref{d2:eq:wavepair} and
\eqref{d2:eq:overlapenergy}, we obtain
\begin{align*}
 E_\b(a_\r,\r)&
 %=\la a_\r,\r\ra -\frac{\b}2 \la a_\r, \Delta^2 a_\r=\\
 = E_\b(a_{0,\r},\r) + t_\b\la \sum_{k,s}a_{s\r},\r\ra - \frac{\b}2 t_\b^2\cE\left(\sum_{k,s} a_{s,\r}\right)\\
 &\ge\frac1{144\b}\|\r\|_2^2
       +t_\b I(\r)-4\b t_\b^2 C_E I(\r)\\
&=\frac1{144\b}\|\r\|_2^2+\frac1{16\b C_E}I(\r).
\end{align*}
This proves~\ref{item:prop5_2.4}.

We now prove properties~\ref{item:prop1_2.4}--\ref{item:prop4_2.4} of the proposition. It is enough to prove these properties for $a_{0,\r}$ and $a_{s,\r}$ for $s\in\cS_k^{sep}(\r)$, $8\le k\le n(\r)$.

By property~\ref{item:roman1} of Lemma~\ref{lem:4.5}, $$\supp(a_{0,\r})\subseteq\supp(\r)\subset \Lda\cap \bx_{10d_\Lda(\r)}(x_\r).$$
For each $8\le k\le n(\r)$ and $s\in\cS_k^{sep}(\r)$, $\supp(\r) \cap s\ne\varnothing$ and
$2^k\le d_\Lda(\r)$. These, together with \eqref{d2:eq:wavesupport}, imply
\begin{equation*}\label{d2:eq:targetsupp}
 \supp( a_{s,\r})\subset \Lda\cap \bx_{10d_\Lda(\r)}(x_\r).
\end{equation*}
Hence~\ref{item:prop1_2.4} follows.
For $\eta\in\cN(\r)$,
\[
 \dist(x_\eta,\supp(\r))
       \ge M2^{-\a}d_\Lda(\eta)^\a> 100d_\Lda(\eta).
\]
Thus $a_{0,\r}=0$ on $\bx_{100d_\Lda(\eta)}(x_\eta)$. This, together with 
\eqref{d2:eq:waveaffine}, proves~\ref{item:prop2_2.4}. Property~\ref{item:prop3_2.4} follows from property~\ref{item:roman1} of Lemma~\ref{lem:4.5}, disjointness of the density supports and
\eqref{d2:eq:activeavoid}.

Next, we verify property~\ref{item:prop4_2.4}. Let $\eta\ne\r$. If $\mom(\eta)\ne 0$, then by property~\ref{item:prop3_2.4}, $\la a_\r,\eta\ra=0.$ If $\mom(\eta)= 0$ and $\eta\in\cN(\r)$, then by property $(2)$, $\la a_\r,\eta\ra=0.$ Now suppose $\eta\notin \cN(\r)$. Then $\dist(\supp(\r),\supp(\eta))\ge Md_\Lda(\r)^\a>10d_\Lda(\r)$ and hence by property~\ref{item:prop1_2.4}, $\la a_\r,\eta\ra=0.$

To prove the second identity, first assume $\r\in\cN(\eta)$. Then, by properties~\ref{item:prop1_2.4}--\ref{item:prop2_2.4}, $a_\eta$ is affine on $\bx_{100 d_\Lda(\r)}(x_\r)$ and 
$\supp(a_\r) \subset \bx_{10 d_\Lda(\r)}(x_\r)$. Hence $\Delta^2a_\eta=0$ on $\supp( a_\r)$, which implies $\la a_\r,\Delta_\Lda^2a_\eta\ra=0.$ Next assume $\r\notin\cN(\eta)$. In this case, $d_\Lda(\r)>2d_\Lda(\eta)$ and hence $\eta\in\cN(\r)$. By an argument similar to that above, we conclude $\la a_\eta,\Delta_\Lda^2a_\r\ra=0.$ But $\la a_\r,\Delta_\Lda^2a_\eta\ra=\la a_\eta,\Delta_\Lda^2a_\r\ra.$ Hence $\la a_\r,\Delta_\Lda^2a_\eta\ra=0$. 
This completes the proof.
\end{proof}

\appendix
\section{Gaussian comparison and the two-dimensional variance bounds}
\label{sec:proof}

\subsection{Proof of Proposition~\ref{thm:upper_bound}}\label{subsec:proof_Garban}
Write $n=|\Lda|$ and $A=\b\Delta_\Lda^2$, a positive-definite
$n\times n$ matrix. For $t\in\R^n$, set
\[
 \Theta_A(t)=\sum_{k\in\Z^n}
        \exp\!\left[-\frac12\la k-t,A(k-t)\ra\right].
\]
Poisson summation for the Schwartz function
$x\mapsto\exp[-\la x,Ax\ra/2]$ gives
\[
 \Theta_A(t)=
 \frac{(2\pi)^{n/2}}{\sqrt{\det A}}
 \sum_{q\in\Z^n}
       \e^{-2\pi^2\la q,A^{-1}q\ra}
       \cos(2\pi\la q,t\ra).
\]
The series is absolutely convergent. Every coefficient is non-negative,
so $\Theta_A(t)\le\Theta_A(0)$. Completing the square in the
discrete Gaussian sum now yields
\[
 \E^{IMM}_{\b,\Lda}\!\left[\e^{\la m,f\ra}\right]
 =
 \e^{\la f,A^{-1}f\ra/2}
 \frac{\Theta_A(A^{-1}f)}{\Theta_A(0)}
 \le \e^{\la\sigma_f,f\ra/(2\b)}.
\]
The last expression is the Gaussian Laplace transform in
Proposition~\ref{thm:upper_bound}. This proves the claim.

\subsection{Proof of Theorem~\ref{thm:ivmm2d}}
We use Theorem~\ref{thm:lower_bound_res}, whose proof is given in
Section~\ref{sec:lower}, with $\eps=1/2$.
Both laws are invariant under a change of sign.
Their finite-volume Gaussian tails justify differentiating
the Laplace transforms twice at zero. Proposition~\ref{thm:upper_bound} and Theorem~\ref{thm:lower_bound_res} therefore imply
\[
 \frac{1}{2\b}G_\Lda(0,0)
 \le \E^{IMM}_{\b,\Lda}[m(0)^2]
 \le \frac1\b G_\Lda(0,0).
\]
The diagonal estimate of
\citet[Theorem~1.1]{muller2019estimates} gives
$cL^2\le G_\Lda(0,0)\le CL^2$ for $L\ge1$ and hence the result follows.

{\section{Discrete Gaussians and positivity of covariance for IVDGFF}}
\subsection{Discrete Gaussians}\label{app:discrete-gaussians}
A \emph{lattice} $K\subset\R^n$ of rank $k$ is the set of all integer linear combinations of $k$ linearly independent vectors in $\R^n$. For any $s>0$, we define the function $\rho_s:\R ^n\to \R$ as
\begin{align*}
  \rho_s(x)=\exp(-\pi\|x\|^2/s^2).  
\end{align*}
For a discrete set $A\subset\mathbb{R}^n$ we define $\rho_s(A)=\sum_{x\in A}\rho_s(x)$. The {\emph{discrete Gaussian $X_{K+w,s}$}} on a lattice coset $K+w$ with parameter $s$ is defined via its probability mass function 
\begin{align*}
\prob(X_{K+w,s}=x)= \frac{\rho_s(x)}{\rho_s(K+w)}, \qquad x\in K+w.  
\end{align*}
When $s=1$, we write $\rho(x)$ and $X_{K+w}$.

The following inequality is stated in \citet[Theorem~2.1]{RSD17-discrete-gaussian-ineq}.

\begin{theorem}[{Lattice Gaussian inequality}]\label{thm:RSD}
For every lattice $K\subset\R^n$ and every $x,y\in\R^n$,
\begin{equation}\label{eq:RSD}
  \rho(K+x)^2\rho(K+y)^2
  \le
  \rho(K)^2\rho(K+x+y)\rho(K+x-y).
\end{equation}
\end{theorem}

We now derive the form that we will use. The following proposition shows that the covariance of a discrete Gaussian decreases in the positive-semidefinite order on symmetric matrices under restriction to a sublattice. A result in this spirit was already obtained in~\citet[Proposition 2.1]{AizenmanHarelPeledShapiro}.

\begin{proposition}[]\label{prop:sublattice-cov}
Let $M\subset L\subset\R^n$ be lattices (their ranks may be different). Then
\begin{equation}\label{eq:cov-order-euclidean}
  \cov(X_M)\preceq\cov(X_L)
\end{equation}
in the positive-semidefinite order on symmetric matrices.

More generally, let $V$ be a finite-dimensional real vector space, let $Q$ be a positive-definite quadratic form on $V$, and let $M\subset L\subset V$ be discrete additive subgroups.  If, for $K\in\{M,\,L\}$,
\begin{align*}
  \prob(Y_K=y)\propto e^{-\b Q(y)/2},\qquad y\in K,
\end{align*}
then for every linear functional $f$ on $V$,
\begin{equation}\label{eq:cov-order-Q}
  \var(f(Y_M))\le\var(f(Y_L)).
\end{equation}
\end{proposition}

\begin{proof}
We first prove \eqref{eq:cov-order-euclidean}.  Choose one representative
$c$ from each coset of $M$ in $L$.  The quotient $L/M$ may be infinite;
all sums below are absolutely convergent because of Gaussian decay.
Apply Theorem~\ref{thm:RSD} to the lattice $M$ with the vectors $x$ and
$c$.  Taking square roots gives
\[
  \frac{\rho(M+x)}{\rho(M)}\,\rho(M+c)
  \le
  \sqrt{\rho(M+c+x)\rho(M+c-x)}.
\]
Using the arithmetic--geometric mean inequality on the right-hand side,
\[
  \frac{\rho(M+x)}{\rho(M)}\,\rho(M+c)
  \le
  \frac{\rho(M+c+x)+\rho(M+c-x)}2.
\]
Sum over all coset representatives $c$.  Since the cosets partition
$L$,
\[
  \sum_c\rho(M+c)=\rho(L),
\]
and similarly, since a lattice is closed under the negation operation $\ell\mapsto -
\ell$
\[
  \sum_c\rho(M+c+x)=\rho(L+x),
  \qquad
  \sum_c\rho(M+c-x)=\rho(L-x)=\rho(L+x).
\]
Therefore
\begin{equation}\label{eq:theta-order}
  f_M(x):=\frac{\rho(M+x)}{\rho(M)}
  \le
  \frac{\rho(L+x)}{\rho(L)}=:f_L(x).
\end{equation}
Both functions are even and smooth and satisfy $f_M(0)=f_L(0)=1$.
Hence $f_L-f_M$ has a global minimum equal to $0$ at the origin, and
therefore
\begin{equation}\label{eq:f_posdef}
  \nabla^2f_L(0)-\nabla^2f_M(0)\succeq0.
\end{equation}
Recall that for a lattice $K\subset \R^n$ we have $f_K(x)=\rho(K+x)/\rho(K)$. Consider now, for {$x\in\R^n$ and $w\in K$},
\[
g_w(x)={\exp(-\pi\|w+x\|^2)}.
\]
Its Hessian is
\[
\nabla^2 g_w(x)=\left(-2\pi \mathbb{I}_n+4\pi^2(x+w)(x+w)^T\right)\exp(-\pi\|w+x\|^2),
\]
where $\mathbb{I}_n$ is the $n\times n$ identity matrix, so

\begin{align*}
 \nabla^2 f_K(0)
 &=\frac{\sum_{w\in K}
   \left.\nabla_x^2\e^{-\pi\|w+x\|^2}\right|_{x=0}}{\rho(K)}\\
 &=-2\pi\mathbb I_n+
    4\pi^2\frac{\sum_{w\in K}ww^T\e^{-\pi\|w\|^2}}{\rho(K)}.
\end{align*}

Since $X_K$ is centred,
$\rho(K)^{-1}\sum_{w\in K}ww^T\e^{-\pi\|w\|^2}=\cov(X_K)$.
Therefore

\begin{equation}\label{eq:hessian-theta}
  \nabla^2f_K(0)
  =-2\pi \mathbb{I}_n+4\pi^2\cov(X_K).
\end{equation}
We substitute \eqref{eq:hessian-theta}, proved for a general $K$, into~\eqref{eq:f_posdef} with $K\coloneqq L$ and $K\coloneqq M$, respectively, and this proves~\eqref{eq:cov-order-euclidean}.

For the general quadratic form, choose an invertible linear map
$T:V\to\R^n$ such that
\[
  \frac\b2Q(v)=\pi\|Tv\|^2.
\]

To see this, choose coordinates in which $Q(v)=v^THv$
with $H$ symmetric positive definite, and take
$T=(\b/(2\pi))^{1/2}H^{1/2}$.

Then $TM\subset TL$ are lattices. For $K\in\{M,L\}$, the
restriction $T:K\to TK$ is a bijection and
\[
 \prob(Y_K=y)
 =\frac{\e^{-\b Q(y)/2}}{\sum_{k\in K}\e^{-\b Q(k)/2}}
 =\frac{\e^{-\pi\|Ty\|^2}}{\rho(TK)}
 =\prob(X_{TK}=Ty).
\]

Applying \eqref{eq:cov-order-euclidean} to the linear functional
$\phi\coloneqq f\circ T^{-1}$ on the lattices $TM$ and $TL$ gives
\begin{align*}   \var(\phi(X_{TM}))\le\var(\phi(X_{TL}))\iff \var(f(Y_M))\le\var(f(Y_L)),
\end{align*}
which proves~\eqref{eq:cov-order-Q}.
\end{proof}

The next lemma is a consequence of Proposition~\ref{prop:sublattice-cov}.

\begin{lemma}\label{lem:add-squares}
Let $U=\Z^N$.  Let
\begin{align*}
  A:U\to\Z^p,
  \qquad
  B:U\to\Z^q
\end{align*}
be integer-linear maps, and assume that $A$ is injective.  Define
\begin{align*}
  \mu_{A,B}(u)
  &\propto
  \exp\!\left[-\frac\b2
      \bigl(\|Au\|^2+\|Bu\|^2\bigr)\right],\\
  \mu_A(u)
  &\propto
  \exp\!\left[-\frac\b2\|Au\|^2\right].
\end{align*}
Then for every linear functional $f:\R^N\to\R$,
\begin{equation}\label{eq:add-squares}
  \var_{\mu_{A,B}}(f(u))
  \le
  \var_{\mu_A}(f(u)).
\end{equation}
\end{lemma}
\begin{proof}
   Consider the two lattices in $\R^p\oplus\R^q$,
\[
  \Gamma:=\{(Au,Bu):u\in U\},
  \qquad
  \Gamma':=AU\times BU.
\]
Because $A$ is injective, the map $u\mapsto(Au,Bu)$ is injective, and
$\Gamma$ is a lattice.  Since $A$ and $B$ have integer coefficients,
$AU$ and $BU$ are discrete additive subgroups, so $\Gamma'$ is also a
lattice.  Clearly $\Gamma\subset\Gamma'$.

The map $A:U\to AU$ is a bijection.  Hence there is a linear functional
$\ell$ on $\operatorname{span}(AU)$ satisfying
\[
  \ell(Au)=f(u).
\]
{Extend $\ell$ linearly to $\R^p$, and then to $\R^p\oplus\R^q$ by ignoring the second coordinate:}
\[
  F(a,b):=\ell(a).
\]
The measure $\mu_{A,B}$ is exactly the pullback under
$u\mapsto(Au,Bu)$ of the centred discrete Gaussian on $\Gamma$ with
quadratic form $\|a\|^2+\|b\|^2$.  On the product lattice
$\Gamma'$, the same Gaussian factorises into an independent product of a
discrete Gaussian on $AU$ and one on $BU$.  Therefore the variance of
$F$ under the Gaussian on $\Gamma'$ is exactly
{$\var_{\mu_A}(f(u))$}.

Applying Proposition~\ref{prop:sublattice-cov} to
$\Gamma\subset\Gamma'$ proves \eqref{eq:add-squares}.
\end{proof}

\subsection{Positive association of the IVDGFF}\label{sec:appendix_B.2}
\begin{proof}[Proof of Lemma~\ref{lem:FKG}]
For $M\in\mathbb N$, let
\(
  \Omega_M:=\{-M,-M+1,\ldots,M\}^{D}.
\)
Let $\prob^{IGFF, (M)}_{\b,D}$ be the restriction of $\prob^{IGFF}_{\b,D}$ to $\Omega_M$, renormalized to have total mass one. That is,
\[
  \prob^{IGFF, (M)}_{\b,D}(h=g)
  =\frac{1}{Z^{IGFF, (M)}_{\b, D}} \exp\left[-\frac{\b}{2} \la g,\, -\Delta g \ra\right] \one_{\{g\in\Omega_M\}} \prod_{x\in\Z^d\setminus D}\one_{\{0\}}(g(x)).
\]
By Claim~$A.1$ of \cite{chen2026limitinglawdiscretegaussian} we conclude that $\prob^{IGFF, (M)}_{\b,D}$ satisfies the FKG lattice condition.  
%The finite-lattice FKG theorem~{\citep{FKG1971}} therefore yields positive association. 
Now, since the coordinate
functions $h\mapsto h(x)$ and $h\mapsto h(y)$ are bounded increasing
functions on $\Omega_M$, we have
\begin{equation}\label{eq:truncated-height-covariance-positive}
  \cov_{\prob^{IGFF, (M)}_{\b,D}}(h(x),h(y))\ge0.
\end{equation}

It remains to pass to $M\to\infty$.  The quadratic form $\la g, -\Delta g\ra$ is positive definite on $\R^D$. Because $D$ is finite, there exists $\lambda_D>0$ such that
\begin{equation}\label{eq:Dirichlet-spectral-gap-finite-D}
  \la h, -\Delta h\ra \ge\lambda_D\sum_{u\in D}h(u)^2,
  \qquad h\in\R^D.
\end{equation}
Consequently,
\begin{equation}\label{eq:FKG-moment-dominating-sum}
  \sum_{h\in\Z^D}
  \left(1+\sum_{u\in D}h(u)^2\right)
  \exp\!\left[-\frac\b2\la h, -\Delta h\ra\right]
  <\infty.
\end{equation}
Indeed, by \eqref{eq:Dirichlet-spectral-gap-finite-D}, the summand is
bounded by a polynomial times
$\exp[-\b\lambda_D\|h\|_2^2/2]$, which is summable over $\Z^D$.

Since $\Omega_M\uparrow\Z^D$, the partition functions converge to
$Z_{\b,D}^{IGFF}$.  The same domination applies to the first moments and to
the mixed second moment, because
\[
  |h(x)|\le 1+h(x)^2,
  \qquad
  |h(x)h(y)|\le\frac12(h(x)^2+h(y)^2).
\]
Therefore
\[
 \lim_{M\to\infty} \E_{\prob^{IGFF, (M)}_{\b,D}}[h(x)]
    = \E_{\prob^{IGFF}_{\b,D}}[h(x)],
  \qquad
 \lim_{M\to\infty}  \E_{\prob^{IGFF, (M)}_{\b,D}}[h(y)]
    =\E_{\prob^{IGFF}_{\b,D}}[h(y)],
\]
and
\[
 \lim_{M\to\infty}  \E_{\prob^{IGFF, (M)}_{\b,D}}[h(x)h(y)]
    = \E_{\prob^{IGFF}_{\b,D}}[h(x)h(y)].
\]
It follows that
\[
  \lim_{M\to\infty} \cov_{\prob^{IGFF, (M)}_{\b,D}}(h(x),h(y))
  =
  \cov_{\prob^{IGFF}_{\b,D}}(h(x),h(y)).
\]
Passing to the limit in
\eqref{eq:truncated-height-covariance-positive} proves the lemma.
\end{proof}

\section{Proofs of the convex expansion and covering bound (Proposition~\ref{prop:2.1} and Theorem~\ref{thm:2.2})}
\label{d2:sec:regrouping}

\subsection{Proof of Proposition~\ref{prop:2.1}}
\label{d2:sec:count}
In this section we prove Proposition~\ref{prop:2.1}. The proof is similar to that of Proposition 20 of \cite{Wirth_2019} and Proposition 2.1 of \cite{K-P_2017}.

To prove Proposition~\ref{prop:2.1} we use the following lemma, which allows us to bound recursively $|\cS_k(\r)|$. Set 
\begin{align*}
    b&:=\log_2\left(3+4M 10^{3\a}\right),\\
    \g(k)&:=\floor{\a^{-1}(k-b)} \text{ for integers } k\ge 0.
\end{align*}
The following lemma is analogous to Lemma~4.1 of \cite{K-P_2017}.
\begin{lemma}[]\label{lem:4.1}
Let $k\ge b$ be an integer. If $|\cS_{\g(k)}(\r)|\ge 2$, then 
\begin{align}\label{eq:l4.1.1}
 |\cS_{k}(\r)| \le \frac12|\cS_{\g(k)}(\r)| + |\cS^{sep}_{\g(k)}(\r)|  .
\end{align}
\end{lemma}
We prove Lemma~\ref{lem:4.1} by using the following claim from \cite{K-P_2017}.
\begin{claim}[Claim 4.6 of \cite{K-P_2017}]\label{claim:4.2}
    Let $G=(V(G), E(G))$ be a finite, non-empty, simple graph with no isolated vertices. Then $V(G)$ is the union of $\floor{|V(G)|/2}$ sets, each of which has size $2$ or $3$ and is connected in $G$.
\end{claim}
Now we prove Lemma~\ref{lem:4.1} and then Proposition~\ref{prop:2.1}.
\begin{proof}[Proof of Lemma~\ref{lem:4.1}]
 Since $|\cS^{"}_{\g(k)}(\r)|\le |\cS_{\g(k)}(\r)|$, it is enough to show 
 \begin{align*}
 |\cS_{k}(\r)| \le \frac12|\cS^{"}_{\g(k)}(\r)| + |\cS^{sep}_{\g(k)}(\r)|.
\end{align*}
We first note that $\supp(\r)$ is covered by $\cS^{sep}_{\g(k)}(\r)\cup\cS^{"}_{\g(k)}(\r)$. Also, $\g(k)<k$ for $k\ge b$. If $\cS^{"}_{\g(k)}(\r)=\emptyset$, then
 \begin{align*}
 |\cS_{k}(\r)| \le |\cS_{\g(k)}(\r)| = \frac12|\cS^{"}_{\g(k)}(\r)| + |\cS^{sep}_{\g(k)}(\r)|.
 \end{align*}
 So, we assume $\cS^{"}_{\g(k)}(\r)\neq\emptyset$. In this case, we complete the proof by showing that $\cS^{"}_{\g(k)}(\r)$ can be covered by at most $|\cS^{"}_{\g(k)}(\r)|/2$ squares of side length $2^k$.

 To see this, consider a graph $G$ in which the vertex set is $\cS^{"}_{\g(k)}(\r)$ and two squares $s_1, s_2$ are adjacent if $\dist(s_1, s_2) < 2M\left(1000\cdot 2^{\g(k)}\right)^{\a}$. The definitions of $\cS^{"}_{\g(k)}(\r)$, $\cS^{sep}_{\g(k)}(\r)$ and the assumption $|\cS_{\g(k)}(\r)|\ge 2$ imply that there are no isolated vertices in $G$. Applying Claim~\ref{claim:4.2} to $G$ reduces our task to showing that if $\{s_1,s_2,s_3\}\subseteq V(G)$ is connected, then $s_1,s_2,s_3$ may be covered by a single square of side length $2^k$.

 Suppose $\{s_1,s_2,s_3\}\subseteq V(G)$ is connected. Without loss of generality, suppose that $\dist(s_1, s_2)$ and $ \dist(s_2, s_3) $ are at most $ 2M\left(1000\cdot 2^{\g(k)}\right)^{\a}$. Then 
 \begin{align*}
     d(s_1\cup s_2\cup s_3) &\le d(s_1) + d(s_2) + d(s_3) + \dist(s_1, s_2) + \dist(s_2, s_3)\\
     &< 3\cdot2^{\g(k)} + 2\cdot 2M\left(1000\cdot 2^{\g(k)}\right)^{\a}\\
     & \le2^{\a\g(k)} \left(3+4M 10^{3\a}\right)\\
     &\le 2^{\a\g(k) + \log_2\left(3+4M 10^{3\a}\right)} \le 2^k.
\end{align*}
This completes the proof of the lemma.
\end{proof}
\begin{proof}[Proof of Proposition~\ref{prop:2.1}]
 Fix a density $\r$. We have
 \begin{align*}
     A(\r) = \sum_{k=0}^{n(\r)} |\cS_k(\r)| \ge 1+ n(\r) = 1+ \ceil{\log_2(M. d_\Lda(\r)^\a)} \ge \log_2\left(d_\Lda(\r)+1\right).
 \end{align*}
 We now prove the upper bound of the statement. 

Define $n^{'}(\r):=\lceil \log_2(M\cdot d(\r)^\a)\rceil$ and $A^{'}(\r):=\sum_{k=0}^{n^{'}(\r)} |\cS_k(\r)|$. We make the following claim: there exists an absolute constant $D>0$ such that
\begin{align}\label{eq:p2.1.7}
    A^{'}(\r) \le D \left( |\cS_0(\r)| + \sum_{k=1}^{n(\r)} |\cS^{sep}_k(\r)| \right).
\end{align}
 We first assume \eqref{eq:p2.1.7} and complete the proof of the proposition. Then we will prove \eqref{eq:p2.1.7}.

 If $d_\Lda(\r)=d(\r)$, then $A(\r)=A^{'}(\r)$. So, let $d_\Lda(\r)=\delta(\r)>d(\r)$. Note that $A(\r)-A^{'}(\r)= n(\r)-n^{'}(\r) \le \log_2(M\cdot d_\Lda(\r)^\a) +1 - \log_2(M\cdot d(\r)^\a) = \a\log_2(\delta(\r)/d(\r)) +1$. Let $m(\r) = \floor{\log_2(\delta(\r)/1000)}$ and $m^{'}(\r) = \ceil{\log_2(d(\r))+1}$. Note that $m(\r)\le n(\r)$ and 
 for $ m^{'}(\r) \le k \le m(\r)$, $|\cS^{sep}_k(\r)|=|\cS_k(\r)|=1$. Since $m(\r)-m^{'}(\r)\ge \log_2(\delta(\r)/d(\r))-\log_2(1000) -1$, we conclude that there exists an absolute constant $D^{'}>0$ such that
 \begin{align}\label{eq:p2.1.8}
    A(\r)-A^{'}(\r) \le D^{'}\left(1+ \sum_{k=m^{'}(\r)}^{m(\r)}|\cS^{sep}_k(\r)| \right) .
 \end{align}
 From \eqref{eq:p2.1.7} and \eqref{eq:p2.1.8} we get
 \begin{align*}
    A(\r) \le D \left( |\cS_0(\r)| + \sum_{k=1}^{n(\r)} |\cS^{sep}_k(\r)| \right).
\end{align*}
 Now the desired upper bound for $A(\r)$ follows from the observation that $|\cS^{sep}_k(\r)| \le|\cS_k(\r)|\le |\cS_0(\r)|$ for any $k\ge 1$, in particular for $1\le k \le 7$.
 
 We now prove \eqref{eq:p2.1.7}. We will use the recursive bound for $\cS_k(\r)$ from Lemma~\ref{lem:4.1}. 
 In order to apply \eqref{eq:l4.1.1}, we first check that $|\cS_{\g(k)}(\r)|\ge 2$ for $b\le k\le n^{'}(\r)$. Indeed,
 \begin{align*}
     \g(k) &\le \a^{-1}(k-b)\\
     & \le \a^{-1}\left(\log_2(M\cdot d(\r)^\a) +1 -\log_2\left(3+4M 10^{3\a}\right)\right)\\
     & \le \a^{-1}\left(\log_2 M +\a\log_2(d(\r))+1 -2 -\log_2 M -3\a\log_2(10)\right)\\
    & =\log_2(d(\r)) -\frac{1}{\a} -3\log_2(10),
 \end{align*}
 and hence $|\cS_{\g(k)}(\r)|\ge 2$.
 For $k\ge b$, we now iterate \eqref{eq:l4.1.1} $\ell(k)$
 times, where $\ell(k)$ is the maximum number such that $\g^{\ell(k)}(k)\ge 0$ ($\g^{m}$ denotes the $m$-fold composition of $\g$ with itself). Thus
 \begin{align}\label{eq:p2.1.1}
     |\cS_k(\r)| &\le \frac12|\cS_{\g(k)}(\r)|+|\cS_{\g(k)}^{sep}(\r)| \nonumber\\
     &\le  \frac1{2^2}|\cS_{\g^2(k)}(\r)|+|\cS_{\g(k)}^{sep}(\r)| + \frac1{2}|\cS_{\g^2(k)}^{sep}(\r)|\nonumber\\
     & \vdots \nonumber\\
     & \le 2^{-\ell(k)} |\cS_{\g^{\ell(k)}(k)}(\r)| + \sum_{m=0}^{\ell(k)-1} 2^{-m}|\cS_{\g^{m+1}(k)}^{sep}(\r)|.
 \end{align}
 We now estimate $\ell(k)$. Fix $k\ge b$. One can check by induction that for $0\le m\le \ell(k)$,
 \begin{align}\label{eq:p2.1.2}
     \a^{-m}k - \frac{\a+b}{\a-1} \le \a^{-m}k-b\sum_{j=1}^m \a^{-j} - \sum_{j=0}^{m-1} \a^{-j} \le \g^m(k) \le \a^{-m}k-b\sum_{j=1}^m \a^{-j}.
\end{align}
Set $k_0:=({\a+b})/({\a-1})$. Then from the lower bound we have
\begin{align}\label{eq:p2.1.3}
    \ell(k) \ge\begin{cases}
        0 \qquad \qquad \qquad 0\le k\le k_0\\
        \frac{\log_2(k/k_0)}{\log_2(\a)} \qquad \quad \text{ otherwise}
    \end{cases}.
\end{align}
Let $m,j\ge 0$. We estimate the cardinality $|N_{m,j}|$ of the set
$$N_{m,j}:=\{k: \g^m(k)=j\}.$$
Let $k_{-}$ and $k_{+}$ be the minimum and maximum integers in $N_{m,j}$, respectively. Then, using \eqref{eq:p2.1.2}, we have
\begin{align*}
  |N_{m,j}|& = k_{+} - k_{-} +1\\
  & \le \a^m \left[ j + b\sum_{\ell=1}^m \a^{-\ell} +\sum_{\ell=0}^{m-1} \a^{-\ell} -j -b\sum_{\ell=1}^m \a^{-\ell}\right] +1 \\
  &\le \a^m \sum_{\ell=0}^{m-1} \a^{-\ell} +1
%  &\le \a^m \sum_{\ell=0}^{\infty} \a^{-\ell} +1\\
%  &\le \a^m \frac{\a}{\a-1} +1\\
  \le \a^m \frac{2\a}{\a-1}.
\end{align*}
Now, using \eqref{eq:p2.1.1} together with the observation that $\ell(k)=0$ for $k<b$, we have
\begin{align}\label{eq:p2.1.4}
    A^{'}(\r)&=\sum_{k=0}^{n^{'}(\r)} |\cS_k(\r)|\nonumber\\
    &\le |\cS_0(\r)|\sum_{k=0}^{n^{'}(\r)} 2^{-\ell(k)} + \sum_{k=0}^{n^{'}(\r)}\sum_{m=0}^{\ell(k)-1} 2^{-m}|\cS_{\g^{m+1}(k)}^{sep}(\r)| .
\end{align}
We now bound both terms above. First, note that using \eqref{eq:p2.1.3} and the assumption that $\a<2$ we have
\begin{align}\label{eq:p2.1.5}
 \sum_{k=0}^{n^{'}(\r)} 2^{-\ell(k)} &\le \sum_{k=0}^{\infty} 2^{-\ell(k)} \nonumber \\
 & \le \lceil k_0 \rceil + 1 + \sum_{k=\lceil k_0 \rceil +1}^{\infty} 2^{-\frac{\log_2(k/k_0)}{\log_2(\a)}} \nonumber\\
 &= \lceil k_0 \rceil + 1 +  \sum_{k= \lceil k_0 \rceil +1}^{\infty} (k_0/k)^{\frac1{\log_2(\a)}} <\infty.
\end{align}
 
Therefore
\begin{align}\label{eq:p2.1.6}
   \sum_{k=0}^{n^{'}(\r)}\sum_{m=0}^{\ell(k)-1} 2^{-m}|\cS_{\g^{m+1}(k)}^{sep}(\r)| & = \sum_{k=0}^{n^{'}(\r)}\sum_{m=0}^{\ell(k)-1} 2^{-m} \sum_{j=0}^{n(\r)} \one_{\g^{m+1}(k)}(j)|\cS_j^{sep}(\r)|\nonumber\\
    %&\le
%\sum_{j=0}^{k(\r)}\sum_{k=0}^{\infty}\sum_{m=0}^{\infty} 2^{-m}\one_{\g^{m+1}(k)}(j)|\cS_j^{sep}(\r)|\nonumber\\
&\le \sum_{j=0}^{n(\r)}\sum_{m=0}^{\infty} 2^{-m} |N_{m+1,j}| |\cS_j^{sep}(\r)|\nonumber\\
& \le \sum_{j=0}^{n(\r)} \a \frac{2\a}{\a-1} |\cS_j^{sep}(\r)| \sum_{m=0}^{\infty} \left(\frac{\a}{2}\right)^m\nonumber\\
%& \le \frac{2\a^2}{\a-1} \frac{1}{1-\frac{\a}{2}} \sum_{j=0}^{k(\r)} |\cS_j^{sep}(\r)|\nonumber\\
& \le \frac{4\a^2}{(\a-1)(2-\a)} \sum_{j=0}^{n(\r)} |\cS_j^{sep}(\r)|.
\end{align}
Substituting \eqref{eq:p2.1.5} and \eqref{eq:p2.1.6} into \eqref{eq:p2.1.4}, we obtain the desired upper bound for $A^{'}(\r)$ in \eqref{eq:p2.1.7}.
\end{proof}

\subsection{Proof of Theorem~\ref{thm:2.2}}
In this section we prove Theorem~\ref{thm:2.2}. We first introduce the following notation. A density $\r_1$ is said to be \textit{compatible} with an ensemble $\cE$ if
\begin{align*}
    \r_1=\sum_{\r\in\cE} \eps(\r_1, \r)\r\qquad \text{ with } \quad \eps(\r_1, \r) \in \{-1, 0, 1\}.
\end{align*}
Note that the coefficients $\eps(\r_1, \r)$ are uniquely determined, as the densities in $\cE$ have disjoint supports.

We say an ensemble $\cE_1$ is a \textit{parent} of an ensemble $\cE_2$, and write $\cE_1\to\cE_2$, when every charge density $\r\in \cE_2$ is compatible with $\cE_1$.

For an integer $k\ge -1$, an ensemble $\cE$ is said to be an \textit{$k$-ensemble} if
\begin{align*}
    \dist(\r_1, \r_2) > 2^k \quad \text{ for all }\r_1,\r_2\in\cE,\,\r_1\neq\r_2.
\end{align*}
Also, let $A_k(\r):=|\cS_k(\r)|$ for $k\ge 0$ and set $A_{-1}(\r):=A_0(\r)=|\supp(\r)|$.

The following lemma is the same as Lemma~4.3 of \cite{K-P_2017}, except that we use the $\ell_\infty$-norm distance on $\Z^2$, rather than the graph distance.
\begin{lemma}\label{lem:4.3}
  Let $k\ge 0$, let $\cE$ be an ensemble and $(K(\r))_{\r\in\cE}$ be real.  There exist a positive absolute constant $C_1$, a finite family of $k$-ensembles $\cF^{'}$, with $\cE\to\cE^{'}$ for every $\cE^{'}\in\cF^{'}$, positive coefficients $(c_{\cE^{'}})_{\cE^{'}\in\cF^{'}}$ summing to $1$, and real coefficients $(K^{'}(\r))_{\r\in\cE^{'},\, \cE^{'}\in\cF^{'}}$, such that for every $\psi:\Lda\to\R$,
  \begin{align}\label{eq:l4.3.1}
      \prod_{\r\in\cE}\left(1+ K(\r) \cos\la\psi,\r\ra\right) = \sum_{\cE^{'}\in\cF^{'}}c_{\cE^{'}}\prod_{\r\in\cE^{'}}\left(1+ K^{'}(\r) \cos\la\psi,\r\ra\right).
  \end{align}
  Moreover, for every $\cE^{'}\in\cF^{'}$ and $\r\in\cE^{'}$ the following properties are satisfied:
  \begin{enumerate}[label=(\arabic*),ref=(\arabic*)]
      \item\label{item:C3_one} If there is a density $\r_1\subset \r$, compatible with $\cE$, then $\dist(\r_1,\r -\r_1) \le 2^k$.
      \item\label{item:C3_two} Let $\eps(\r, \r^{'}) \in \{-1, 0, 1\}$ be such that $\r=\sum_{\r^{'}\in\cE}\eps(\r, \r^{'})\r^{'}$. Then
      \begin{align*}
          |K^{'}(\r)| \le 3^{|\{\r_1\in\cE:\dist(\r,\r_1)\le 2^k\}|}\prod_{\r^{'}\in\cE} |K(\r^{'})|^{|\eps(\r, \r^{'})|}.
      \end{align*}
      Moreover, if $\cE$ is a $(k-1)$-ensemble, then 
      \begin{align}\label{eq:l4.3.2}
          |K^{'}(\r)| \le \e^{C_1A_{k-1}(\r)}\prod_{\r^{'}\in\cE} |K(\r^{'})|^{|\eps(\r, \r^{'})|}.
      \end{align}
  \end{enumerate}
\end{lemma}

\begin{proof}[Proof of Lemma~\ref{lem:4.3}]
Observe that in \eqref{eq:l4.3.1}, the left-hand side is expressed as a convex combination of products of expressions of the form $1+ K\cos\la\psi, \r\ra$. To prove \eqref{eq:l4.3.1}, we iteratively apply the following trigonometric identity:
  \begin{align}\label{eq:l4.3.3}
      \left(1+K_1 \cos\la\psi, \r_1\ra \right)&\left(1+K_2 \cos\la\psi, \r_2\ra \right) = \frac13 \left(1+ 3K_1 \cos\la\psi, \r_1\ra \right) + \frac13\left(1+ 3K_2 \cos\la\psi, \r_2\ra \right)\nonumber\\ 
      & + \frac16\left(1+3 K_1K_2 \cos\la\psi, \r_1-\r_2\ra \right) + \frac16\left(1+3 K_1K_2 \cos\la\psi, \r_1+\r_2\ra \right).
  \end{align}
  The above identity expresses the left-hand side as a convex combination of expressions of the form $1+ K\cos\la\psi, \r\ra$. Note that all densities on the right-hand side of \eqref{eq:l4.3.3} are compatible with $\{\r_1, \r_2\}$.
  If there are no densities $\r_1,\r_2 \in \cE$ such that $\dist(\r_1, \r_2) \le 2^k$, then $\cE$ is a $k$-ensemble and hence the lemma holds trivially. 
  
  Otherwise, we start by applying \eqref{eq:l4.3.3} to two arbitrary factors on the left-hand side of \eqref{eq:l4.3.1} corresponding to two densities $\r_1, \r_2\in\cE$, for which 
  \begin{align}\label{eq:l4.3.4}
      \dist(\r_1, \r_2) \le 2^k.
  \end{align}
  The right-hand side of \eqref{eq:l4.3.3} is inserted on the left-hand side of \eqref{eq:l4.3.1}, replacing these two factors by a sum of four terms. Then, by expanding, we obtain a sum of four products. If one of the resulting products contains two factors corresponding to two charge densities, $\r_1^{'}, \r_2^{'}$, satisfying \eqref{eq:l4.3.4}, we apply the identity \eqref{eq:l4.3.3} again and expand the resulting expression into a sum of products. We repeat this operation until we obtain a sum over products indexed by ensembles $\cE^{'}$ with the property that  $\dist(\r_1, \r_2) > 2^k$ for any $\r_1, \r_2\in\cE^{'}$. Note that after every application of \eqref{eq:l4.3.3}, every product still corresponds to an ensemble, i.e. all the densities in the product have disjoint supports.

  Clearly, the process yields $c_{\cE^{'}}>0$, $\sum_{\cE^{'}} c_{\cE^{'}} =1$, and $\cE \to \cE^{'}$ for each $\cE^{'}$. Also, property~\ref{item:C3_one} follows from \eqref{eq:l4.3.4}. We now prove property~\ref{item:C3_two}.

  At intermediate stages of the above procedure, we have
  \begin{align*}
      \prod_{\r\in\cE}\left(1+ K(\r) \cos\la\psi,\r\ra\right) = \sum_{\cI}c_{\cI}\prod_{\r\in\cI}\left(1+ K^{'}(\r) \cos\la\psi,\r\ra\right),
  \end{align*}
  with $\cE\to\cI$. Let $\cE^{'}\in\cF^{'}$ and $\r\in\cE^{'}$. Consider an intermediate ensemble $\cI$. Note that $\r$ must be compatible with $\cI$ (if $\r$ is not compatible with $\cI$, further operations on the factors corresponding to densities in $\cI$ can never produce $\r$). We consider all possible applications of identity \eqref{eq:l4.3.3} to pairs $\{\r_1, \r_2\}\subset \cI$:
\begin{enumerate}
    \item $\pm\r_1\subset \r,\, \r_2\cap\r=\emptyset$, or
    \item $\pm\r_2\subset \r,\, \r_1\cap\r=\emptyset$, or
    \item $\pm\r_1\subset \r, \,\pm\r_2\subset \r$, or
    \item $\r_1\cap\r=\emptyset,\, \r_2\cap\r=\emptyset$.
\end{enumerate}
Clearly, in case $(4)$ the applications of \eqref{eq:l4.3.3} have no effect on $K^{'}(\r)$.

In case $(1)$, the term on the right-hand side of identity \eqref{eq:l4.3.3}, for which $\r$ is still compatible with its ensemble, is one in which $\r_2$ is eliminated. In fact, in this case $\r$ is not compatible with the ensemble $\cI_2$ of the term for which $\r_1$ is eliminated, as $\supp(\r_1)\cap (\cup_{\r\in\cI_1}\supp(\r)) = \emptyset$. $\r$ is also not compatible with the ensembles $\cI_{\pm}$, for which $\r_1 \pm \r_2\in\cI_{\pm}$, since $\r_1\pm\r_2$ cannot be separated any more. Thus, the factor $3K_1$ will be a factor of $K^{'}(\r)$. Note that in every application of \eqref{eq:l4.3.3} of type $(1)$, $\dist(\r,\r^{'}) \le 2^k$ for some $\r^{'}\subset \r_2, \, \r^{'}\in\cE$. Therefore, the total number of applications of identity \eqref{eq:l4.3.3} of type $(1)$ which affect $K^{'}(\r)$ is at most $$|\{\r^{'}\in\cE: \supp(\r)\cap\supp(\r^{'}) =\emptyset, \dist (\r,\r^{'}) \le 2^k\}|.$$
Case $(2)$ is the same as case $(1)$, with $\r_1$ and $\r_2$ interchanged.

In case $(3)$, the term on the right-hand side of identity \eqref{eq:l4.3.3}, for which $\r$ is still compatible with its ensembles, is either the one with a term of the density $\r_1 +\r_2$, or the one with a term of the density $\r_1-\r_2$. In any case, the factor $3K_1 K_2$ will be a factor of $K^{'}(\r)$. The number of such possible applications of \eqref{eq:l4.3.3} is exactly $|\{\r^{'}\in\cE: \r^{'}\subset \r\}|$.

From the discussion of these cases, we now conclude
\begin{align*}
    |K^{'}(\r)| 
    %&\le \left(3^{|\{\r^{'}\in\cE: \r\cap\r^{'} =\emptyset, \dist (\r,\r^{'}) \le 2^k\}|} 3^{|\{\r^{'}\in\cE: \r^{'}\subset \r\}|}\right)\prod_{\r^{'}\in\cE}|K(\r^{'})|^{|\eps(\r, \r^{'})|}\\
    \le 3^{n_{\cE}(\r)}\prod_{\r^{'}\in\cE}|K(\r^{'})|^{|\eps(\r, \r^{'})|},
\end{align*}
where $$n_{\cE}(\r):=|\{\r^{'}\in\cE: \dist (\r,\r^{'}) \le 2^k\}|.$$

To finish the proof, we show that if $\cE$ is a $(k-1)$-ensemble, then $n_{\cE}(\r) \le C_1 A_{k-1}(\r)$, for some positive absolute constant $C_1$.

First, note that the case $k=0$ holds trivially.

For $\ell\ge 0$, let $\overline{A_{\ell}}(\r)$ be the minimal number of squares of side length $2^\ell$ needed to cover $\{x\in\Lda: \dist(x, \supp(\r))\le 2^{\ell+1}\}$. Then $\overline{A_{\ell}}(\r) \le 25 A_{\ell}(\r)$.

Let $k\ge 1$. Note that for any two distinct densities $\r_1, \r_2\in\cE$, it holds that $\dist(\r_1, \r_2)>2^{k-1}$. Therefore, any square of side length $2^{k-1}$ cannot intersect more than one density in $\cE$ (if more than one density intersects such a square, then the distance between two of these densities is less than or equal to $2^{k-1}$, which is a contradiction). Therefore, we obtain
\begin{align*}
   n_{\cE}(\r) \le \overline{A_{k-1}}(\r) \le 25 A_{k-1}(\r). 
\end{align*}
For $k=0$, we have $n_{\cE}(\rho)\le9A_0(\rho)$.
\end{proof}

\begin{remark}
 The final $\cF^{'}$ of Lemma~\ref{lem:4.3} may depend on the order in which we apply \eqref{eq:l4.3.3}.   
\end{remark}

\begin{proof}[Proof of Theorem~\ref{thm:2.2}]
 
Recall 
$$F_N(x)= 1+ \sum_{q=1}^N 2(1-\frac{q}{N})\cos(qx).$$ 
Let $$C(N):=\sum_{q=1}^N\e^{-q^2}.$$
We write
\begin{align*}
    F_N(x)&= \sum_{q=1}^N \left[\frac{\e^{-q^2}}{C(N)} + 2\left(1-\frac{q}{N}\right)\cos(qx)\right]\\
    &=\sum_{q=1}^N \frac{\e^{-q^2}}{C(N)}\left[1 + 2C(N)\e^{q^2}\left(1-\frac{q}{N}\right)\cos(qx)\right].
\end{align*}
Hence 
\begin{align}\label{eq:t2.2.1}
    \prod_{x\in\Lda} F_N(\psi(x)) &= \prod_{x\in\Lda} \sum_{q=1}^N \frac{\e^{-q^2}}{C(N)}\left[1 + 2C(N)\e^{q^2}\left(1-\frac{q}{N}\right)\cos(q\psi(x))\right]\nonumber\\
    &= \sum_{\mathbf{q}} \xi(\mathbf q)\prod_{x\in\Lda} \left[1+ z(q(x))\cos(q(x)\psi(x))\right],
\end{align}
with
\begin{align}\label{eq:t2.2.2}
    \mathbf q=(q(x))_{x\in\Lda}\in\{1,2,\ldots, N\}^{\Lda},\quad \xi(\mathbf q)= \prod_{x\in\Lda}\frac{\e^{-q(x)^2}}{C(N)}, \quad z(q(x))=2C(N)\e^{q(x)^2}\left(1-\frac{q(x)}{N}\right).
\end{align}
Note that $\sum_{\mathbf q} \xi(\mathbf q) =1$.

Therefore, to prove Theorem~\ref{thm:2.2}, it is enough to show that there exist a positive absolute constant $D_2$ and a family $\cF_{\mathbf q}$ of ensembles satisfying properties \ref{item:one_6.4}--\ref{item:three_6.4} of Theorem~\ref{thm:2.2}, such that
\begin{align}\label{eq:t2.2.3}
    \prod_{x\in\Lda} \left[1+ z(q(x))\cos(q(x)\psi(x))\right]=\sum_{\cN\in\cF_{\mathbf q}} c_{\cN}^{'} \prod_{\r\in\cN}\left[1+ K(\r)\cos\la\psi, \r\ra\right],
\end{align}
where $c_{\cN}^{'}>0,\,\cN\in\cF_{\mathbf q}$ and $\sum_{\cN\in\cF_{\mathbf q}}c_{\cN}^{'}=1$.

To prove this, we apply Lemma~\ref{lem:4.3} iteratively to appropriately chosen $k$-ensembles $\cE_k$. Roughly speaking, we shall apply the lemma to ensembles which do not satisfy the assertions of Theorem~\ref{thm:2.2}. We define the iterative process by induction. Let $\cQ_{-1}:=\{\r^x\}_{x\in\Lda}$, where $\r^x:=q(x)\cdot\delta_x$. Also let $\cG_{-1}:=\emptyset$. Note that $\cQ_{-1}$ is a $(-1)$-ensemble. For $k\ge 0$, assume that we have defined $\cQ_{k-1}$ and $\cG_{k-1}$ such that $\cE_{k-1}=\cQ_{k-1}\setminus \cG_{k-1}$ is a $(k-1)$-ensemble. Given $\cQ_{k-1},\,\cG_{k-1},\,k\ge 0$, define $\cQ_k$ to be a $k$-ensemble in the following way. First, if $|\cE_{k-1}|\le 1$, we set $\cQ_k:=\cQ_{k-1},\,\cG_k:=\cG_{k-1}$. Otherwise, we apply Lemma~\ref{lem:4.3} to $\cE_{k-1}$, which is a $(k-1)$-ensemble. Thus we obtain
\begin{align*}
   \prod_{\r\in\cQ_{k-1}}\left[1+ K(\r)\cos\la\psi, \r\ra\right]&=\prod_{\r\in\cG_{k-1}}\left[1+ K(\r)\cos\la\psi, \r\ra\right] \prod_{\r\in\cQ_{k-1}\setminus \cG_{k-1}}\left[1+ K(\r)\cos\la\psi, \r\ra\right]\\
   &=\prod_{\r\in\cG_{k-1}}\left[1+ K(\r)\cos\la\psi, \r\ra\right] \sum_{\cE^{'}\in\cF^{'}} c_{\cE^{'}} \prod_{\r\in\cE^{'}}\left[1+ K^{'}(\r)\cos\la\psi, \r\ra\right]\\
   &=\sum_{\cE^{'}\in\cF^{'}} c_{\cE^{'}} \prod_{\r\in\cE^{'}\cup \cG_{k-1}}\left[1+ K^{'}(\r)\cos\la\psi, \r\ra\right],
\end{align*}
where $K'(\rho)=K(\rho)$ for $\r\in\cG_{k-1}$.
For each $\cE^{'}\in\cF^{'}$ we set $\cQ_k:=\cQ_{k,\cE^{'}}:=\cE^{'}\cup \cG_{k-1}$, and continue the iterations with each $\cQ_k=\cQ_{k,\cE^{'}}$ separately. Given $\cQ_k$, we define $\cG_k$ as follows. We set $\cG_{k,0}=\cG_{k-1}$. Then we order the densities in $\cQ_k\setminus \cG_{k-1}$ in ascending order of $d_\Lda(\r)$. We let $(\r_n)_{n=1}^{m_k}$ be the sequence, where $m_k=|\cQ_k\setminus \cG_{k-1}|$. For $n\ge 1$, we set $\cG_{k, n}:=\cG_{k,n-1}\cup\{\r_n\}$ if
\begin{align}\label{eq:t2.2.5}
    \dist(\r^{'}, \r_n) \ge M d_\Lda(\r_n)^\a, \qquad \r^{'}\in\cQ_k\setminus(\cG_{k,n-1}\cup\{\r_n\}).
\end{align}
Otherwise, we set $\cG_{k,n}=\cG_{k,n-1}$. Finally, we let $\cG_k:=\cG_{k,m_k}$. With this definition, $\cG_{k-1}\subset \cG_k$ and $\cE_k:=\cQ_k\setminus \cG_k$ is a $k$-ensemble (since $\cE_k\subset \cE^{'}$). Moreover, if $\r, \r^{'}$ are distinct densities such that 
\begin{align*}
   \dist(\r, \r^{'}) < M \min\{d_\Lda(\r), d_\Lda(\r^{'}) \}^\a, 
\end{align*}
then $\r, \r^{'}\in\cE_k$. Thus, $|\cE_k|\le 1$ implies that $\cQ_k$ satisfies property~\ref{item:one_6.4} of Theorem~\ref{thm:2.2}.

There exists $k^{*}=k^{*}(L)$ such that any $k$-ensemble $\cE$ with $k\ge k^{*}$ satisfies $|\cE|\le 1$. Therefore, we can terminate the above process after $k^{*}$ iterations to obtain a family $\cF$, positive coefficients $(C_\cN)_{\cN\in\cF}$ summing to $1$ and real coefficients $(K(\r))_{\r\in\cN\in\cF}$ satisfying property~\ref{item:one_6.4} of Theorem~\ref{thm:2.2} such that for every $\psi:\Lda \to \R$
\begin{align*}
    \prod_{x\in\Lda} F_N(\psi(x)) = \sum_{\cN\in\cF} C_{\cN} \prod_{\r\in\cN}\left[1+ K(\r)\cos(\la \psi, \r \ra)\right].
\end{align*}
We now verify the other properties of Theorem~\ref{thm:2.2}.

To verify property $(2)$, let $\cN\in\cF$ and $\r_1\subseteq \r\in\cN$, $\r_1\neq \r$, such that $$\dist(\r_1, \r-\r_1) \ge 2Md(\r_1)^\a=:R.$$
For each integer $1\le k \le k^{*}$, let $\cQ_k$ be the (unique) ensemble generated at the $k$th stage of the
iterative process by which we constructed $\cF$ such that $\cN$ was constructed by successive applications
of Lemma~\ref{lem:4.3} to $\cQ_k$. Note that $\cQ_k \to \cN$. Let $k$ be the smallest integer such that there exists $\r^{*}\in\cQ_k$ satisfying $\supp(\r_1)\cap \supp(\r^{*}) \neq \emptyset$ and $\supp(\r-\r_1)\cap \supp(\r^{*})\neq \emptyset$. Then there exists $\r_\mu, \r_\nu \in \cE_{k-1}$ such that $\r_\mu, \r_\nu \subset \r^{*}$, $\r_\mu\subset \r_1$ and $\r_\nu \subset \r-\r_1$. We have $\dist(\r_\mu, \r_\nu) \ge R$ and by property~\ref{item:C3_one} of Lemma~\ref{lem:4.3}, $R\le 2^k$. This implies $d(\r_1) < 2^{k-1}$ which gives $\r_1=\r_\mu$ since $\r_1$ is compatible with $\cE_{k-1}$ (which contains no densities within distance $2^{k-1}$ of $\r_\mu$). Moreover, $\cE_{k-1}$ is a $(k-1)$-ensemble and by assumption there exists $\r_\nu\in\cE_{k-1}$ such that $$2^{k-1} < \dist(\r_1, \r_\nu) < M d_\Lda(\r_1)^\a.$$ Thus we have $2M d(\r_1)^\a =R\le 2^k < 2Md_\Lda(\r_1)^\a$.
It follows that $d_\Lda(\r_1)= \dist (\r_1, \partial_2\Lda)$ and $\mom(\r_1)\neq 0$. Also, $$2M \dist(\r_1, \partial_2\Lda)^\a > 2^k \ge \dist(\r_\mu, \r_\nu)\ge \dist(\r_1, \r-\r_1).$$

Next, we verify property $(3)$ of Theorem~\ref{thm:2.2}. Let $\cN\in\cF$ and $\r\in\cN$. For any $1\le k\le k^{*}$, let $\cQ_k$ be as above. We claim that $\r\in\cG_{n(\r)}$, where $n(\r)=\ceil{\log_2(M\cdot d_\Lda(\r)^\a)}$ as in \eqref{eq:def_nrho}. Indeed, since $2^{n(\r)} \ge M d_\Lda(\r)^\a$ and $\cE_{n(\r)}$ is an $n(\r)$-ensemble, we must have $\r\in\cG_{n(\r)}$ (otherwise, there must exist $\r_1, \r_2\in \cE_{n(\r)}$ such that $\r_1, \r_2 \subset \r$ and hence $d(\r) > 2^{n(\r)}$).

Let $m\le n(\r)$ be the smallest index such that $\r\in\cG_{m}$. Then by part $(2)$ of Lemma~\ref{lem:4.3}
\begin{align*}
    |K(\r)| \le \e^{C_1 A_{m-1}(\r)} \prod_{\r_{\g}}|K(\r_\g)|,
\end{align*}
where $\r=\sum \eps(\r,\r_\g)\r_\g$ with $\eps(\r,\r_\g)\in\{\pm 1\}$ and all densities $\r_\g$ belong to some $(m-1)$-ensemble. 

If $m\ge 1$, we apply part $(2)$ of Lemma~\ref{lem:4.3} again, which yields
\begin{align*}
   |K(\r)| \le \exp\left[C_1\left( A_{m-1}(\r) + \sum_{\r_\g\subset\r} A_{m-2}(\r_\g)\right)\right] \prod_{\r_{\delta}}|K(\r_\delta)|, 
\end{align*}
for densities $\r_\delta$ in some $(m-2)$-ensemble.

Note that 
\begin{align*}
   \sum_{\r_\g\subset\r} A_{m-2}(\r_\g) \le A_{m-2}(\r),
\end{align*}
since every square of side length $2^{m-2}$ can intersect at most one density $\r_\g$ and since a cover of $\r$ with squares of side length $2^{m-2}$ is also a cover of $\r_\g$ for all $\g$.

Iterating this argument, we obtain
\begin{align*}
   |K(\r)| \le \exp\left[C_1\sum_{m=-1}^{n(\r)-1} A_{m}(\r) \right] \prod_{\r_{\delta}}|K(\r_\delta)|, 
\end{align*}
for densities $\r_\delta$ in some $(-1)$-ensemble.
But in that case, 
\begin{align*}
\prod_{\r_{\delta}}|K(\r_\delta)|\le \prod_{x\in\supp(\r)} |z(|\r(x)|)|.   
\end{align*}
Thus we have
\begin{align*}
   |K(\r)| &\le \exp\left[C_1\sum_{m=-1}^{n(\r)-1} A_{m}(\r) \right] \prod_{x\in\supp(\r)} \left|2C(N)\e^{\r(x)^2}\left(1-\frac{|\r(x)|}{N}\right)\right|\\
   &\le \e^{D_2 A(\r)}\prod_{x\in\supp(\r)} \e^{\r(x)^2}\left|\left(1-\frac{|\r(x)|}{N}\right)\right|.
\end{align*}

\end{proof}

%\clearpage
\bibliographystyle{abbrvnat}
\bibliography{bib_ivmm}
\end{document}